\documentclass[12pt]{article}
\usepackage{amsmath,amsfonts,amssymb,theorem}
\usepackage{type1cm}
\newcommand{\simby}[1]{\buildrel{#1}\over\sim}    % for eqce relns
\newcommand{\simT}{\simby{\T}}                    % TT-equivalence

\newcommand{\comment}[1]{}

\newcommand{\lank}{\ \cdot\ }
\newcommand{\dth}{{\textstyle\frac1d}}
\newcommand{\eeth}{{\textstyle\frac1e}}
\newcommand{\ka}{\kappa}

\newcommand{\Aff}{\mathbb A}

\newcommand{\C}{\mathbb C}
\newcommand{\Cstar}{\C^\times}
\newcommand{\Gm}{\mathbb G_m}

\newcommand{\M}{\mathbb M}
\newcommand{\bMbar}{\overline{\mathbb M}}

\newcommand{\PP}{\mathbb P}
\newcommand{\Q}{\mathbb Q}
\newcommand{\R}{\mathbb R}
\newcommand{\T}{\mathbb T}
\newcommand{\Z}{\mathbb Z}

\newcommand{\sI}{\mathcal I}

\newcommand{\Oh}{\mathcal O}

\newcommand{\ep}{\varepsilon}
\newcommand{\fie}{\varphi}
\newcommand{\al}{\alpha}
\newcommand{\be}{\beta}
\newcommand{\ga}{\gamma}
\newcommand{\de}{\delta}

\newcommand{\si}{\sigma}
\newcommand{\la}{\lambda}
\newcommand{\La}{\Lambda}
\newcommand{\om}{\omega}
\newcommand{\Om}{\Omega}
\newcommand{\Ga}{\Gamma}

\newcommand{\1}{^{-1}}
\newcommand{\iso}{\cong}
\newcommand{\bij}{\leftrightarrow}

\newcommand{\broken}{\dasharrow}

\newcommand{\Span}[1]{\left<#1\right>}
\newcommand{\dbk}[1]{[\kern-.15em[#1]\kern-.15em]} % to give [[ ]]
\newcommand{\rdup}[1]{\lceil#1\rceil}

\DeclareMathOperator{\divi}{div}

\DeclareMathOperator{\hcf}{hcf}
\DeclareMathOperator{\Hom}{Hom}
\DeclareMathOperator{\id}{id}
\DeclareMathOperator{\Pf}{Pf}

\DeclareMathOperator{\PSL}{PSL}
\DeclareMathOperator{\SL}{SL}
\DeclareMathOperator{\Spec}{Spec}
\DeclareMathOperator{\sHom}{\mathcal H \mathit{om}}

\newtheorem{theorem}{Theorem}[section]
\newtheorem{claim}[theorem]{Claim}
\newtheorem{cla}{Claim}
 
\newtheorem{lem}[theorem]{Lemma}
\newtheorem{prop}[theorem]{Proposition}

\newtheorem{cor}[theorem]{Corollary}
\newtheorem{indass}[theorem]{Inductive assumption}

{
\theorembodyfont{\rmfamily}

\newtheorem{exa}[theorem]{Example}

\newtheorem{remark}[theorem]{Remark}

}

\renewcommand{\mod}{\ \mathrm{mod}\ }
\newcommand{\dd}{\mathrm d\,}

\newcommand{\QED}{Q.E.D.}
\newcommand{\qed}{\quad\ensuremath{\square}}
\newenvironment{pf}{\paragraph{Proof}}{{\hfill \QED}\par\medskip}

\newcommand{\step}[1]{\paragraph{\sc Step #1}}
\newcommand{\case}[1]{\paragraph{\sc Case #1}}

\numberwithin{equation}{section}
\numberwithin{figure}{section}
\numberwithin{table}{section}

\begin{document}
\setcounter{tocdepth}{2}

\title{Diptych varieties. I}
\author{Gavin Brown and Miles Reid
\thanks{Partially funded by Korean Government WCU Grant
R33-2008-000-10101-0}}
\date{}
\maketitle

\begin{abstract} We present a new class of affine Gorenstein 6-folds
obtained by smoothing the 1-dimensional singular locus of a reducible
affine toric surface; their existence is established using explicit
methods in toric geometry and serial use of Kustin--Miller Gorenstein
unprojection. These varieties have applications as key varieties in
constructing other varieties, including local models of Mori flips of
Type~A. \end{abstract}

%%%%%%%%%%%%%%%%%%%%%%%%%%%%%%%%%%%%%%%%%%%%%%%%%%%%%%%%%

We introduce a large class of remarkable 6-folds called {\em diptych
varieties}. Each is an affine 6-fold $V_{ABLM}$ constructed starting
from two toric 4-fold panels $V_{AB}\cup V_{LM}$ hinged along a
reducible toric surface $T=V_{AB}\cap V_{LM}$ (compare the Wilton
diptych \cite{W}). The construction depends on discrete toric data
called a {\em diptych of long rectangles}, that describe the monomial
cone of the two toric panels $V_{AB}$ and $V_{LM}$. It is equivariant
under a big torus $\T=(\Gm)^4=(\Cstar)^4$. Apart from easy initial
cases, diptych varieties are indexed by 3 natural numbers $d,e,k$, or by
a 2-step recurrent continued fraction $[d,e,d,\dots,(\hbox{$d$ or
$e$})]$ to $k$ terms (Classification Theorem~\ref{th!d-e}). Once this
combinatorial data is set up, Main Theorem~\ref{th!main} guarantees the
existence of the diptych variety. The worked example \ref{ex!intr}
illustrates almost all the main features of our construction. This paper
is backed up by a website
\begin{quote}
\verb!http://www-staff.lboro.ac.uk/~magdb/aflip.html!
\end{quote}
that contains current drafts of Parts~II--IV, together with computer
algebra calculations, links to other papers and further auxiliary
material.

Diptych varieties $V_{ABLM}$ are designed for use as ambient spaces or
{\em key varieties} in constructing other spaces, much as toric
varieties. As discussed briefly in the final Section~\ref{s!fin}, our
main motivation is their relation with the ``continued division''
algorithm \cite{Mo}, that Mori used to prove the existence of flips of
Type~A. Our work also overlaps with the more recent Gross--Hacking--Keel
deformations of cycles of planes \cite{GHK} in some cases where these
lead to algebraic varieties.

\tableofcontents

%%%%%%%%%%%%%%%%%%%%%%%%%%%%%%%%%%%%%%%%%%%%%%%%%%%%%%%%%

\section{Introduction}

This section gives rough statements of our main results and an outline
plan of the paper. The extended example of \ref{ex!intr} illustrates all
the main ideas. We write $\Aff^n=\C^n$ for affine space, $\Gm=\Cstar$
for the multiplicative group and $\T=(\C^\times)^4$ for the
4-dimensional torus. Our main interest is in varieties over $\C$,
although in the final analysis, our diptych varieties are defined as
schemes over $\Z$.

\subsection{Main results and overview of the paper}
\label{intro!I}

A {\em tent} is a reducible affine surface $T=S_0\cup S_1\cup S_2\cup
S_3$ as in Figure~\ref{f!T}. Its four irreducible components are
$S_0,S_2\iso\Aff^2$ and $S_1,S_3$ cyclic quotient singularities of type
$\frac1r(\al,1)$ and $\frac1s(\be,1)$, where $r, \al$ are coprime
natural numbers, and similarly for $s,\be$. We glue the four toric
surfaces transversally along their toric strata, giving $T$ four
1-dimensional singular axes of transverse ordinary double points; the
two axes on $S_2$ are the {\em top} axes of $T$, and the two on $S_0$
its {\em bottom} axes.

Section~\ref{s!trc} recalls basic facts on toric geometry and studies
certain deformations of tents. Our first result is Theorem~\ref{l!be}:
an extension $T\subset V_{AB}$ of a tent $T$ to an affine toric 4-fold
$V_{AB}$ that smooths the top axes is determined by a matrix
$\left(\begin{smallmatrix} r&a\\b&s \end{smallmatrix}\right)
\in\SL(2,\Z)$ with $a,b\ge0$ and $a\equiv\al\mod r$, $b\equiv\be\mod s$.
Corollary~\ref{cor!VAB} gives an alternative statement in terms of {\em
continued fraction expansions of\/~$0$}, obtained by concatenating with
a 1 the expansions of {\em complementary} fractions $\frac r\be$ and
$\frac r{r-\al}$. This is routine material in toric geometry, but the
basic results and detailed notation for the monomial cone $\si_{AB}$
introduced here are in use throughout the paper.

Section~\ref{s!Cl} treats our first substantial result, Classification
Theorem~\ref{th!d-e}, classifying diptychs of toric extensions $T\subset
V_{AB}$ and $T\subset V_{LM}$ that smooth respectively the top and
bottom axes of $T$. By Lemma~\ref{l!2nd}, the numerical conditions on
$T$ for the second smoothing to exist is a second matrix
$\left(\begin{smallmatrix} r&g\\h&s \end{smallmatrix}\right)
\in\SL(2,\Z)$ with $ag\equiv1\mod r$ and $bh\equiv1\mod s$.
Theorem~\ref{th!d-e} classifies all solutions to this problem: with
simple initial exceptions, each corresponds to a \hbox{2-step} recurrent
continued fraction $[d,e,d,\dots,\hbox{($d$ or $e$)}]$.
Theorem~\ref{th!d-e} is proved by a simple descent argument.

At this point we introduce a case division (we discuss the necessity for
this briefly in Section~\ref{s!fin}). The {\em main case} is $d,e\ge2$
and $de>4$; we concentrate our efforts primarily on this case in the
rest of the current paper. The other cases involve some new features,
and their proofs require minor modifications; they are as follows:

\begin{itemize}
\item $de\le3$. This involves only a small number of quite small cases,
and we deal with them in an appendix to \cite{BR2}.

\item The cases $d=e=2$ and $d=1$, $e=4$ are treated in \cite{BR2}.
There are two infinite series of varieties with a convincing standard
quasihomogeneous structure.

\item $d$ or $e=1$ and  $de>4$. This case requires a proof that is
basically on the same scale as the main case; we relegate the details to
\cite{BR3} to avoid excessive repetition, bulky notation, and many case
divisions.
\end{itemize}

Diptychs serve as the input to our Main Theorem, the existence of
diptych varieties:

\begin{theorem}\label{th!main}
A diptych of $4$-fold toric panels $T\subset V_{AB}$ and $T\subset
V_{LM}$ that smooth respectively the top and bottom axes of $T$ extends
to a $6$-fold $V_{ABLM}$:
\begin{equation}
\renewcommand{\arraystretch}{1.5}
\renewcommand{\arraycolsep}{.6em}
\begin{array}{ccc}
T & \subset & V_{AB} \\
\bigcap && \bigcap \\
\kern2mm V_{LM} & \subset & V_{ABLM}\kern-4mm
\end{array}
\end{equation}
The {\em diptych variety} $V_{ABLM}$ is an affine variety with an action
of the torus $\T=(\Gm)^4$. It has a regular sequence $A,B,L,M$
consisting of eigenfunctions of the $\T$-action such that $V_{AB}$ and
$V_{LM}$ are the sections given by $L=M=0$ and $A=B=0$, with $T$ their
intersection $A=B=L=M=0$.
\end{theorem}

It follows that $V_{ABLM}$ is a Gorenstein affine $6$-fold and is a flat
$4$-parameter deformation of the tent $T$. The $\T$-action restricts to
the big torus of both $4$-fold panels $V_{AB}$ and $V_{LM}$; the
original tent $T$ is a union of toric strata in each, with the
$\T$-action inducing the natural $(\Gm)^2$ action on each of its four
toric components.

Section~\ref{s!pp} lays the groundwork for the proof in
Section~\ref{s!pf}. The main idea is to exploit the relation between
the monomial lattices and the monomial cones of the two different toric
varieties $V_{AB}$ and $V_{LM}$ to deduce important consequences for
monomials in the coordinate ring of the diptych variety $V_{ABLM}$. Our
proof of Main Theorem~\ref{th!main} in Section~\ref{s!pf} makes
essential use of convexity properties of these monomials (illustrated in
the Pretty Polytope of Figure~\ref{f!pp}) and congruence properties (the
Padded Cell of Figure~\ref{f!pad}).

Section~\ref{s!pf} proves Theorem~\ref{th!main} in the main case by
{\em serial unprojection}. We start from two equations defining a
codim\-ension~2 complete intersection
$V_0\subset\Aff^8_{\Span{x_0,x_1,y_0,y_1,A,B,L,M}}$, and adjoin the
remaining variables one at a time by unprojection $V_{\nu+1}\to V_\nu$.
Section~\ref{ss!pseq} determines the {\em unprojection order} in which we
must adjoin the variables $x_2,\dots,y_l$. It is inverse to the order of
elimination (or projection) of variables from the toric panel $V_{AB}$,
corresponding to the concatenated continued fraction
$[a_2,\dots,a_k,b_l,\dots,b_1]=0$. Serial use of the Kustin--Miller
unprojection theorem of \cite{PR} provides most of what we need.

Extended Example~\ref{ex!intr} is the case corresponding to the
recurrent continued fraction $[2,4,2]$ or the expansion of zero
$[4,2,1,3,2,2]=0$. We use a beautiful trick with Pfaffians to compute
the sequence of unprojection variables $[x_2,y_2,y_3,x_3,y_4]$ as
rational functions with specified poles, the geometric interpretation of
Kustin--Miller unprojection. This example illustrates all but one of the
main points, and exemplifies our strategy of handling a diptych variety
$V_{ABLM}$ as an explicit object, but without necessarily writing down
all the relations for its coordinate ring, much as for a toric variety.

The extended example glosses over one logical point that is the key
issue for most of Sections~\ref{s!pp}--\ref{s!pf}. Each step
$V_{\nu+1}\to V_\nu$ of the induction must set up a new unprojection
divisor $D_\nu\subset V_\nu$. The divisor $D_\nu$ itself is the product
of a monomial curve $A^\al B^\be=0$ with an affine space
$\Aff^4_{\Span{x_i,y_j,L,M}}$, but we still have to prove it is a
subscheme of $V_\nu$.

\subsection{Extended example} \label{ex!intr}

\subsubsection{Background and notation}
For $r>0$ and $a$ coprime to $r$, we write $\frac1r(1,a)$ for
the action of $\Z/r$ on $\Aff^2$ given by $(u,v)\mapsto(\ep u,\ep^av)$
where $\ep=\exp\frac{2\pi i}r\in\C$ is a chosen primitive $r$th root
of~1. We use the same notation for the cyclic quotient singularity
$\Aff^2/(\Z/r)=\Spec\C[u,v]^{\Z/r}$. We focus here on concrete cases,
starting with $\frac17(1,2)$; the ring of invariants $\C[u,v]^{\Z/7}$ is
generated by the monomials
\begin{equation}
y_0=u^7, \quad y_1=u^5v, \quad y_2=u^3v^2, \quad y_3=uv^3, \quad y_4=v^7,
\end{equation}
with relations between them determined by the {\em tag equations}
\begin{equation}
y_0y_2=y_1^2,\quad y_1y_3=y_2^2,\quad y_2y_4=y_3^3.
\label{eq!tag72}
\end{equation}
These are of the general form $v_{i-1}v_{i+1}=v_i^{a_i}$ for any 3
consecutive monomials $v_{i-1},v_i,v_{i+1}$ on the Newton boundary. The
exponents or {\em tags} $a_i$ are the entries in the Jung--Hirzebruch
continued fraction expansion of $\frac r{r-a}$; here
$\frac7{7-2}=2-\frac{1}{2-\frac13}=[2,2,3]$. The quotient $\Aff^2\to
S\subset\Aff^5_{\Span{y_{0\dots4}}}$ is thus the morphism
$(u,v)\mapsto(y_{0\dots4})$, and the image $S$ is uniquely determined by
\eqref{eq!tag72}: the complete intersection \eqref{eq!tag72} consists of
$S$ plus the $(y_0,y_4)$-plane with a ``fat'' nonreduced structure. To see
actual generators of the ideal $I_S$ we also need the ``long equations''
$y_0y_3=y_1y_2$, $y_1y_4=y_2y_3^2$ and $y_0y_4=y_1y_3^2$, that derive
from \eqref{eq!tag72} using easy syzygy manipulations. In what follows,
we write $S=S_3$ for the quotient $\frac17(1,2)$.

In the same way, the quotient singularity $\frac17(1,3)$ is
\begin{equation}
S_1\subset\Aff^4_{\Span{x_0,x_1,x_2,x_3}} \quad\hbox{given by}\quad
x_0x_2=x_1^2, \quad x_1x_3=x_2^4
\label{eq!tag73}
\end{equation}
with $[2,4]=2-\frac14=\frac74$.

\subsubsection{The tent $T$}
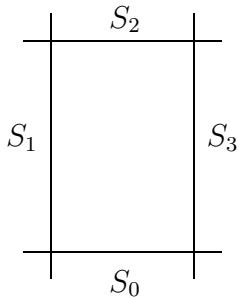
\begin{figure}[ht]
\begin{picture}(150,110)(-112,0)
\put(48,5){\line(0,1){100}}
\put(102,5){\line(0,1){100}}
\put(37.5,15){\line(1,0){75}}
\put(37.5,95){\line(1,0){75}}
\put(31,55){$S_1$}
\put(107,55){$S_3$}
\put(70,0){$S_0$}
\put(70,100){$S_2$}
\end{picture}
\caption{The tent $T=S_0\cup S_1\cup S_2\cup S_3\subset
\Aff^{k+l+2}_{\Span{x_{0\dots k},y_{0\dots l}}}$ is obtained by glueing $S_0\cup S_1$
transversally along the $x_0$-axis, $S_1\cup S_2$ along the $x_k$-axis,
$S_2\cup S_3$ along the $y_l$-axis, and $S_3\cup S_0$ along the
$y_0$-axis}
\label{f!T}
\end{figure}

The starting point for our example is the reducible affine surface or
{\em tent} of Figure~\ref{f!T} (with $k=3$, $l=4$ and
$k+l+2=9$ in our case). It consists of a cycle of 4 components, with
vertical sides the surface quotient singularities
$S_1\subset\Aff^4_{\Span{x_{0\dots k}}}$ and
$S_3\subset\Aff^5_{\Span{y_{0\dots l}}}$ of types $\frac17(1,3)$ and
$\frac17(1,2)$ as just described, and top and bottom the coordinate
planes $S_2=\Aff^2_{\Span{x_k,y_l}}$ and $S_0=\Aff^2_{\Span{x_0,y_0}}$.
In equations, $T\subset\Aff^9$ is the reducible variety defined by
\begin{equation}
I_{S_1}, I_{S_3} \quad \hbox{and} \quad
x_iy_j=0 \quad\hbox{for all $i,j$ with $(i,j)\ne(0,0),(k,l)$.}
\label{eq!T}
\end{equation}

\subsubsection{First toric extension $T\subset V_{AB}$}
We now seek to embed $T$ into a toric variety $V$ (irreducible and
normal) so that $T$ is both a regular section of $V$ and a union of
toric strata.

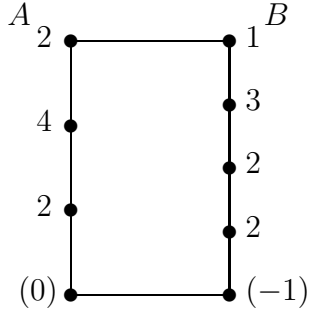
\begin{figure}[ht]
\begin{picture}(150,110)(-112,0)
\put(45,7){\line(0,1){96}}
\put(105,7){\line(0,1){96}}
\put(45,7){\line(1,0){60}}
\put(45,103){\line(1,0){60}}
\put(45,7){\circle*{5}}
\put(45,39){\circle*{5}}
\put(45,71){\circle*{5}}
\put(45,103){\circle*{5}}
\put(25,5){(0)}
\put(32,37){2}
\put(32,69){4}
\put(32,101){2}
\put(21,109){$A$}
\put(105,7){\circle*{5}}
\put(105,31){\circle*{5}}
\put(105,55){\circle*{5}}
\put(105,79){\circle*{5}}
\put(105,103){\circle*{5}}
\put(111,5){$(-1)$}
\put(111,29){2}
\put(111,53){2}
\put(111,77){3}
\put(111,101){1}
\put(118,109){$B$}
\end{picture}
\caption{The long rectangle for $V_{AB}$}
\label{f!lrAB}
\end{figure}

One solution is the affine toric 4-fold $V_{AB}$ with monomial cone
schematically represented in Figure~\ref{f!lrAB}, our first {\em long
rectangle}.
It is a schematic representation of a cone $\si(V_{AB})$, the Newton
polygon of $V_{AB}$ in the monomial lattice $\M=\Z^4$.
We read $\si(V_{AB})$ and the toric
variety $V_{AB}$ automatically from the figure as follows: the dots
around the boundary (clockwise from bottom left) are the generators
$x_{0\dots k}$, $y_{l\dots0}$; the two remaining generators $A,B$ are
shown as {\em annotations} at the top corners. We also draw them in
their correct geo\-metric position in the 4-dimensional lattice $\M$ in
Figure~\ref{f!siAB}, but this long rectangle shorthand is usually more
convenient. The relations \eqref{eq!tag72} and \eqref{eq!tag73} continue
to hold, as represented by the tags down the long sides. These constrain
$x_{0\dots k}$ to a plane face of $\si(V_{AB})$, and in that plane they
generate the Newton boundary of $\frac17(1,4)$; ditto $y_{0\dots l}$.
The new ingredients are the tags and annotations $^A2,1^B$ at the top
corners, that say how we intend to deform the reducible equations
$x_2y_4=0$ and $x_3y_3=0$ for $T$ appearing in \eqref{eq!T} to usual
binomial equations of toric geometry:
\begin{equation}
x_2y_4=x_3^2A, \quad x_3y_3=y_4B.
\label{eq!AB}
\end{equation}
We view $A$ and $B$ as deformation parameters, and interpret
\eqref{eq!AB} as smoothing the reducible double locus along the $x_3$-
and $y_4$-axes, the top corners $S_1\cap S_2$ and $S_2\cap S_3$ of
Figure~\ref{f!T}.

On the other hand, equations \eqref{eq!AB} and the original tag
equations (\ref{eq!tag72}--\ref{eq!tag73}) now completely determine the
cone $\si(V_{AB})$ in a monomial lattice $\M=\Z^4$. Indeed, $x_3,y_4,A,B$
is a $\Z$-basis of $\M$, and the remaining generators $x_2,\dots,x_0$,
$y_3,\dots,y_0$ are Laurent monomials in this basis, obtained by {\em
continued division} from \eqref{eq!AB} together with
(\ref{eq!tag72}--\ref{eq!tag73}):
\begin{equation}
\renewcommand{\arraystretch}{1.3}
\begin{array}{l}
x_2=x_3^2(Ay_4\1),\\
x_1=x_3^7(Ay_4\1)^4,\\
x_0=x_3^{12}(Ay_4\1)^7,
\end{array} \quad
\renewcommand{\arraystretch}{1.1}
\begin{array}{l}
y_3=y_4(Bx_3\1),\\
y_2=y_4^2(Bx_3\1)^3,\\
y_1=y_4^3(Bx_3\1)^5,\\
y_0=y_4^4(Bx_3\1)^7.
\end{array}
\label{eq!Laur}
\end{equation}

A rational polyhedral cone $\si$ in the monomial lattice $\M$ defines a
irreducible, normal toric variety $V_{\M,\si}=\Spec\C[\M\cap\si]$. We
claim more: our monomials $x_{0\dots3},y_{0\dots4},A,B$ in $\M$ generate
$\M\cap\si_{AB}$, and the resulting toric variety
$V_{AB}=\Spec\C[\M\cap\si_{AB}]$ is a flat deformation of $T$. When we
say deformation, we mean the total space of the deformation; in fact
$A,B$ define a flat morphism $V_{AB}\to\Aff^2_{\Span{A,B}}$ with fibre
$T:(A=B=0)$ over $0$, although this morphism does not figure prominently
in our considerations.

The relations satisfied by our monomials come implicitly from their
inclusion in $\M$. We are usually not interested in writing them all out,
but we want to find enough equations to justify our claim. By
substituting from \eqref{eq!Laur}, we find the relation
\begin{equation}
x_1y_0=A^4B^7
\label{eq!tag0}
\end{equation}
that deforms the original equation $x_1y_0=0$ in $T$; this is the {\em
corner tag} $(0)$ of Figure~\ref{f!lrAB}, indicating a tag equation at
$x_0$, with tag~0 derived from the other tags (the annotation $A^4B^7$
is left implicit). We view it as a partial smoothing of the reducible
double locus of $T$ along the $x_0$-axis -- \eqref{eq!tag0} of course
defines a normal hypersurface in $\Aff^4_{\Span{x_1,y_0,A,B}}$

Now, how does the relation $x_0y_1=0$ deform? {From} \eqref{eq!Laur}
we write out $x_0y_1=x_3^7y_4^{-4}A^7B^5$, hence
\begin{equation}
x_0y_1=y_0\1A^7B^{12} \quad\hbox{or} \quad x_0y_1=x_1A^3B^5.
\label{eq!tag-1}
\end{equation}
The first equality is a tag equation for $y_0$, with negative tag $-1$;
this is the $(-1)$ at the bottom right of Figure~\ref{f!lrAB}. Along
the $y_0$-axis of $T$, where $y_0\ne0$, \eqref{eq!tag-1} ensures that
the $A,B$ deformation is also a partial smoothing of the singularity,
making it irreducible and normal. However, \eqref{eq!tag-1} with its
negative tag is anomalous in that it is not a polynomial equation, so we
are not really allowed to use it as a generator of the ideal of the
affine variety $V_{AB}$. We thus replace it by the second expression,
which in view of \eqref{eq!tag0} is equivalent to it where $y_0\ne0$.
The relation $x_0y_1=x_1A^3B^5$ is also anomalous as a tag equation for
$y_0$, since it involves the ``opposite'' generator $x_1$ in place of
$y_0$. Now the equations of $V_{AB}$ include
(\ref{eq!tag0}--\ref{eq!tag-1}); these define an irreducible normal
complete intersection in $\Aff^6_{\Span{x_0,x_1,y_0,y_1,A,B}}$.

Since $V_{AB}$ is a toric 4-fold, it is Cohen--Macaulay; we see in
Lemma~\ref{l!tnt} that it is also Gorenstein. (Or one checks directly from
the description above that the semigroup ideal of interior monomials of
$\si(V_{AB})$ is generated by $AB$; compare \ref{ss!bet} and
Figure~\ref{f!siAB}.) One checks that the locus $(A=B=0)$ inside $V_{AB}$
equals $T$ at the general point of each component, and in particular
each component is 2-dimensional. Therefore $A,B$ is a regular sequence
and $T\subset V_{AB}$ is a flat deformation.

\subsubsection{Conclusion} In this example we found the $A,B$
deformation $T\subset V_{AB}$ in a more-or-less inevitable way starting
from the new tag equations $x_2y_4=x_3^2A$ and $x_3y_3=y_4B$, that
naturally smooth the double locus of $T$ along the $x_3$- and
$y_4$-axes. After a monomial calculation that is birationally forced,
our rectangle closed up neatly to give the tag equations
(\ref{eq!tag0}--\ref{eq!tag-1}), so that this deformation also leads to
partial smoothings of the $x_0$ and $y_0$-axes, giving an irreducible
and normal variety $V_{AB}$ such that $A=B=0$ contains $S_0$ as a
reduced component. Corollary~\ref{cor!VAB} explains that this miracle
works precisely because the concatenation $[4,2,1,3,2,2]$ is a {\em
continued fraction expansion of\ $0$}. These numbers are the tags at
$x_2,x_3$, $y_4,\dots,y_1$; the asymmetry ($x_1$ omitted but $y_1$
included) is significant, and relates to the anomalous tag equations
\eqref{eq!tag-1}.

\subsubsection{Second toric extension $T\subset V_{LM}$} 
\begin{figure}[ht]
\begin{picture}(150,111)(-112,0)
\put(45,7){\line(0,1){96}}
\put(105,7){\line(0,1){96}}
\put(45,7){\line(1,0){60}}
\put(45,103){\line(1,0){60}}
\put(45,7){\circle*{5}}
\put(45,39){\circle*{5}}
\put(45,71){\circle*{5}}
\put(45,103){\circle*{5}}
\put(21,-3){$L$}
\put(32,4){4}
\put(32,37){2}
\put(32,69){4}
\put(24,101){(0)}
\put(105,7){\circle*{5}}
\put(105,31){\circle*{5}}
\put(105,55){\circle*{5}}
\put(105,79){\circle*{5}}
\put(105,103){\circle*{5}}
\put(111,4){1}
\put(118,-3){$M$}
\put(111,29){2}
\put(111,53){2}
\put(111,77){3}
\put(111,101){$(-3)$}
\end{picture}
\caption{The long rectangle for $V_{LM}$}
\label{f!lrLM}
\end{figure}
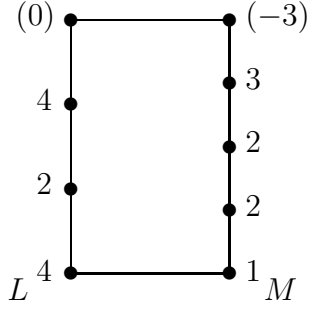

As hinted
above, $T$ has more than one deformation to a toric 4-fold. We now write
down the second long rectangle Figure~\ref{f!lrLM} and the resulting
deformation $T\subset V_{LM}$. The calculations are just as for
$V_{AB}$, except that we start from the bottom and work up. Hindsight
based on Corollary~\ref{cor!VAB} and $[3,2,2,1,4,2]=0$ tells us that
this will work. The new tag equations that smooth out the $x_0$- and
$y_0$-axes of $T$ are represented by the ${}_L4,1_M$ at the bottom:
\begin{equation}
x_1y_0=x_0^4L, \quad x_0y_1=y_0M.
\label{eq!LM}
\end{equation}
This time $x_0,y_0,L,M$ base the monomial lattice and \eqref{eq!LM}
together with (\ref{eq!tag72}--\ref{eq!tag73}) give the remaining
variables as Laurent monomials:
\begin{equation}
\renewcommand{\arraystretch}{1.3}
\begin{array}{l}
x_1=x_0^4(Ly_0\1), \\
x_2=x_0^7(Ly_0\1)^2, \\
x_3=x_0^{24}(Ly_0\1)^7,
\end{array} \quad
\renewcommand{\arraystretch}{1.1}
\begin{array}{l}
y_1=y_0(Mx_0\1), \\
y_2=y_0(Mx_0\1)^2, \\
y_3=y_0(Mx_0\1)^3, \\
y_4=y_0^2(Mx_0\1)^7.
\end{array}
\end{equation}
As before, we deduce the tag equations for $x_3$ and $y_4$:
\begin{equation}
x_2y_4=L^2M^7, \quad x_3y_3=y_4^{-3}L^7M^{27}=x_2^3LM^3.
\label{eq!tLM}
\end{equation}
The latter is anomalous as before: the partial smoothing along the
$y_4$-axis is specified either by the Laurent monomial $y_4^{-3}$ or by a
polynomial equation $x_2^3$ in the ``opposite'' variable $x_2$.

\subsubsection{The 6-fold $V_{ABLM}$} \label{ss!2xt}
We now have two deformations $T \subset V_{AB}$ and $T \subset V_{LM}$
of our tent $T$ to toric 4-folds; we call this a {\em diptych of toric
deformations}. The two panels are quite different: $V_{AB}$ is smooth
along the $x_3$- and $y_4$-axes by \eqref{eq!AB}, but has hypersurface
singularities along the $x_0$- and $y_0$-axes of transverse type
$x_1y_0=A^4B^7$ and $x_0y_1=y_0\1A^7B^{12}$ by \eqref{eq!tag0} and
\eqref{eq!tag-1}. In contrast, $V_{LM}$ smooths the $x_0$- and
$y_0$-axes by \eqref{eq!LM}, but leaves the $x_3$- and $y_4$-axes with
the transverse hypersurface singularities $x_2y_4=L^2M^7$ and
$x_3y_3=y_4^{-3}L^7M^{27}$ of \eqref{eq!tLM}.

Theorem~\ref{th!main} now asserts that these two toric panels fit
together in a 4-parameter deformation $T\subset V_{ABLM}$:
\begin{equation}
\renewcommand{\arraystretch}{1.5}
\renewcommand{\arraycolsep}{.6em}
\begin{array}{ccc}
T & \subset & V_{AB} \\
\bigcap && \bigcap \\
\kern2mm V_{LM} & \subset & V_{ABLM}\kern-4mm
\end{array}
\label{eq!dip}
\end{equation}
More precisely, we build an affine 6-fold $V_{ABLM}$ with a regular
sequence $A,B,L,M$ such that the section $L=M=0$ is $V_{AB}$ and $A=B=0$
is $V_{LM}$. The idea is amazingly naive: starting at the top, we simply
merge the tag equations \eqref{eq!AB} and \eqref{eq!tLM} for $x_3$ and
$y_4$ from $V_{AB}$ and $V_{LM}$, obtaining
$W\subset\Aff^8_{\Span{x_2,x_3,y_4,y_3,A,B,L,M}}$ defined by
\begin{equation}
x_2y_4=x_3^2A+L^2M^7, \quad x_3y_3=y_4B+x_2^3LM^3.
\label{eq!mix}
\end{equation}
It is a codimension~2 complete intersection, $A,B,L,M$ is a regular
sequence for $W$, and the section $L=M=0$ is birational to $V_{AB}$ by
the Laurent monomial argument of \eqref{eq!Laur}.

The plan is now to adjoin $x_1,x_0,y_2,y_1,y_0$ as rational functions on
$W$, so $V_{ABLM}$ will be birational to $W$. In commutative algebra
terms, the coordinate ring of $V_{ABLM}$ is constructed from the
complete intersection \eqref{eq!mix} by {\em serial unprojection}. We
run through the construction as a pleasant narrative; the reasons it all
works include some detailed tricks that we explain later when we treat
the material more formally. Suffice it to say that we add the new
variables $x_1,x_0,y_2,y_1,y_0$ one at a time, {\em and in that order}.
Adding them in a different order does not work.

\subsubsection{First pentagram} We construct $x_1$ as a rational
function on $W$ \eqref{eq!mix} with divisor of poles the codimension~3
complete intersection
\begin{equation}
D:(x_3=y_4=LM^3=0)\subset W,
\label{eq!D}
\end{equation}
where $LM^3$ is the hcf of the two terms $L^2M^7$ and $x_2^3LM^3$ in
\eqref{eq!mix}. The new variable $x_1$ appears in three equations
\begin{equation}
x_1x_3=\cdots, \quad x_1y_4=\cdots, \quad x_1LM^3=\cdots,
\label{eq!x1un}
\end{equation}
that express the rational function $x_1$ as a homomorphism
$\sI_D\to\Oh_W$. More intrinsically, $x_1$ is an {\em unprojection
variable} $x_1\in\sHom(\sI_D,\om_W)$ with Poincar\'e residue a basis of
$\om_D\iso\Oh_D$; see \cite{PR} and \cite{Ki} for the theory and
practice of unprojection. In our calculation we take as input the
equations \eqref{eq!mix} and \eqref{eq!D} of $W$ and $D$, and use them
to fix up a $5\times5$ skew matrix $A=\{a_{ij}\}$ whose five $4\times4$
Pfaffians are the two input equations \eqref{eq!mix} and the three new
unprojection equations \eqref{eq!x1un} for $x_1$. This calculation is
repeated serially in what follows, and we make it systematic with {\em
magic pentagrams}:
\begin{equation} \begin{array}{cl}
\begin{picture}(80,36)(0,25)
\put(0,0){$x_1$} \put(16,3){\circle*5} \put(16,3){\line(1,1){50}}
\put(0,25){$x_2$} \put(16,28){\circle*5} \put(16,28){\line(2,1){50}}
\put(16,28){\line(1,0){50}}
\put(0,50){$x_3$} \put(16,53){\circle*5} \put(16,53){\line(2,-1){50}}
\put(66,28){\circle*5} \put(73,25){$y_3$}
\put(66,53){\circle*5} \put(73,50){$y_4$}
\qbezier(16,3)(36,28)(16,53)
\end{picture}
&\begin{pmatrix}
y_3 & x_2^3 & -B & -x_1 \\
& y_4 & LM^3 & -x_3A \\
&& x_3 & LM^4 \\
&&& x_2
\end{pmatrix}
 \\[1cm]
\kern-1mm \begin{array}{ll}
23.45 & x_2y_4=x_3^2A+L^2M^7, \\
12.34 & x_3y_3=y_4B+x_2^3LM^3,
\end{array}
& \begin{array}{ll}
12.35 & x_1y_4=x_2^3x_3A+y_3LM^4, \\
13.45 & x_1x_3=x_2^4+BLM^4, \\
12.45 & x_2y_3=x_3AB+x_1LM^3.
\end{array} \kern-1mm
\end{array}
\label{penta!1}
\end{equation}
The array is a skew $5\times5$ matrix $A=\{a_{ij}\}$; we only write the
10 upper-triangular entries $a_{12}=y_3,\dots,a_{15}=-x_1$, etc. Its
$4\times4$ Pfaffians are
\begin{equation}
\Pf_{ij.kl}=a_{ij}a_{kl}-a_{ik}a_{jl}+a_{il}a_{jk}
\quad\hbox{for any distinct $i,j,k,l$}
\end{equation}
(as with minors and cofactors, with an overall choice of $\pm1$; in long
calculations we abbreviate $\Pf_{ij.kl}$ to $ij.kl$). In
\eqref{penta!1}, viewing $y_3,y_4,x_3,x_2$ and the two equations
$x_2y_4=\cdots$, $x_3y_3=\cdots$ as given, we seek to add $x_1$ and
three new equations $x_1x_3=\cdots$, $x_1y_4=\cdots$ and
$x_2y_3=\cdots+x_1LM^3$. These trinomial equations play a role for
$V_{ABLM}$ similar to the binomial tag equations
$v_{i-1}v_{i+1}=v_i^{a_i}$ for the cyclic quotient singularities $S_i$
and the tent $T$. The array is written out automatically from the
pentagram and the given equations \eqref{eq!mix}: we write the given
variables $y_3,y_4,x_3,x_2$ down the superdiagonal, the new unprojection
variable $x_1$ in the top right, and the given
$LM^3=\hcf(L^2M^7,x_2^3LM^3)$ as the entry $a_{24}$. Requiring
$\Pf_{12.34}$ and $\Pf_{23.45}$ to give \eqref{eq!mix} determines the
remaining entries. The output is the three equations involving $x_1$ as
the three remaining Pfaffians in \eqref{penta!1}.

\subsubsection{Serial pentagrams} \label{ss!ser}
The remaining variables $x_0,y_2,y_1,y_0$ are adjoined likewise to give
the codimension~7 variety $V_{ABLM}$ (see Section~\ref{s!pf} for a
formal treatment). We write out the calculations without further comment
for your delight.
\begin{equation*} \begin{array}{cl}
\begin{picture}(80,36)(0,21)
\put(0,-7){$x_0$} \put(16,-4){\circle*5}
\put(0,12){$x_1$} \put(16,15){\circle*5} \put(16,15){\line(5,2){50}}
\put(0,31){$x_2$} \put(16,34){\circle*5}
\put(16,34){\line(1,0){50}}
\put(0,50){$x_3$} \put(16,53){\circle*5}
\put(64,34){\circle*5} \put(71,31){$y_3$}
\qbezier(16,15)(34,34)(16,53)
\qbezier(16,-4)(34,15)(16,34)
\qbezier(16,-4)(50,24.5)(16,53)
\end{picture}
& \begin{pmatrix}
y_3 & x_1 & -AB & -x_0 \\
& x_3 & LM^3 & -x_2^3 \\
&& x_2 & BM \\
&&& x_1
\end{pmatrix}
\\[1cm]
\begin{array}{ll}
23.45 & x_1x_3=x_2^4+BLM^4, \\
12.34 & x_2y_3=ABx_3+LM^3x_1,
\end{array}
\enspace
& \begin{array}{ll}
12.35 & x_0x_3=x_1x_2^3+BMy_3, \\
13.45 & x_0x_2=x_1^2+AB^2M, \\
12.45 & x_1y_3=ABx_2^3+LM^3x_0.
\end{array}
\end{array}
\end{equation*}

\begin{equation*}
\renewcommand{\arraycolsep}{0.3em}
\begin{array}{cc}
\begin{picture}(80,36)(0,25)
\put(0,0){$x_0$} \put(16,3){\circle*5} \put(16,3){\line(1,1){50}}
\put(0,25){$x_1$} \put(16,28){\circle*5} \put(16,28){\line(2,1){50}}
\put(16,28){\line(1,0){50}}
\put(0,50){$x_2$} \put(16,53){\circle*5} \put(16,53){\line(2,-1){50}}
\put(66,28){\circle*5} \put(73,25){$y_2$}
\put(66,53){\circle*5} \put(73,50){$y_3$}
\qbezier(16,3)(36,28)(16,53)
\end{picture}
& \begin{pmatrix}
y_3 & LM^2x_0 & -ABx_2^2 & -y_2 \\
& x_2 & M & -x_1 \\
&& x_1 & AB^2 \\
&&& x_0
\end{pmatrix}
\\[1cm]
\begin{array}{ll}
23.45 & x_0x_2=x_1^2+AB^2M, \\
12.34 & x_1y_3=ABx_2^3+LM^3x_0,
\end{array}
& \begin{array}{ll}
12.35 & x_2y_2=AB^2y_3+LM^2x_0x_1, \\
13.45 & x_1y_2=A^2B^3x_2^2+LM^2x_0^2, \\
12.45 & x_0y_3=ABx_1x_2^2+My_2.
\end{array}
\end{array}
\end{equation*}

\begin{equation*} \begin{array}{cc}
\renewcommand{\arraycolsep}{0.4em}
\begin{picture}(80,36)(0,25)
\put(0,0){$x_0$} \put(16,3){\circle*5} \put(16,3){\line(5,4){50}}
\put(0,25){$x_1$} \put(16,28){\circle*5} \put(16,27.5){\line(3,1){50}}
\put(16,27.5){\line(6,-1){50}}
\put(0,50){$x_2$} \put(16,53){\circle*5} \put(16,53){\line(3,-2){50}}
\put(67,19){\circle*5} \put(74,16){$y_1$}
\put(67,43.5){\circle*5} \put(74,40.5){$y_2$}
\qbezier(16,3)(40,28)(16,53)
\end{picture}
& \begin{pmatrix}
y_2 & LMx_0^2 & -A^2B^3x_2 & -y_1 \\
& x_2 & M & -x_1 \\
&& x_1 & AB^2 \\
&&& x_0
\end{pmatrix}
\\[1cm]
\begin{array}{ll}
23.45 & x_0x_2=x_1^2+AB^2M, \\
12.34 & x_1y_2=A^2B^3x_2^2+LM^2x_0^2,
\end{array}
\enspace
& \begin{array}{ll}
12.35 & x_2y_1=AB^2y_2+LMx_0^2x_1, \\
13.45 & x_1y_1=A^3B^5x_2+LMx_0^3, \\
12.45 & x_0y_2=A^2B^3x_1x_2+My_1.
\end{array}
\end{array}
\end{equation*}

\begin{equation*} \begin{array}{cc}
\renewcommand{\arraycolsep}{0.4em}
\begin{picture}(80,36)(0,25)
\put(0,0){$x_0$} \put(16,3){\circle*5} \put(16,3){\line(2,1){50}}
\put(0,25){$x_1$} \put(16,28){\circle*5} \put(16,28){\line(2,-1){50}}
\put(16,28){\line(1,0){50}}
\put(0,50){$x_2$} \put(16,53){\circle*5} \put(16,53){\line(1,-1){50}}
\put(66,28){\circle*5} \put(73,25){$y_1$}
\put(66,3){\circle*5} \put(73,0){$y_0$}
\qbezier(16,3)(36,28)(16,53)
\end{picture} \quad
& \begin{pmatrix}
y_1 & Lx_0^3 & -A^3B^5 & -y_0 \\
& x_2 & M & -x_1 \\
&& x_1 & AB^2 \\
&&& x_0
\end{pmatrix}
\\[1cm]
\begin{array}{ll}
23.45 & x_0x_2=x_1^2+AB^2M, \\
12.34 & x_1y_1=A^3B^5x_2+LMx_0^3,
\end{array}
\enspace
& \begin{array}{ll}
12.35 & x_2y_0=AB^2y_1+Lx_0^3x_1, \\
13.45 & x_1y_0=A^4B^7+Lx_0^4, \\
12.45 & x_0y_1=A^3B^5x_1+My_0.
\end{array}
\end{array}
\end{equation*}

The final two equations $x_1y_0=\cdots$ and $x_0y_1=\cdots$ merge the
tag equations (\ref{eq!tag0}--\ref{eq!tag-1}) and \eqref{eq!LM} for
$x_0$ and $y_0$ at the bottom of the two long rectangles in exactly the
same way as \eqref{eq!mix} merged the tag equations at the top. In other
words, the whole calculation could have been done starting with these
two equations and working up -- if you liked the puzzle, you will enjoy
turning it upside down and doing it all over again.

\section{Toric partial smoothings of tents} \label{s!trc}
This chapter centres around the combinatorics of continued fractions.
After recalling standard facts, we define a tent $T$, and, under
appropriate assumptions, construct a toric extension $T\subset V_{AB}$
that smooths its top two axes. The toric variety $T\subset V_{AB}$ can
be treated in terms of a matrix $\left(\!
\begin{smallmatrix}r&a\\b&s\end{smallmatrix} \!\right) \in\SL(2,\Z)$, or
equivalently, in terms of a certain continued fraction expansion of $0$.
We use the latter treatment in \ref{ss!pseq} to understand $V_{AB}$ by a
sequence of Gorenstein projections.

\subsection{Jung--Hirzebruch continued fractions} \label{ss!HJ}

A {\em continued fraction expansion} is a formal expression
\begin{equation}
\begin{split}
[c_1,\dots,c_n] &= c_1-1/(c_2-1/(c_3-\cdots-1/c_n)\cdots) \\
&= c_1-\frac1{c_2-\frac1{\dots-\frac1{c_n}}}
= c_1-\frac1{[c_2,\dots,c_n]}
\end{split}
\label{eq!HJ}
\end{equation}
The entries $c_i$ are called {\em tags}. If $c_1,\dots, c_n$ are
integers, the righthand side is a rational number, provided that the
expression makes sense, that is, division by zero does not
occur. (The notation is explained in Riemenschneider \cite{Rie}
\S3, pp.~220--3.)

The next proposition discusses four aspects of continued fractions. We
spell out this material, because we use it often and with large
multi\-plicity in what follows: we invert continued fractions and pass
to complementary fractions, we ``top and tail'' them by cutting off a
tag at one end and adding one at the other, say:
\begin{equation}
[a_0,\dots,a_{k-1}] \mapsto [a_k,a_{k-1},\dots,a_1], \quad\hbox{etc.,}
\end{equation}
and we concatenate the resulting fractions.

\begin{prop} \label{prop!hj}
\begin{description}

\item[(a) Factoring a matrix:] The formal identity
\begin{equation}
\begin{pmatrix} 0 & 1 \\ -1 & c_1 \end{pmatrix}
\begin{pmatrix} 0 & 1 \\ -1 & c_2 \end{pmatrix}
\cdots \begin{pmatrix} 0 & 1 \\ -1 & c_n \end{pmatrix}
= \begin{pmatrix} -q' & q \\ -p' & p \end{pmatrix}.
\label{eq!prod}
\end{equation}
holds in indeterminates or variables $c_1,\dots,c_n$, where $p,q,p',q'$
are polynomials, the numerators and denominators of
$p/q=[c_1,\dots,c_n]$ and $p'/q'=[c_1,\dots,c_{n-1}]$. (No cancellation
occurs in the fraction $p/q$, whatever the nature or values of the
quantities $c_i$, because $p$ and $q$ satisfy an $\hcf$ identity $\al
p+\be q=1$.) The fraction $p'/q'$ is the {\em first convergent} of $p/q$.

\item[(b) Blowdown:] $[c_1,\dots,c_{n-1},1]=[c_1,\dots,c_{n-1}-1]$ and
\begin{equation}
[c_1,\dots,c_{i-1},1,c_{i+1},\dots,c_n] =
[c_1,\dots,c_{i-1}-1,c_{i+1}-1,\dots,c_n].
\label{eq!bld}
\end{equation}
This is just the identity
$\left(\begin{smallmatrix} 0 & 1 \\ -1 & a \end{smallmatrix}\right)
\left(\begin{smallmatrix} 0 & 1 \\ -1 & 1 \end{smallmatrix}\right)
\left(\begin{smallmatrix} 0 & 1 \\ -1 & b \end{smallmatrix}\right)
= \left(\begin{smallmatrix} 0 & 1 \\ -1 & a-1 \end{smallmatrix}\right)
\left(\begin{smallmatrix} 0 & 1 \\ -1 & b-1 \end{smallmatrix}\right)$.
\item[] Two notions of ``inverse'' of a continued fraction play a role
in our theory:
\item[(c) Reciprocal:] $[c_1,\dots,c_n]=p/q$ and its {\em
reciprocal} continued fraction
\begin{equation}
[c_n,\dots,c_1]=p/q^*
\end{equation}
share the same numerator $p$, and their denominators are inverse modulo
$p$. More precisely, there is a formal identity
\begin{equation}
qq^*=N(c_2,\dots,c_{n-1})\cdot p+1,
\end{equation}
where $N(c_2,\dots,c_{n-1})$ is the numerator of $[c_2,\dots,c_{n-1}]$.
In particular, if $c_i\in\Z$ and the expressions are meaningful then
$[c_n,\dots,c_1]=p/q^*$, where $qq^*\equiv 1\mod p$. See \eqref{eq!opp}
for what this means in our context.

\item[(d) Complement:] Let $p/q=[c_1,\dots,c_n]$ with $c_i\in\Z$ and
$c_i\ge2$. Then the {\em complementary} continued fraction is
$[b_1,\dots,b_m]=p/(p-q)$, and satisfies
\begin{equation}
[c_n,\dots,c_1,1,b_1,\dots,b_m]=0.
\end{equation}
Moreover, serial blowdown reduces the expansion to $[1,1]=[0]=0$; in
particular, $\sum (c_i-1)=\sum (b_j-1)$, and one of $b_1,c_1\le2$.
For example,
\begin{equation}
[4,\underline{2,1,3},2,2]=[\underline{4,1,2},2,2]
=[\underline{3,1,2},2]=[\underline{2,1,2}]=[1,1]=0.
\end{equation}
\end{description}
\end{prop}

\begin{remark} \label{rk!JHc}
Traditionally, one uses Jung--Hirzebruch continued fractions to
write a fraction $\frac ra$ with $r>a\ge1$ and $a,r$ coprime
integers as
\begin{equation}
\frac ra=[b_1,\dots,b_{n-1}]=b_1-\frac1{b_2-\cdots}\,.
\notag
\end{equation}
Then $b_1$ is the round-up $b_1=\rdup{\frac ra}$, and is $\ge2$, because
$\frac ra>1$, and for the same reason all subsequent $b_i\ge2$ (to the
end of the algorithm). Here we do something slightly bigger, with
$a\ge1$, but $r\in\Z$ any integer coprime to $a$: for example,
$\frac{-24}7=-3-\frac37=[-3,3,2,2]$. This means that $b_1=\rdup{\frac
ra}\in\Z$; however, from the second step onwards and to the end of the
algorithm, $1/(b_1-\frac ra)>1$ is a conventional fraction, so that
$b_i\ge2$ for each $i$ with $2\le i\le n-1$.

In traditional use, \eqref{eq!prod} identifies 3 types of data: a
rational fraction $p/q>1$, a continued fraction $[c_1,\dots,c_n]$ with all
$c_i\ge2$, and a matrix $\bigl(\!\begin{smallmatrix} -q'&q \\-p' &
p\end{smallmatrix}\bigr)\in\SL(2,\Z)$ with $p>q>0$. However, we relax
these restrictions, considering things like $[5,1,3]= 5-\frac32=
\frac72= [4,2]$ (a blowdown) or $[2,0,2]=4$, with
\begin{equation}
\label{eq!202}
\begin{pmatrix} 0 & 1 \\ -1 & 2 \end{pmatrix}
\begin{pmatrix} 0 & 1 \\ -1 & 0 \end{pmatrix}
\begin{pmatrix} 0 & 1 \\ -1 & 2 \end{pmatrix}
= \begin{pmatrix} 0 & 1 \\ -1 & 2 \end{pmatrix}
\begin{pmatrix} -1 & 2 \\ 0 & -1 \end{pmatrix}
= \begin{pmatrix} 0 & -1 \\ 1 & -4 \end{pmatrix}.
\end{equation}
The matrix product \eqref{eq!prod} is meaningful even when \eqref{eq!HJ}
involves division by zero. More generally, the sequence of integer tags
$[c_1,\dots,c_n]$ contains more information than the matrix
\eqref{eq!prod}, which contains more information than the fraction
$\frac pq$: while $\bigl(\! \begin{smallmatrix} -q'&q\\-p'& p
\end{smallmatrix} \bigr)\in\SL(2,\Z)$, the fraction $\frac pq$ (when
defined) is its image in the quotient group $\PSL(2,\Z)$, whereas the
expression $[c_1,\dots,c_n]$ is a lift to the ``universal cover'' of
$\SL(2,\Z)$ inside the universal cover of $\SL(2,\R)$, keeping track of
winding number. For example, $[0,0,0,0]$ is the composite of 4 rotations
by $\pi/2$, or $\bigl(\! \begin{smallmatrix} 0&1\\-1& 0\end{smallmatrix}
\bigr)^4=\id$. Running around one of our long rectangles below always gives
winding number~1.
\end{remark}

\paragraph{Notation for the quotient $\frac1r(\al,1)$}
As in Example \ref{ex!intr}, for $r\ge 1$ and $0<\al\le r$ coprime to $r$
we write $\frac1r(\al,1)$ for the $\Z/r$ action on $\Aff^2_{\Span{u,v}}$ given by
$(u,v)\mapsto(\ep^\al u,\ep v)$, and for the quotient
$S=\Aff^2/\frac1r(\al,1)$ by this action. We allow $\Aff^2$ as the case
$r=1$, without worrying unduly about the value of $\al$ (of course,
$\al=0$); it corresponds to the identity matrix or the empty continued
fraction $[\emptyset]$. The lattice $\La$ of invariant Laurent monomials consists
of $u^iv^j$ with $\al i+j\equiv0$ \hbox{mod $r$}; it is a lattice $\La\cong\Z^2$,
but with no preferred basis. The coordinate ring of $S$, based by
$\Z/r$-invariant monomials, is minimally generated by monomials on the
Newton boundary of the positive quadrant $\si\subset \La_\R$.
Setting $0<\be\le r$ with $\al\be=-1$ modulo~$r$, these monomials
are $x_0=u^r$, $x_1=u^\be v$, etc. Either continued fraction
$\frac r{\be}=[a_1,\dots,a_{k-1}]$ or
$\frac r{r-\al}=[a_{k-1},\dots,a_1]$ provides the generators $x_{0\dots k}$ and
the {\em tag equations} holding between them:
\begin{equation}
x_{i-1}x_{i+1}=x_i^{a_i} \quad\hbox{for $i=1,\dots,k-1$.}
\label{eq!tag}
\end{equation}
In particular,
\begin{equation}
\begin{array}{rcl}
\displaystyle\frac r{\be}=[a_1,\dots,a_{k-1}] &\mapsto &
x_0=u^r,\ x_1=u^{\be}v,\ x_2=x_1^{a_1}x_0\1,\dots \\[10pt]
\displaystyle\frac r{r-\al}=[a_{k-1},\dots,a_1] &\mapsto &
x_k=v^r,\ x_{k-1}=uv^{r-\al},\dots.
\end{array}
\label{eq!opp}
\end{equation}
The tag equations \eqref{eq!tag} determine $S$ completely: they express any $x_j$ as a Laurent
monomial in any two consecutive monomials $x_i,x_{i+1}$. The complete
intersection in $\Aff^{k+1}_{\Span{x_{0\dots k}}}$ given by
\eqref{eq!tag} is $S$ plus $\Aff^2_{\Span{x_0,x_k}}$ (usually with a
nonreduced structure). The other generators of $I_S$ are ``long
equations'' $x_ix_j=$ monomial for $|i-j|>2$, that can be deduced from
\eqref{eq!tag} via syzygies.

\subsection{Tents and fans}\label{ss!tnt}

\begin{figure}[h] \[
\setlength{\unitlength}{1.5mm}
\begin{picture}(27,26)
\thicklines
\put(5,1){\line(1,0){17}}
\put(5,1){\line(0,1){24}}
\put(5,25){\line(1,0){17}}
\put(22,1){\line(0,1){24}}
\thinlines
\put(5,1){\circle*{1.3}} %% x0
\put(5,6){\circle*{1.3}} %% x1
\put(5,20){\circle*{1.3}} %% xl-1
\put(5,25){\circle*{1.3}} %% xl
\put(22,1){\circle*{1.3}} %% xn
\put(22,6){\circle*{1.3}} %% y1
\put(22,25){\circle*{1.3}} %% xl+2
\put(22,20){\circle*{1.3}} %% xl+1
\put(1,0.5){$x_0$}
\put(1,5.5){$x_1$}
\put(7,5.5){$a_1$}
\put(-2,19.5){$x_{k-1}$}
\put(7,19.5){$a_{k-1}$}
\put(1,24.5){$x_k$}
\put(24,24.5){$y_l$}
\put(24,19.5){$y_{l-1}$}
\put(16,19.5){$b_{l-1}$}
\put(24,5.5){$y_1$}
\put(18,5.5){$b_1$}
\put(24,0.5){$y_0$}
\put(3,12){\vdots}
\put(6,12){\vdots}
\put(20,12){\vdots}
\put(23,12){\vdots}
\end{picture} \]
\caption{Coordinates and tags for a tent $T$ \label{f!tnt}}
\end{figure}
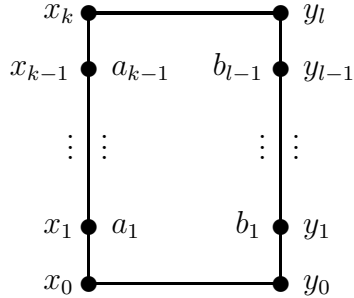

A {\em tent} $T=S_0\cup S_1\cup S_2\cup S_3 \subset
\Aff^{k+l+2}_{\Span{x_{0\dots k},y_{0\dots l}}}$ is the union of the
four affine toric surfaces of Figure~\ref{f!T}, with horizontal sides
$S_0=\Aff^2_{\Span{x_0,y_0}}$ and $S_2=\Aff^2_{\Span{x_k,y_l}}$ and
vertical sides the cyclic quotient singularities
\begin{align*}
&\hbox{$S_1=\frac1r(\al,1)$ with coordinates $x_{k\dots 0}$ from
$\frac{r}{r-\al}=[a_{k-1},\dots,a_1]$, and} \\
&\hbox{$S_3=\frac1s(\be,1)$ with coordinates $y_{l\dots 0}$ from
$\frac{s}{s-\be}=[b_{l-1},\dots,b_1]$},
\end{align*}
where $\al\le r$ are coprime natural numbers, and similarly for $\be\le
s$ (there are no other conditions on $\al$, $\be$ at this stage). The
coordinates $x_{0\dots k},y_{0\dots l}$ of the ambient space
$\Aff^{k+l+2}$ and the equations for $T$ are shown schematically in
Figure~\ref{f!tnt}; once we have added {\em corner tags} in \ref{ss!fan}
and {\em annotations} in \ref{ss!bet}, we refer to such arrays as {\em
long rectangles}, and use them as a shorthand for certain toric 4-folds.
The components glue transversally along their toric strata (= coordinate
axes), giving $T$ four singular axes of transverse ordinary double
points; the two axes on $S_2$ are the {\em top} axes, and the two on
$S_0$ the {\em bottom} axes.

\subsubsection{Tents without embeddings}
Our definition expresses $T$ embedded in $\Aff^{k+l+2}$ by explicit
coordinates; its ideal $I_T$ is generated by $I_{S_1}$ and $I_{S_3}$,
determined by the tags down the sides as in \ref{ss!HJ}, together with
the cross-equations $x_iy_j=0$ for all pairs $(i,j)\ne(0,0),(k,l)$.

However, $T$ can be viewed abstractly as an identification scheme as
studied more generally in Reid \cite{dP}: write $\Ga'_i\cup\Ga''_i$ for
the toric 1-strata of the $S_i$ and
$C=\bigsqcup_{i=1}^4(\Ga'_i\cup\Ga''_i)$. Let $D$ be the four axes
$\Aff^1$ with coordinates $x_0,x_k,y_l,y_0$ glued transversally at a
common origin (as coordinate axes in $\Aff^4$); write $\fie\colon C\to
D$ for the morphism given by $x_0$ on the $x_0$-axes of $S_0$ and $S_1$,
and so on, to perform the identifications of Figure~\ref{f!T}. Then
\begin{equation}
T=\bigl(S_0\sqcup S_1\sqcup S_2\sqcup S_3\bigr)/\fie.
\end{equation}
There are no parameters or moduli in this glueing.

\begin{lem} \label{l!tnt}
Let $T$ be the tent as above. Then $T$ is a
Gorenstein scheme. Moreover, $T$ has an action of $(\Gm)^4$ that
restricts to the toric structure on each component. \end{lem}

\begin{pf} We use elementary results of \cite{dP}, Section~2. $T$ is
Cohen--Macaulay because all the glueing happens in codimension~1
(\cite{dP}, 2.2). We prove it is Gorenstein using the criterion of
\cite{dP}, Corollary~2.8.

Each component $S_i$ is a toric surface; on each, choose a $\Z$-basis
$m_1,m_2$ for the monomial lattice, oriented clockwise (e.g., on $S_1$,
take $x_0,x_1$ or $x_{k-1},x_k$; on $S_2$, take $x_k,y_l$). The 2-form
$s=\frac{\dd m_1}{m_1}\wedge\frac{\dd m_2}{m_2}\in\Om^2_\T$ on the big
torus is a basis for $\Om^2_\T$, is defined over $\Z$, independent of
the choice of oriented basis, and has log poles along each stratum of
$S$, with residue along each stratum $\Aff^1$ equal to $\pm$ times the
natural basis $\frac{\dd m}{m}$ of $\Om^1_{\T'}$. We take this basis
element $s$ on each component. Under the identification $\fie\colon C\to
D$ of the double locus, over the general point of each component of $D$,
the residues from the two components are $\pm\frac{\dd m}{m}$, and
therefore cancel out; thus $s$ satisfies the conditions of \cite{dP},
Corollary~2.8.ii and is a basis of the dualising sheaf $\om_T$.

Each component of $T$ is a toric variety, so $(\Gm)^8$ acts on the
disjoint union of the components. Each glueing imposes one linear
condition on the action; we think of
$\T_{S_0}=(\Gm)^2=\{(\la_0,1,1,\la_3)\}$ as the big torus of $S_0$ and
$\T_{S_1}=\{(\la_0,\la_1,1,1)\}$ that of $S_1$, etc. \end{pf}

\subsubsection{The fan $\Phi(\protect\begin{smallmatrix}
r&a\protect\\b&s \protect\end{smallmatrix})$ in the plane given by
$(\protect\begin{smallmatrix} r&a\protect\\b&s
\protect\end{smallmatrix})\in\SL(2,\Z)$} \label{ss!fan}

Jung--Hirzebruch continued fractions factor a base change in $\SL(2,\Z)$
into elementary moves (Proposition~\ref{prop!hj}(a)); in our case, the
base change goes from the monomials $x_0,y_0$ at the bottom of our long
rectangle to $x_k,y_l$ at the top (up to sign and orientation).
\ref{ss!bet} constructs the toric variety $V_{AB}$ and the first
extension $T\subset V_{AB}$ generalising \eqref{eq!Laur}, using a matrix
in $\SL(2,\Z)$ to generate the monomial cone $\si_{AB}$ of
Figure~\ref{f!siAB} in the 4-dimensional lattice $\M=\Z^4$.

We start by analysing the combinatorics of this construction in a
stripped-down 2-dimensional setting $\bMbar=\Z^2$.
Consider two oriented bases $x_0,y_0$ and $\eta,\xi$ of
$\bMbar$ related by inverse base changes
\begin{equation}
x_0=\eta^{-r} \xi^a, \quad y_0=\eta^b\xi^{-s}
\quad\hbox{and}\quad
\eta=x_0^{-s}y_0^{-a},\quad \xi=x_0^{-b}y_0^{-r}.
\label{eq!x0}
\end{equation}
Here $r,s,a,b\ge0$ are integers with $rs-ab=1$, so
\begin{equation}
\begin{pmatrix} r&a\\b&s
\end{pmatrix}
\quad\hbox{and}\quad
\begin{pmatrix} s&-a\\-b&r \end{pmatrix}\in\SL(2,\Z)
\label{eq!mats}
\end{equation}
are a pair of inverse elements. (If $a$ or $b=0$ then $r=s=1$, and one
or two points in what follows need minor restatement. Rather than do
that systematically, it is easier simply to list all these initial cases,
as in~\ref{ss!init}.)

The vectors $x_0,y_0,\eta,\xi$ subdivide the plane $\bMbar_\R$ into the
fan $\Phi(\begin{smallmatrix} r&a\\b&s \end{smallmatrix})$ of
Figure~\ref{f!Phi}.a consisting of 4 cones $\Span{x_0,y_0}$,
$\Span{x_0,\xi}$, $\Span{\xi,\eta}$, $\Span{y_0,\eta}$. It determines a
tent $T$, with coordinate ring generated by the 4 monomial cones and
related by $m_1m_2=0$ if $m_1,m_2$ are not in a common cone.
The next lemma computes the affine toric surfaces that make up the tent $T$
corresponding to $\Phi(\begin{smallmatrix}r&a\\b&s\end{smallmatrix})$;
compare with the first long rectangle of Example~\ref{ex!intr} for which
$\left(\begin{smallmatrix} r&a\\b&s
\end{smallmatrix}\right)=\left(\begin{smallmatrix} 7&12\\4&7
\end{smallmatrix}\right)$.

\begin{lem} \label{l!wT}
Suppose that $r,s,a,b\ge 1$. Consider the cone
$\Span{x_0,\xi}$ (marked $S_1$ in Figure~\ref{f!Phi}(a)). The lattice
$\bMbar$ is generated by the monomials $x_0,\xi$ together with either of
\begin{equation}
y_0\1=(x_0^b\xi)^{1/r} \quad\hbox{or}\quad
\eta=(x_0\1\xi^a)^{1/r}.
\notag
\end{equation}
Therefore $\Span{x_0,\xi}$ is the monomial cone $\frac1r(\al,1)$ or
$\frac1r(1,r-\be)$, where $\al$ is the least residue of\/ $a$ {\rm
mod}~$r$, and $\be$ that of\/ $b$ (note that $rs-ab=1$ implies $\al$ and
$r-\be$ are inverse {\rm mod}~$r$).

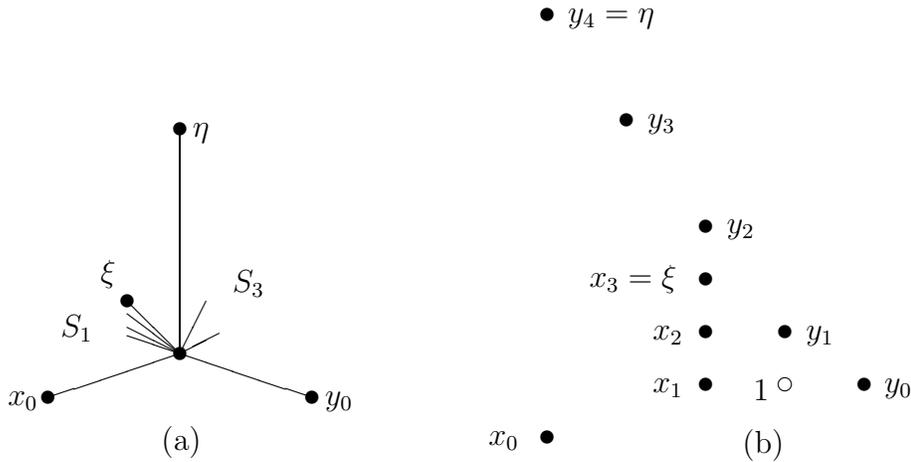
\begin{figure}[ht]
\begin{picture}(100,100)(-40,-20)
\put(-15,-3){$x_0$} \put(0,0){\circle*{5}} \put(0,0){\line(3,1){50}}
\put(105,-3){$y_0$} \put(100,0){\circle*{5}} \put(100,0){\line(-3,1){50}}
\put(50,16.667){\circle*{5}} \put(50,16.667){\line(-1,1){20}}
\put(20,43.667){$\xi$} \put(30,36.667){\circle*{5}}
\put(50,16.667){\line(1,2){10}}
\put(50,16.667){\line(2,1){15}}
\put(50,16.667){\line(0,1){85}}
\put(55,98.667){$\eta$} \put(50,101.667){\circle*{5}}
\put(70,40){$S_3$}
\put(50,16.667){\line(-2,1){20}}
\put(50,16.667){\line(-3,1){20}}
\put(50,16.667){\line(-4,3){20}}
\put(5,23){$S_1$}
\put(43,-20){(a)}
\end{picture}
\begin{picture}(150,168)(-215,-25)
\put(-12,-6){1} \put(0,0){\circle{5}}
\put(-16,-26){(b)}
\put(37.5,-3){$y_0$} \put(30,0){\circle*{5}}
\put(8,17){$y_1$} \put(0,20){\circle*{5}}
\put(-22,57){$y_2$} \put(-30,60){\circle*{5}}
\put(-52,97){$y_3$} \put(-60,100){\circle*{5}}
\put(-82,137){$y_4=\eta$} \put(-90,140){\circle*{5}}
\put(-74,37){$x_3=\xi$} \put(-30,40){\circle*{5}}
\put(-50,17){$x_2$} \put(-30,20){\circle*{5}}
\put(-50,-3){$x_1$} \put(-30,0){\circle*{5}}
\put(-112,-23){$x_0$} \put(-90,-20){\circle*{5}}
\end{picture}
\caption{The fan $\Phi(\protect\begin{smallmatrix} r&a\protect\\b&s
\protect\end{smallmatrix})\in\SL(2,\Z)$ defined by $x_0,y_0,\eta,\xi$}
\label{f!Phi}
\end{figure}
Write $x_0,x_1,\dots,x_{k-1},x_k=\xi$ for the successive monomials along
the Newton boundary of $\Span{x_0,\xi}$. The number $k$ and the
monomials themselves come from factoring the base change \eqref{eq!x0}
into elementary moves:
\begin{equation}
\begin{pmatrix} -r&a\\b&-s \end{pmatrix}
=\begin{pmatrix} 0&1\\-1&a_0 \end{pmatrix}
\begin{pmatrix} 0&1\\-1&a_1 \end{pmatrix}\cdots
\begin{pmatrix} 0&1\\-1&a_{k-1} \end{pmatrix}
\begin{pmatrix} 0&1\\-1&a_k \end{pmatrix},
\end{equation}
in which each of $a_1,\dots,a_{k-1}\ge2$.
More concretely, they are given by the continued fraction expansions
\begin{equation}
[a_0,a_1,\dots,a_{k-1}]=\frac{-b}r \quad\hbox{and}\quad
[a_k,\dots,a_1]=\frac ar
\label{2eq!b/s}
\end{equation}
by either of the following constructions:
\begin{enumerate}
\renewcommand{\labelenumi}{(\arabic{enumi})}
\item From the bottom, $x_0$ is given, and $x_1=(x_0^\be\xi)^{1/r}$,
where $\be$ is the least residue of $b$ {\rm mod} $r$. Thus
$a_0=\rdup{\frac {-b}r}=\frac{-b+\be}r\le0$ and
\begin{equation}
x_1=(x_0^\be\xi)^{1/r}=y_0\1x_0^{a_0},
\quad\hbox{that is,}\quad
x_1y_0=x_0^{a_0}.
\label{eq!ctg}
\end{equation}
If $\be=0$ then $r$ divides $b$, whereas $rs-ab=1$ implies that $r,b$
are coprime; thus $r=1$, so that $k=1$ and $x_1=\xi$. Otherwise
$x_2,\dots,x_k$ are determined as usual by tag equations
\begin{equation}
x_{i-1}x_{i+1}=x_i^{a_i} \quad\hbox{for $i=1,\dots,k-1$,}
\notag
\end{equation}
where $[a_1,\dots,a_{k-1}]=\frac r\be$ (see Remark~\ref{rk!JHc}).

\item From the top, $x_k=\xi$ is given; if $r\mid a$ then, as before,
$r=1$ and the only monomials are $x_0,x_1=\xi$. Otherwise, set
$x_{k-1}=(x_0\xi^{r-\al})^{1/r}$, where $\al$ is the least residue of
$a$ {\rm mod} $r$. Then $r-\al=a_kr-a$ where $a_k=\rdup{\frac ar}\ge1$, and
\begin{equation}
x_{k-1}=\xi^{a_k}(x_0\xi^{-a})^{1/r}=x_k^{a_k}\eta\1
\quad\hbox{that is,}\quad
x_{k-1}\eta=x_k^{a_k}.
\notag
\end{equation}
The remaining monomials are determined by
\begin{equation}
x_{i-1}x_{i+1}=x_i^{a_i}, \quad\hbox{where $[a_{k-1},\dots,a_1]=\frac
r{r-\al}$.}
\notag
\end{equation}
\end{enumerate}

In the same way, the sequence $[b_0,b_1,\dots,b_l]$ factors the inverse
trans\-formation of \eqref{eq!x0} into elementary moves:
\begin{equation}
\begin{pmatrix} -s&-a\\-b&-r \end{pmatrix}
=\begin{pmatrix} 0&1\\-1&b_l \end{pmatrix}
\begin{pmatrix} 0&1\\-1&b_{l-1} \end{pmatrix}\cdots
\begin{pmatrix} 0&1\\-1&b_1 \end{pmatrix}
\begin{pmatrix} 0&1\\-1&b_0 \end{pmatrix}.
\end{equation}
More concretely, $\Span{y_0,\eta}$ is the monomial cone $\frac1s(b,1)$
or $\frac1s(1,-a)$, and the tags and monomials on the $S_3$ side are
$b_0,\dots,b_l$ and $y_0,\dots,y_k$, given by
\begin{equation}
[b_0,b_1,\dots,b_{l-1}]=\frac{-a}s \quad\hbox{and}\quad
[b_l,\dots,b_1]=\frac bs
\end{equation}
and $y_1=x_0\1y^{b_0}=(y_0^\be\eta)^{1/s}$ where $\be$ is the least
residue of $b$ {\rm mod} $s$.
\end{lem}

Not every tent $T=S_0\cup S_1\cup S_2\cup S_3$ is given by a fan
$\Phi(\begin{smallmatrix} r&a\\b&s \end{smallmatrix})$. Which are?
And in how many ways? What extra data does the fan know about
beyond $T$? The tent $T$ knows the 4 monomial cones up to $\SL(2,\Z)$
isomorphism, but does not know how they fit together in $\Z^2$; it
knows the fractions $\frac1r(\al,1)$ and $\frac1s(\be,1)$, but not the
corner tags $a_0,b_0,a_k,b_l$.

\begin{cor}
The fan $\Phi(\begin{smallmatrix} r&a\\b&s \end{smallmatrix})$ gives $T$
with $S_1=\frac1r(\al,1)$, $S_3=\frac1s(\be,1)$ by the construction of
\ref{ss!fan} if and only if $a\equiv\al$ \hbox{mod}~$r$ and $b\equiv\be$
\hbox{mod}~$s$.

For fixed $T$, except for initial cases with $r=s=1$ (see
\ref{ss!init}), there are $0$, $1$ or $2$ matrixes for which
$\Phi(\begin{smallmatrix} r&a\\b&s \end{smallmatrix})$ gives $T$:
\begin{itemize}
\renewcommand{\itemsep}{0pt}
\item if neither $\al$ nor $\be$ divides $rs-1$, there are none;
\item if $\al$ divides $rs-1$ then $a=\al$, $b=(rs-1)/\al$ provides a
solution;
\item similarly, if $\be\mid(rs-1)$ then $a=(rs-1)/\be$, $b=\be$
provides a solution.
\end{itemize}
\end{cor}

\begin{remark} Whereas Figure~\ref{f!Phi}(a)  sketches the division of
the plane into 4 cones $\Span{x_0,y_0}$, $\Span{x_0,\xi}$,
$\Span{\xi,\eta}$, $\Span{y_0,\eta}$, \ref{f!Phi}(b) accurately plots the
monomials in the case $(\begin{smallmatrix} r&a\\b&s \end{smallmatrix})=
(\begin{smallmatrix} 7&24\\2&7\end{smallmatrix})$, with tags
$[a_0,\dots,a_3]= [0,4,2,4]$ at the $x_i$ and $[b_0,\dots,b_4]=
[-3,3,2,2,1]$ at the $y_i$. The comparison of the rich and messy reality
of \ref{f!Phi}(b) with our square-cut projective pictures such as
Figures~\ref{f!lrAB}--\ref{f!lrLM} and \ref{f!tnt} is startling but enlightening:
it reveals, for example,
\begin{align*}
\hbox{$a_0=0$ at $x_0$} &\implies
\hbox{$x_1,0,y_0$ are in arithmetic progression;} \\
\hbox{$b_4=1$ at $y_4$} &\implies
\hbox{$0x_3y_4y_3$ is a parallelogram;} \\
\hbox{$b_0=-3$ at $y_0$} &\implies
\hbox{1 is in the affine convex hull $1\in\Span{x_0,y_1,y_0}$,}
\end{align*}
and so on. The figure and its monomials have other convexity and
collinearity properties to which we return later (compare the Scissors
of Figure~\ref{f!sciss}).
\end{remark}

\subsubsection{Big end, little end, and attitude of a long rectangle}
\label{ss!at}

In the $\SL(2,\Z)$ geometry of the plane, all basic cones are
equivalent, so there is of course no notion of the {\em size} of an
angle. Despite this, the bottom cone $\Span{x_0,y_0}$ is clearly the
{\em big end} of the fan $\Phi$ in Figure~\ref{f!Phi}: if we view $\Phi$ as a
pie chart, $\Span{x_0,y_0}$ occupies the lion's share of the plane,
practically 50\%. The issue is not size, but convexity. Our choice of
signs in \eqref{eq!x0} is equivalent to
\begin{equation}
-\Span{\xi,\eta}\subseteq\Span{x_0,y_0}.
\label{eq!big}
\end{equation}
Even more holds: every monomial appearing as a minimal generator in the
other cones has inverse in $\Span{x_0,y_0}$.

Our choices in $\Phi$ have already decided that the bottom
$S_0=\Aff^2_{\Span{x_0,y_0}}$ is its big end and the top
$S_2=\Aff^2_{\Span{\xi,\eta}}$ its little end. (The two players will swap
ends for the second half of the game.) Once this choice is out of the
way, there are still two dichotomies for the corner tags, forming a
division into 4 cases, the {\em attitude} of the long rectangle and of
the panel $V_{AB}$. Treating this carefully here will save many
headaches later.

\begin{cor} \label{cor!at}
Except for initial cases with $r$ or $s=1$ (see \ref{ss!init})
$r,s\ne a,b$ and
\begin{equation}
r<a \iff b<s \quad\hbox{and}\quad r<b \iff a<s.
\notag
\end{equation}
The long rectangle $\si_{AB}$ thus has {\em attitude:}
\begin{tabbing}
{\bf Top tags:} \qquad \= either\qquad\= $a_k\ge2$ and $b_l=1$ if $r<a$ and $b<s$; or\\
 \>\> $a_k=1$ and $b_l\ge2$ if $r>a$ and $b>s$; and\\
  \\
{\bf Bottom tags:} \>either\> $a_0\le-1$ and $b_0=0$ if $r<b$ and $a<s$; or\\
 \>\> $a_0=0$ and $b_0\le-1$ if $r>b$ and $a>s$.
\end{tabbing}
\end{cor}

\begin{cor} \label{cor!VAB}
If $a_0<0$ and $b_0=0$ then $[a_2,\dots,a_k,b_l,\dots,b_2,b_1]=0$.
If $a_0=0$ and $b_0<0$ then $[a_1,\dots,a_k,b_l,\dots,b_2]=0$.

Conversely, given $\frac1r(\al,1)$ and $\frac1s(\be,1)$,
the tent $T$ is given by a fan
$\Phi(\begin{smallmatrix} r&a\\b&s \end{smallmatrix})$ with
big end $S_0=\Aff^2_{\Span{x_0,y_0}}$ if and only if
the continued fractions
\begin{equation}
\frac{r}{r-\al}=[a_{k-1},\dots,a_1]
\quad\hbox{and}\quad
\frac{s}{s-\be}=[b_{l-1},\dots,b_1]
\notag
\end{equation}
can be concatenated with $a_k$ and $b_l$ such that
\begin{equation}
\hbox{either \quad $[a_2,\dots,a_k,b_l,\dots,b_2,b_1]=0$} \quad\hbox{or
\quad $[a_1,a_2\dots,a_k,b_l,\dots,b_2]=0$.}
\notag
\end{equation}
\end{cor}

\begin{pf} $x_1$ and $y_0$ are opposite vectors in Figure~\ref{f!Phi},
so $\Span{x_1,x_2,\dots,y_1,y_0}$ is a half-space with a basic
subdivision. \end{pf}

\subsubsection{Initial cases} \label{ss!init}

We list here all the cases with $r$ or $s\le1$, treating all
cases with attitude not covered by Corollary~\ref{cor!at}.
\[
\begin{tabular}{|c|l|l|l|}
\hline
$\begin{pmatrix}1&0\\0&1\end{pmatrix}$ &
\begin{picture}(55,28)(-37,15)
\put(10,5){\circle*{5}}
\put(10,5){\line(1,0){20}}
\put(10,5){\line(0,1){28}}
\put(-2,0){$0$}
\put(30,5){\circle*{5}}
\put(35,0){$0$}
\put(10,33){\circle*{5}}
\put(-2,30){$0$}
\put(30,33){\circle*{5}}
\put(30,33){\line(-1,0){20}}
\put(30,33){\line(0,-1){28}}
\put(35,30){$0$}
\end{picture} &
$\begin{array}{l} x_0y_1=A, \\ x_1y_0=B. \end{array}$
\\[20pt]
\hline

$\begin{pmatrix}1&0\\b&1\end{pmatrix}$ &
\begin{picture}(55,28)(-37,15)
\put(10,5){\circle*{5}}
\put(10,5){\line(1,0){20}}
\put(10,5){\line(0,1){28}}
\put(-2,0){$0$}
\put(30,5){\circle*{5}}
\put(35,0){$-b$}
\put(10,33){\circle*{5}}
\put(-2,30){$b$}
\put(30,33){\circle*{5}}
\put(30,33){\line(-1,0){20}}
\put(30,33){\line(0,-1){28}}
\put(35,30){$0$}
\end{picture} &
$\begin{array}{l} x_0y_1=x_1^bA, \\ x_1y_0=B. \end{array}$ \\[20pt]
\hline

$\begin{pmatrix}1&a\\0&1\end{pmatrix}$ &
\begin{picture}(55,28)(-37,15)
\put(10,5){\circle*{5}}
\put(10,5){\line(1,0){20}}
\put(10,5){\line(0,1){28}}
\put(-10,0){$-a$}
\put(30,5){\circle*{5}}
\put(35,0){$0$}
\put(10,33){\circle*{5}}
\put(-2,30){$0$}
\put(30,33){\circle*{5}}
\put(30,33){\line(-1,0){20}}
\put(30,33){\line(0,-1){28}}
\put(35,30){$a$}
\end{picture} &
$\begin{array}{l} x_0y_1=A, \\ x_1y_0=y_1^aB. \end{array}$ \\[20pt]
\hline

$\begin{pmatrix}1&1\\s-1&s\end{pmatrix}$ &
\begin{picture}(95,28)(-37,15)
\put(10,5){\circle*5}
\put(10,5){\line(1,0){20}}
\put(10,5){\line(0,1){28}}
\put(-2,0){$0$}
\put(-2,15){$s$}
\put(10,20){\circle*{5}}
\put(30,5){\circle*{5}}
\put(35,0){$-(s-1)$}
\put(10,33){\circle*{5}}
\put(-2,30){$1$}
\put(30,33){\circle*{5}}
\put(30,33){\line(-1,0){20}}
\put(30,33){\line(0,-1){28}}
\put(35,30){$1$}
\end{picture} &
$\begin{array}{l} x_1y_1=x_2A, x_2y_0=y_1B, \\
x_0x_2=x_1^s, \\
x_1y_0=AB, x_0y_1=x_1^{s-1}A. \end{array}$ \\[20pt]
\hline

$\begin{pmatrix}r&1\\r-1&1\end{pmatrix}$ &
\begin{picture}(115,28)(-37,15)
\put(10,5){\circle*{5}}
\put(10,5){\line(1,0){20}}
\put(10,5){\line(0,1){28}}
\put(-39,0){$-(r-1)$}
\put(30,5){\circle*{5}}
\put(35,0){$0$}
\put(10,33){\circle*{5}}
\put(-2,30){$1$}
\put(30,19){\circle*{5}}
\put(35,19){$r$}
\put(30,33){\circle*{5}}
\put(30,33){\line(-1,0){20}}
\put(30,33){\line(0,-1){28}}
\put(35,30){$1$}
\end{picture} &
$\begin{array}{l} x_0y_1=x_2A, x_1y_1=y_1B, \\
y_0y_2=y_1^r, \\
x_1y_0=y_1^{r-1}B, x_0y_1=AB. \end{array}$ \\[20pt]
\hline
\end{tabular}
\]

The cases with $a$ or $b=1$ and $r,s\ge2$ are not exceptional; rather,
they serve as the first regular example of our construction:
\begin{equation}
\begin{array}{cc}
\begin{picture}(120,54)(-20,-8)
\put(-75,15){$\begin{pmatrix}r&rs-1\\1&s\end{pmatrix}$}
\put(5,0){\circle*{5}}
\put(5,0){\line(1,0){25}}
\put(5,0){\line(0,1){38}}
\put(-7,-5){$0$}
\put(-7,15){$r$}
\put(5,20){\circle*{5}}
\put(25,17){$\equiv$}
\put(37,16){$2^{s-1}$}
\put(30,0){\circle*{5}}
\put(37,-5){$-(r-1)$}
\put(5,38){\circle*{5}}
\put(-7,35){$s$}
\put(30,38){\circle*{5}}
\put(30,38){\line(-1,0){25}}
\put(30,38){\line(0,-1){38}}
\put(35,35){$1$}
\end{picture} &
\begin{picture}(120,54)(-100,-8)
\put(-108,15){$\begin{pmatrix}r&1\\rs-1&s\end{pmatrix}$}
\put(5,0){\circle*{5}}
\put(5,0){\line(1,0){25}}
\put(5,0){\line(0,1){38}}
\put(-48,-5){$-(s-1)$}
\put(-24,15){$2^{r-1}$}
\put(0,17){$\equiv$}
\put(30,20){\circle*{5}}
\put(35,16){$s$}
\put(30,0){\circle*{5}}
\put(35,-5){$0$}
\put(5,38){\circle*{5}}
\put(-7,35){$1$}
\put(30,38){\circle*{5}}
\put(30,38){\line(-1,0){25}}
\put(30,38){\line(0,-1){38}}
\put(35,35){$r$}
\end{picture}
\end{array}
\label{eq!k=1}
\end{equation}

\subsection{Construction of $T\subset V_{AB}$ from
$(\protect\begin{smallmatrix} r&a\protect\\b&s
\protect\end{smallmatrix})\in\SL(2,\Z)$} \label{ss!bet}

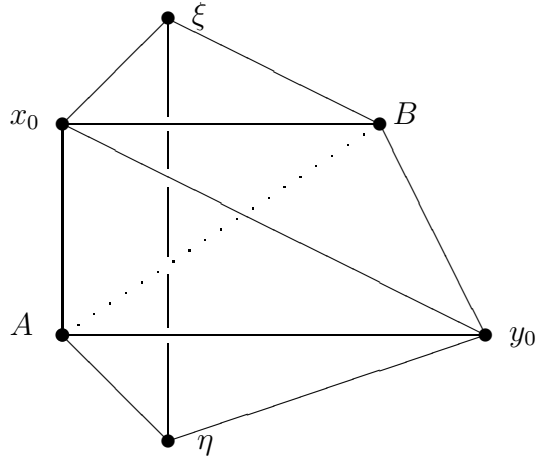
\begin{figure}[ht]
\[
\begin{picture}(190,165)(-20,-42)
\put(0,0){\circle*{5}} \put(-20,0){$A$}
\put(0,0){\line(0,1){80}}
\put(0,0){\line(1,0){160}}
\put(0,0){\line(1,-1){40}}
\put(0,0){\circle*{5}}
\put(0,80){\circle*{5}} \put(-20,80){$x_0$}
\put(0,80){\line(1,0){120}}
\put(0,80){\line(1,1){40}}
\put(0,80){\line(2,-1){160}}
\put(40,120){\circle*{5}} \put(49,118){$\xi$}
\put(40,120){\line(2,-1){80}}
\put(40,-40){\circle*{5}}
\put(120,80){\circle*{5}} \put(125,80){$B$}
\qbezier[20](0,0)(60,40)(120,80)
\put(120,80){\line(1,-2){40}}
\put(160,0){\circle*{5}} \put(169,-2){$y_0$}
\put(40,-40){\line(3,1){120}} \put(51,-43){$\eta$}
\put(40,-40){\line(0,1){37}}
\put(40,3){\line(0,1){20}}
\put(40,30){\line(0,1){26}}
\put(40,63){\line(0,1){14}}
\put(40,83){\line(0,1){37}}
\end{picture}
\]
\caption{The monomial cone $\si_{AB}$\label{f!siAB}}
\end{figure}

To construct the deformation $T\subset V_{AB}$, we pump up the fan
$\Phi(\begin{smallmatrix} r&a\\b&s \end{smallmatrix})$ of \ref{ss!fan}
out of the plane $\bMbar$ to the cone $\si_{AB}$ of Figure~\ref{f!siAB}
in the 4-space of $\M=\Z^4$, using the new variables $A,B$ respectively
to bend along the $\xi$ and $\eta$ axes. In more detail, consider the
monomial lattice $\M\iso\Z^4$ based by $\xi,\eta,A,B$, and the cone
$\si_{AB}$ in $\M_\R$ spanned by
\begin{equation}
\xi,\eta,A,B \quad\hbox{together with}\quad
x_0=(A\eta\1)^r \xi^a, \quad y_0=\eta^b(B\xi\1)^s
\label{eq!x0y0}
\end{equation}
(compare these with the equations of \eqref{eq!x0}).
We draw $\si_{AB}$ projectively, so that it has two quadrilateral faces
$\xi\eta Ax_0$ and $\xi\eta y_0B$ (the ``back''), and four triangles
$\xi Bx_0$, $Bx_0y_0$, $x_0y_0A$ and $y_0A\eta$ (the ``front''). The
primitive vectors orthogonal to these faces are, in order,
\begin{equation}
(0,0,0,1) \enspace (0,0,1,0) \enspace (0,1,1,0) \enspace (rb,rs,1,0)
\enspace (rs,as,0,1) \enspace (1,0,0,1).
\label{eq!orthv}
\end{equation}
We get $\si_{AB}$ from the simplex $\Span{\xi,\eta,x_0,y_0}$ by pulling
out each of the two back faces to quadrilaterals, adding a vertex $A$ in
the plane of $x_0,\xi,\eta$ such that $\xi,\eta,A$ is basic and the
quadrilateral $x_0,\xi,\eta,A$ is convex as shown, and likewise for $B$
in the plane of $y_0,\eta,\xi$. The picture can be viewed from different
perspectives (we use some below, see Figure~\ref{f!pp}), and trying to
read metric properties from these can be misleading.

The dual cone $\si_{AB}^\vee$ is the convex hull of the orthogonal
vectors \eqref{eq!orthv}. Since these are all in the hyperplane of
weights $w$ with $w(AB)=1$, the dotted line from $A$ to $B$ is interior
to $\si_{AB}$, and $AB$ generates the ideal of interior monomials. This
is Danilov's criterion for the toric variety
$V_{AB}=\Spec(\C[\si_{AB}\cap \M])$ to be Gorenstein. The unextended
simplex $\Span{\xi,\eta,x_0,y_0}$ itself does not in general determine
$A,B$ or the matrixes \eqref{eq!mats}.

The 2-faces $\Span{x_0,\xi}$ and $\Span{y_0,\eta}$ of the simplicial
cone $\Span{\xi,\eta,x_0,y_0}$ are also faces of $\si_{AB}$, basic in
$\M$ if and only if $r=1$, respectively $s=1$; they are the monomial
cones of toric surfaces $S_1$ and $S_3$, and are determined exactly
as in \ref{ss!fan}. The new feature is the relations \eqref{eq!x0y0} and
their inverses
\begin{equation}
\eta=(A^rB^ax_0\1)^sy_0^{-a},\quad \xi=x_0^{-b}(A^bB^sy_0\1)^r
\label{eq!xi,eta}
\end{equation}
that determine tag relations at the corners.
Indeed \eqref{eq!x0y0} and \eqref{eq!xi,eta} give
\begin{equation}
\begin{array}{ccc}
(x_0\xi^{-a})^{1/r}=A\eta\1 &\hbox{and}& (y_0\eta^{-b})^{1/s}=B\xi\1 \\[6pt]
(x_0^b\xi)^{1/r}=A^bB^sy_0\1 &\hbox{and}&
(y_0^a\eta)^{1/s}=A^rB^ax_0\1
\end{array} \in M.
\label{eq!4mon}
\end{equation}
and the analogue of Lemma~\ref{l!wT} follows as in \ref{ss!fan}.
\begin{lem}  \label{l!mon} The face
$\Span{x_0,\xi}$ spans a $2$-dimensional vector space in $\M_\R$, that
intersects $\M$ in the sublattice generated as a $\Z$-module by $x_0,\xi$
together with {\em either} of
\begin{equation}
(x_0^b\xi)^{1/r}=A^bB^sy_0\1 \quad\hbox{or}\quad
(x_0\xi^{-a})^{1/r}=A\eta\1.
\notag
\end{equation}
Write $x_0,x_1,\dots,x_{k-1},x_k=\xi$ for the successive monomials along
the Newton boundary of $\Span{x_0,\xi}$. The number $k$ and the
monomials themselves come from either of the continued fraction
expansions
\begin{equation}
[a_0,a_1,\dots,a_{k-1}]=\frac{-b}r \quad\hbox{and}\quad
[a_k,\dots,a_1]=\frac ar
\label{eq!b/s2}
\end{equation}
by the following constructions:
\begin{enumerate}

\item From the bottom, $x_0$ is given, and $x_1=(x_0^\be\xi)^{1/r}$,
where $\be$ is the least residue of $b$ modulo $r$. Thus
$a_0=\rdup{\frac {-b}r}=\frac{-b+\be}r\le0$ and
\begin{equation}
x_1=(x_0^\be\xi)^{1/r}=A^bB^sx_0^{a_0}y_0\1,
\quad\hbox{that is,}\quad
x_1y_0=A^bB^sx_0^{a_0}.
\end{equation}
If $\be=0$ then $r$ divides $b$, whereas $rs-ab=1$ implies that $r,b$
are coprime; thus $r=1$, so that $k=1$ and $x_1=\xi$. Otherwise
$x_2,\dots,x_k$ are determined as usual by tag equations
\begin{equation}
x_{i-1}x_{i+1}=x_i^{a_i} \quad\hbox{for $i=1,\dots,k-1$,}
\notag
\end{equation}
where $[a_1,\dots,a_{k-1}]=\frac r\be$ (see Remark~\ref{rk!JHc}).

\item From the top, $x_k=\xi$ is given; if $r\mid a$ then, as before,
$r=1$ and the only monomials are $x_0,x_1=\xi$. Otherwise, set
$x_{k-1}=(x_0\xi^{r-\al})^{1/r}$, where $\al$ is the least residue of
$a$ mod $r$. Then $r-\al=a_kr-a$ where $a_k=\rdup{\frac ar}\ge1$, and
\begin{equation}
x_{k-1}=\xi^{a_k}(x_0\xi^{-a})^{1/r}=x_k^{a_k}A\eta\1
\quad\hbox{that is,}\quad
x_{k-1}\eta=Ax_k^{a_k}.
\notag
\end{equation}
The remaining monomials are determined by
\begin{equation}
x_{i-1}x_{i+1}=x_i^{a_i}, \quad\hbox{where $[a_{k-1},\dots,a_1]=\frac
r{r-\al}$.}
\notag
\end{equation}
\end{enumerate}

The ring $\C[\M\cap\Span{x_0,\xi}]=\C[S_1]$ is isomorphic to the
invariant ring of the cyclic quotient
singularity $\frac1r(a,1)\iso\frac1r(1,-b)$; here $ab=rs-1$, so that
$ab\equiv-1$ {\rm mod} $s$.

In the same way, $\Span{y_0,\eta}\iso\frac1s(b,1)\iso\frac1s(1,-a)$ with
initial monomials $y_1$ and $y_{l-1}$ determined by the corner tag
equations
\begin{equation}
x_0y_1=A^rB^ay_0^{b_0} \quad\hbox{and}\quad
\xi y_{l-1}=B\eta^{b_l},
\notag
\end{equation}
with $b_0=\rdup{\frac{-a}s}\le0$ and $b_l=\rdup{\frac bs}\ge1$, and the
remaining monomials for $S_3$ are $y_0,y_1,\dots,y_{l-1},y_l$ tagged by
\begin{equation}
[b_0,b_1,\dots,b_{l-1}]=\frac{-a}s \quad\hbox{and}\quad
[b_l,\dots,b_1]=\frac bs\,.
\label{eq!a/r}
\end{equation}
\end{lem}

In conclusion, the following theorem states the complete solution to toric
deformation of tents that smooth the axes at one end.

\begin{theorem} \label{l!be}
Let
\begin{equation}
T=S_0\cup S_1\cup S_2\cup S_3
\notag
\end{equation}
be a tent with two given cyclic quotient singularities in reduced form
$S_1=\frac1r(\al,1)$ and $S_3=\frac1s(\be,1)$. Then toric deformations
$T\subset V_{AB}$ that smooth the $\xi$ and $\eta$ axes correspond
one-to-one with matrixes
\begin{equation}
\begin{pmatrix} r&a\\b&s
\end{pmatrix}
\in\SL(2,\Z)
\quad\hbox{with}\quad
\begin{matrix}
a \equiv \al \mod r, \\
b \equiv \be \mod s.
\end{matrix}
\notag
\end{equation}

Since $ab=rs-1$ obviously implies that $a<r$ or $b<s$, this means that
\begin{equation}
\begin{array}{rl}
\hbox{either} & a=\al\mid rs-1 \quad\hbox{and}\quad b=\frac{rs-1}\al\,,
\\[6pt]
\hbox{or} & b=\be\mid rs-1 \quad\hbox{and}\quad a=\frac{rs-1}\be\,.
\end{array}
\notag
\end{equation}
There may be $0$, $1$ or $2$ solutions.
\end{theorem}

%%%%%%%%%%%%%%%%%%%%%%%%%%%%%%%%%%%%%%%%%%%%%%%%%%%%%%%%%

\section{Classification of diptychs} \label{s!Cl}

A {\em diptych}, for a tent T, is a pair of toric deformations
\[
T\subset V_{AB}\quad\text{and}\quad
T\subset V_{LM}
\]
(the two {\em panels} of the diptych), in which the first smooths the
top axes and the second smooths the bottom axes.

Our construction in \ref{ss!bet} of $T\subset V_{AB}$ is given, already
at the level of $T$, by the fan $\Phi(\begin{smallmatrix} r&a\\b&s
\end{smallmatrix})$ dividing the plane $\bMbar$ into the four cones of
Figure~\ref{f!Phi}. Its key properties are that its four cones give the four
sides of $T$, and the union of its three top cones is {\em one step
beyond convex}; by this we mean that shaving either $x_0$ or $y_0$ off
the two side cones makes the union of the three top cones convex,
which we express by saying that the cone $\left< x_0, y_0\right>$
corresponding to $S_0$ is the {\em big end} of the fan.

\subsection{A second fan $\Phi'(\protect\begin{smallmatrix}
r&g\protect\\h&s \protect\end{smallmatrix})$ and a second panel
$V_{LM}$} \label{ss!2Ph}
For the right panel $V_{LM}$ of our diptych, we need a second fan
$\Phi'$ in a plane $\bMbar'$ (not identified with $\bMbar$), defining
the same tent $T$, but this time the big end of $\Phi'$ is the top
$\Span{\xi,\eta}$ corresponding to $S_2$, and its little end the bottom
$\Span{x_0,y_0}$ corresponding to $S_0$. For this, replace \eqref{eq!x0}
with the base change
\begin{equation}
x_0=\eta^{-r} \xi^{-g}, \quad y_0=\eta^{-h}\xi^{-s}
\quad\hbox{and}\quad
\eta=x_0^{-s}y_0^g,\quad \xi=x_0^hy_0^{-r}
\end{equation}
based on the inverse pair
$(\begin{smallmatrix}-r&-g\\-h&-s\end{smallmatrix})$ and
$(\begin{smallmatrix}-s&g\\h&-r\end{smallmatrix})$, with $g,h\ge0$. As
before, $x_0,\xi,\eta,y_0$ define a fan $\Phi'$ of 4 cones, but with
signs giving the inclusion $-\Span{x_0,y_0}\subseteq\Span{\xi,\eta}$
opposite to \eqref{eq!big}, so that $\Span{\xi,\eta}$ is the big end.

\begin{lem} In $\Phi'$ the cone $\Span{x_0,\xi}$ corresponding is
$\frac1r(1,h)\iso\frac1r(-g,1)$; the cone $\Span{y_0,\eta}$ is
$\frac1s(1,g)\iso\frac1s(-h,1)$.

Hence $\Phi'$ defines the same tent $T$ as $\Phi$ of \ref{ss!fan} if and
only if $-g\equiv\al$ {\rm mod}~$r$ and $-h\equiv\be$ {\rm mod}~$s$.
\end{lem}

We say that $\Phi$ and $\Phi'$ related in this way are {\em partners}.
\ref{ss!Clp} classifies all partner pairs. The analysis of the coordinate
ring of $V_{AB}$ in Lemma~\ref{l!mon} can be applied, with the ends
exchanged, to $V_{LM}$ to prove immediately:

\begin{lem} \label{l!2nd} From $V_{AB}$, the cone
$\Span{x_0,\xi}$ is $\frac1r(a,1)\iso\frac1r(1,-b)$ and from $V_{LM}$ it
is $\frac1r(1,g)\iso\frac1r(-h,1)$. The cone $\Span{y_0,\eta}$ is
$\frac1s(b,1)\iso\frac1s(1,-a)$ and also
$\frac1s(1,h)\iso\frac1s(-g,1)$. Therefore $ag\equiv1\mod r$ and
$bh\equiv1\mod s$; together with $rs-ab=rs-gh=1$, these imply that
\begin{equation}
a+h\equiv b+g\equiv0\mod r\ \hbox{and\/} \mod s.
\label{eq!9}
\end{equation}
\end{lem}
We draw the two monomial cones $\si_{AB}$ and $\si_{LM}$ together in
Figure~\ref{f!pp}; it is easy to see that the union
$\si_{AB}\cup\si_{LM}$ has convex hull a cone with a vertex.

As an example and sanity check, it is a fun exercise to run through
$\left(\begin{smallmatrix} r&a\\b&s
\end{smallmatrix}\right)=\left(\begin{smallmatrix} 7&12\\4&7
\end{smallmatrix}\right)$ and $\left(\begin{smallmatrix} r&g\\h&s
\end{smallmatrix}\right)=\left(\begin{smallmatrix} 7&24\\2&7
\end{smallmatrix}\right)$ to recover the two long
rectangles of Example~\ref{ex!intr}.

%%%%%%%%%%%%%%%%%%%%%%%%%%%%%%%%%%%%%%%%%%%%%%%%%%%%%%%%%

\subsection{Classification of partner pairs} \label{ss!Clp}
Classifying all partner pairs $\Phi$, $\Phi'$ of fans is an elementary
``infinite descent''.

\paragraph{Rules of the game:} Given integers
\begin{equation}
\begin{gathered}
r,s\ge1,\quad a,b,g,h\ge0,\quad\hbox{with}\quad ab=gh=rs-1
\\
\hbox{and}\quad
a+h\equiv b+g\equiv0\mod r\ \hbox{and\/} \mod s.
\end{gathered}
\label{eq!10}
\end{equation}
Use the congruences to define two integers $d\ge1$ and $e\ge1$:
\begin{equation}
\label{eq!10a}
a+h=d s \quad\hbox{and}\quad b+g=e r.
\end{equation}

\begin{theorem}[Classification Theorem I] \label{th!d-e}
Each solution of (\ref{eq!10}--\ref{eq!10a}) is one of the exceptional
solutions \eqref{eq!ex} below, or is given either by
\begin{equation}
\begin{aligned}
\begin{pmatrix} r & a\\ b & s \end{pmatrix} &=
\begin{pmatrix} d & -1\\ 1 & 0 \end{pmatrix}
\begin{pmatrix} e & -1\\ 1 & 0 \end{pmatrix}\cdots
\begin{pmatrix} e\hbox{ or }d & -1\\ 1 & 0 \end{pmatrix}
\begin{pmatrix} 0 & 1\\ -1 & 0 \end{pmatrix},
\\[6pt]
\begin{pmatrix} r & g\\ h & s \end{pmatrix}
&=
\begin{pmatrix} 0 & 1\\ -1 & 0 \end{pmatrix}
\begin{pmatrix} 0 & -1\\ 1 & d\hbox{ or }e \end{pmatrix}
\cdots\begin{pmatrix} 0 & -1\\ 1 & d \end{pmatrix}
\begin{pmatrix} 0 & -1\\ 1 & e \end{pmatrix}
\end{aligned}
\label{eq!d-e1}
\end{equation}
or the same with the two lefthand sides exchanged, or by
\begin{equation}
\begin{aligned}
\begin{pmatrix} r & a\\ b & s \end{pmatrix} &=
\begin{pmatrix} 0 & 1\\ -1 & d \end{pmatrix}
\begin{pmatrix} 0 & 1\\ -1 & e \end{pmatrix}\cdots
\begin{pmatrix} 0 & 1\\ -1 & e\hbox{ or }d \end{pmatrix}
\begin{pmatrix} 0 & -1\\ 1 & 0 \end{pmatrix},
\\[6pt]
\begin{pmatrix} r & g\\ h & s \end{pmatrix}
&=
\begin{pmatrix} 0 & -1\\ 1 & 0 \end{pmatrix}
\begin{pmatrix} d\hbox{ or }e & 1\\ -1 & 0 \end{pmatrix}
\cdots\begin{pmatrix} d & 1\\ -1 & 0 \end{pmatrix}
\begin{pmatrix} e & 1\\ -1 & 0 \end{pmatrix}
\label{eq!d-e}
\end{aligned}
\end{equation}
or the same with the two lefthand sides exchanged.

In each case, the values $d,e\ge1$
alternate, the two lines have the same number $k+1$ of factors for some
$k\ge1$, and the values of $d$, $e$ and $k$ that are allowed are
constrained only by the following table:
\begin{equation}
\label{eq!lt6}
\begin{array}{c||c|c|c|c|c}
d e & 0 & 1 & 2 & 3 & \ge4 \\
\hline
k & 1 & \le2 & \le3 & \le5 & \mathrm{any}
\end{array}
\end{equation}

\paragraph{Exceptional solutions} The cases $b=g=0$ or $a=h=0$, the
matrixes
\begin{equation}
\label{eq!ex}
\begin{array}{ccccc}
\begin{pmatrix} r & a\\ b & s \end{pmatrix} &=
\begin{pmatrix} 1 & a\\ 0 & 1 \end{pmatrix}
&\quad
\text{and}
\quad&
\begin{pmatrix} r & g\\ h & s \end{pmatrix} &=
\begin{pmatrix} 1 & 0\\ h & 1 \end{pmatrix},
\end{array}
\end{equation}
for any $a,g\ge0$, or the same with both matrixes transposed.
\end{theorem}

\begin{remark}
(1) In the statement, exchanging the two lefthand sides amounts
to exchanging the roles of the two long rectangles, so exchanges
$V_{AB}$ and $V_{LM}$ in the diptych (and turns them upside
down if one draws them as long rectangles).
Whether the first or second factorisation occurs depends on the
attitude of the long rectangles, which is determined by whether
$b<r$ or $b>r$; this becomes clear in the proof.

(2) The computation of a pair of long rectangles from these two matrices
is implicit from Lemma~\ref{l!wT}, but we spell it out.
The tags on the long rectangle of $V_{AB}$ are given by the
tags of the continued fraction expansion
\[
-b/r = [a_0,\dots,a_{k-1}]
\quad\text{and}\quad
a/r = [a_k,\dots,a_1].
\]
If $b<r$ and $d,e\ge2$, then the alternating $d,e$ tags run up
the lefthand side, and the first of these will be of the form $[0,d,e,d,\dots]$.
If either $d=1$ or $e=1$, the tags one computes are those after
blowdown of the $1$s, as in Proposition~\ref{prop!hj}(b); one can
reintroduce them by blowup as redundant generators to see the
alternating $d,e$ sequence.
The tags down the righthand side are
\[
-a/s = [b_0,\dots,b_{l-1}]
\quad\text{and}\quad
b/s = [b_l,\dots,b_1].
\]
The tags on the long rectangle for $V_{LM}$ are
\[
g/r = [a_0',a_1,\dots,a_{k-1}],
\qquad
-h/r = [a_k',a_{k-1},\dots,a_1]
\]
and
\[
h/s = [b_0',b_1,\dots,b_{l-1}]
\qquad
-g/s = [b_l',b_{l-1},\dots,b_1]
\]
where all but the corner tags are of course common to both
long rectangles.

(3) The exceptional cases correspond to the not-very-long
rectangles and not-very-surprising diptych varieties:
\begin{equation}
\begin{picture}(55,28)(-37,15)
\put(10,5){\circle*{5}}
\put(10,5){\line(1,0){20}}
\put(10,5){\line(0,1){28}}
\put(-2,0){$0$}
\put(30,5){\circle*{5}}
\put(35,0){$-a$}
\put(10,33){\circle*{5}}
\put(-2,30){$a$}
\put(30,33){\circle*{5}}
\put(30,33){\line(-1,0){20}}
\put(30,33){\line(0,-1){28}}
\put(35,30){$0$}
\end{picture}
\qquad
\begin{picture}(55,28)(-37,15)
\put(10,5){\circle*{5}}
\put(10,5){\line(1,0){20}}
\put(10,5){\line(0,1){28}}
\put(-2,0){$0$}
\put(30,5){\circle*{5}}
\put(35,0){$h$}
\put(10,33){\circle*{5}}
\put(-12,30){$-h$}
\put(30,33){\circle*{5}}
\put(30,33){\line(-1,0){20}}
\put(30,33){\line(0,-1){28}}
\put(35,30){$0$}
\end{picture}
\qquad
\qquad
\begin{aligned}
x_0y_1 &= Ax_1^a + My_0^h \\
x_1y_0 &= B + L
\end{aligned}
\notag
\end{equation}
and we do not mention them again.

(4) The cases $b=h=0$ or $a=g=0$ are regular solutions in
Theorem~\ref{th!d-e} with $k=1$ and (say)
$d=a$, $e=g$:
\begin{equation}
\label{eq!in}
\begin{pmatrix} 1 & a\\ 0 & 1 \end{pmatrix} =
\begin{pmatrix} d & -1\\ 1 & 0 \end{pmatrix}
\begin{pmatrix} 0 & 1\\ -1 & 0 \end{pmatrix},
\quad
\begin{pmatrix} 1 & g\\ 0 & 1 \end{pmatrix} =
\begin{pmatrix} 0 & 1\\ -1 & 0 \end{pmatrix}
\begin{pmatrix} 0 & -1\\ 1 & e \end{pmatrix}.
\end{equation}
They provide the endpoint of our infinite descent:
\begin{equation}
\begin{picture}(55,28)(-37,15)
\put(10,5){\circle*{5}}
\put(10,5){\line(1,0){20}}
\put(10,5){\line(0,1){28}}
\put(-2,0){$0$}
\put(30,5){\circle*{5}}
\put(35,0){$-a$}
\put(10,33){\circle*{5}}
\put(-2,30){$a$}
\put(30,33){\circle*{5}}
\put(30,33){\line(-1,0){20}}
\put(30,33){\line(0,-1){28}}
\put(35,30){$0$}
\end{picture}
\qquad
\begin{picture}(55,28)(-37,15)
\put(10,5){\circle*{5}}
\put(10,5){\line(1,0){20}}
\put(10,5){\line(0,1){28}}
\put(-2,0){$g$}
\put(30,5){\circle*{5}}
\put(35,0){$0$}
\put(10,33){\circle*{5}}
\put(-2,30){$0$}
\put(30,33){\circle*{5}}
\put(30,33){\line(-1,0){20}}
\put(30,33){\line(0,-1){28}}
\put(35,30){$-g$}
\end{picture}
\qquad
\qquad
\begin{aligned}
x_0y_1 &= Ax_1^a + M \\
x_1y_0 &= B + Lx_0^g
\end{aligned}
\notag
\end{equation}
The equations can be used to eliminate variables $B$ and $M$,
so the diptych varieties in these cases are simply isomorphic to $\C^6$.

(5) The restriction on $k$ when $de\le3$ in \eqref{eq!lt6} arises
because the product in \eqref{eq!d-e} no longer satisfies $r,s,a,b\ge0$
for bigger values of $k$. Thus
\begin{equation}
\begin{pmatrix} d & -1\\ 1 & 0 \end{pmatrix}
\begin{pmatrix} e & -1\\ 1 & 0 \end{pmatrix}
\begin{pmatrix} 0 & 1\\ -1 & 0 \end{pmatrix}
=\begin{pmatrix} d & de-1\\ 1 & e \end{pmatrix}
\notag
\end{equation}
has top righthand entry $<0$ for $de=0$ and $k=2$. For
$de=1,2,3$ and $k=3,4,6$ respectively, the product of
$k$ factors is $-1$:
\begin{equation}
\begin{pmatrix} d & -1\\ 1 & 0 \end{pmatrix}
\begin{pmatrix} e & -1\\ 1 & 0 \end{pmatrix}
\cdots
{\renewcommand{\arraycolsep}{.6em}
\begin{pmatrix} \hbox{$d$ or $e$} & -1\\ 1 & 0 \end{pmatrix}
}
=\begin{pmatrix} -1 & 0\\ 0 & -1 \end{pmatrix},
\notag
\end{equation}
\end{remark}
so we are basically into elements of finite order in $\SL(2,\Z)$.

\paragraph{Proof of the Classification Theorem}

The following two operations preserve all the equalities and congruences
in the rules of the game
while interchanging the roles of $e$ and~$d$:
\begin{equation}
\begin{array}{ccccc}
\begin{pmatrix} r & a\\ b & s \end{pmatrix}
&\mapsto&
\begin{pmatrix} 0 & 1\\ -1 & d \end{pmatrix}
\begin{pmatrix} r & a\\ b & s \end{pmatrix} &=&
\begin{pmatrix} b & s\\ d b-r & h \end{pmatrix}
\\ [10pt]
\begin{pmatrix} r & g\\ h & s \end{pmatrix}
&\mapsto&
\begin{pmatrix} r & g\\ h & s \end{pmatrix}
\begin{pmatrix} e & 1\\ -1 & 0 \end{pmatrix}
&=&
\begin{pmatrix} b & r\\ e h-s & h \end{pmatrix}
\end{array}
\label{eq!22}
\end{equation}
and (``its inverse with $d$, $e$ interchanged'')
\begin{equation}
\begin{array}{ccccc}
\begin{pmatrix} r & a\\ b & s \end{pmatrix}
&\mapsto&
\begin{pmatrix} e & -1\\ 1 & 0 \end{pmatrix}
\begin{pmatrix} r & a\\ b & s \end{pmatrix} &=&
\begin{pmatrix} g & e a - s\\ r & a \end{pmatrix}
\\ [10pt]
\begin{pmatrix} r & g\\ h & s \end{pmatrix}
&\mapsto&
\begin{pmatrix} r & g\\ h & s \end{pmatrix}
\begin{pmatrix} 0 & -1\\ 1 & d \end{pmatrix}
&=&
\begin{pmatrix} g & d g-r\\ s & a \end{pmatrix}
\end{array}
\label{eq!23}
\end{equation}
Indeed, under operation \eqref{eq!22} transforms the
equalities for the sums of opposing off-diagonal terms
$a+h=d s$ and $b+g=e r$ into
\[
s + (e h-s) = e h
\quad\textrm{and}\quad
(d b-r) + r = d b.
\]
The inequalities in the rules of the game need not be
preserved, but their failure is a termination condition.

It turns out that a series of these operations (say using \eqref{eq!22}
to result in \eqref{eq!d-e1}) with alternating $e$,
$d$ reduces to the {\em initial case}
$(\begin{smallmatrix}1&e\\0&1\end{smallmatrix})$,
$(\begin{smallmatrix}1&d\\0&1\end{smallmatrix})$
(or the other way round) and then down to
$(\begin{smallmatrix}0&1\\-1&0\end{smallmatrix})$,
$(\begin{smallmatrix}0&1\\-1&0\end{smallmatrix})$, so that
inverting the procedure proves the theorem.
The only point is to show that these operations, or combinations of them,
decrease the entries of both matrixes; the claim then follows.
When $d,e\ge2$, which operation works is a matter of the attitude
of the long rectangles; when $d$ or $e=1$, either operation
decreases some entries and increases others, but composing the two,
in an order determined by attitude, decreases them all.
We treat the attitude in terms of the relative sizes of $r,\dots,h$.

Consider an initial pair
\[
\begin{pmatrix} r&a\\b&s \end{pmatrix}
\quad\hbox{and}\quad
\begin{pmatrix} r&g\\h&s \end{pmatrix}
\]
satisfying the rules of the game.

\paragraph{The case $d\ge 2$ and $e\ge2$}

Suppose provisionally that $b<r$.
We apply the reduction operation \eqref{eq!22} to get
\[
\begin{pmatrix} r&a\\b&s \end{pmatrix}
\mapsto
\begin{pmatrix} b&s\\db-r&h \end{pmatrix}
\quad\hbox{and}\quad
\begin{pmatrix} r&g\\h&s \end{pmatrix}
\mapsto
\begin{pmatrix} b&r\\eh-s&h \end{pmatrix}.
\]
We claim that every entry of the two resulting matrices is strictly smaller
than the corresponding entry of the initial pair: this holds in the
top left entry of either matrix by the case assumption.

Since $rs-ab=1$, we get $s<a$. 
By \eqref{eq!10a}, $b<r$ implies that $g>r$, and again it is immediate that $s>h$.
It remains to consider the two larger entries in the bottom left of the pair.

To see that $eh-s<h$, it is enough to check $hb-1<hr$:
indeed multiplying by $r$ and substituting for $er=b+g$ and $rs=gh-1$ gives
\[
ehr - rs = h(b+g) - (gh+1) = bh - 1 < hr.
\]
But the inequality $bh - 1 < rh$ holds by the initial assumption.

Similarly we check $db-r<b$ by observing that the equivalent inequality
\[
bh-1=b(a+h)-(ab+1)=bds-rs<bs,
\]
holds since we already know that $h<s$.

The inequality $db-r<b$ also implies that the
resulting matrices have the same attitude, so that if $b,h\ge 1$,
the same operation \eqref{eq!22} will be applied at the next step,
but with $e,d$ exchanged, and the descent continues.

The termination condition is that $r=0$ or $s=0$, since the inequalities
for $r,s$ are the only rules of the game that the reduction operation can break.
In either case $ab=-1$, so that $b<r$ and its friends imply
$(\begin{smallmatrix}r&a\\b&s\end{smallmatrix})=(\begin{smallmatrix}0&1\\-1&0\end{smallmatrix})$
and
$(\begin{smallmatrix}r&g\\h&s\end{smallmatrix})=(\begin{smallmatrix}0&1\\-1&0\end{smallmatrix})$.
Multiplying by the inverse matrices gives the factorisation \eqref{eq!d-e1}.

Finally, notice that if instead we have $b>r$, then we must have $g<r$ (otherwise
both $a<s$ and $h<s$, implying $d=1$, contrary to the case assumptions),
in which case the operation \eqref{eq!23} performs the required reduction.
This gives the factorisation \eqref{eq!d-e}.

\paragraph{The case $d>4$ and $e=1$}

The definitions \eqref{eq!10a} imply that $b<r$ and $g<r$.
Suppose provisionally that $b<g$.

In this case we apply the reduction operation \eqref{eq!22} twice, alternating
$d$ and $e$, to see that it reduces the pair. Thus we compute a new pair
\[
\begin{pmatrix} 0&1\\-1&1 \end{pmatrix}
\begin{pmatrix} 0&1\\-1&d \end{pmatrix}
\begin{pmatrix} r&a\\b&s \end{pmatrix}
=
\begin{pmatrix} db-r & h \\ (d-1)b-r & h-s \end{pmatrix}
\]
and
\[
\begin{pmatrix} r&g\\h&s \end{pmatrix}
\begin{pmatrix} 1&1\\-1&0 \end{pmatrix}
\begin{pmatrix} d&1\\-1&0 \end{pmatrix}
=
\begin{pmatrix} db-r & b \\ (d-2)h-a & h-s \end{pmatrix}.
\]

We start knowing $b, g<r$, and so $a,h>s$, together with the case assumption $b<g$,
or equivalently $h<a$.
So we have $bh-1 < gh-1 = rs$. Substituting for $h$ from \eqref{eq!10a} gives
$dbs-rs=dbs-ab-1<rs$, so $db-r<r$.
Similarly $bh-1<gh-1=rs$, so substituting for $b$ from \eqref{eq!10a} gives
$hr-rs=hr-gh-1<rs$, so $h-s<s$.

The two longer inequalities remain: $(d-1)b-r<b$ and $(d-2)h-a<h$.
For the first, note that $hb-1<2bs$ since $h-s<s$. Substituting for $h$ gives
$(d-1)bs-rs<(d-1)bs-ab-1=(ds-a)b-bs-1<bs$, and dividing by $s$ concludes.

Substituting for $h$ in $h-s<s$ gives $a>(d-2)s$.
Since $rs-ab=1$, $r/b=1/(bs)+a/s>1/(bs)+d-2>d-2$ and we have $r>(d-2)b$.
Substituting for $r$ now gives $g>(d-2)b$.
Since $g/b=a/h$, we get $(d-2)h-a<h$ as required for the second longer inequality.

The same calculations show that $d b-r\ge 0$, so that the analogue of the provisional
supposition $b<g$ holds again after the two reduction steps, and the descent continues
unless we have reached a terminal stage where the inequalities don't hold any more.
(Since we jumped straight in with two reduction steps,
we should also check whether the inequalities already fail
after just one of the steps: by the same calculation, this
would only happen if $b=h=0$, in which case the theorem follows despite
the fact that not all matrix entries reduce.)

Finally, if $b>g$ then operation \eqref{eq!23} applied twice makes the reduction
following a similar analysis (in this case a terminal state cannot arise after just one of the steps).

\paragraph{The case $d=1$, $e>4$}

The definitions \eqref{eq!10a} imply that $a<s$ and $h<s$.
Suppose provisionally that $h<a$.

In this case we apply the reduction operation \eqref{eq!22} twice, alternating
$d$ and $e$, to see that it reduces the pair. Thus we compute a new pair
\[
\begin{pmatrix} 0&1\\-1&e \end{pmatrix}
\begin{pmatrix} 0&1\\-1&1 \end{pmatrix}
\begin{pmatrix} r&a\\b&s \end{pmatrix}
=
\begin{pmatrix} db-r & h \\ (d-1)b-r & h-s \end{pmatrix}
\]
and
\[
\begin{pmatrix} r&g\\h&s \end{pmatrix}
\begin{pmatrix} e&1\\-1&0 \end{pmatrix}
\begin{pmatrix} 1&1\\-1&0 \end{pmatrix}
=
\begin{pmatrix} db-r & b \\ (d-2)h-a & h-s \end{pmatrix}.
\]
The analysis is now virtually identical to the other cases, and we omit it.

%%%%%%%%%%%%%%%%%%%%%%%%%%%%%%%%%%%%%%%%%%%%%%%%%%%%%%%%%

\section{Combining monomial cones $\si_{AB}$ and $\si_{LM}$}
\label{s!pp}

Here we spell out how the factorisations in the Classification
Theorem~\ref{th!d-e} imply growth conditions and congruences on the
generators of the varieties $V_{AB}$ and $V_{LM}$; these are the
conditions (i)--(v) of Corollaries~\ref{c!ming} and~\ref{c!Q}. {\bf From
Corollary~\ref{c!ming}(ii) onwards we restrict to the case $d,e\ge 2$}.
In Section~\ref{s!pf} we also impose $de>4$, so that we are in the main
case of the introduction~\ref{intro!I}. The other cases are treated in
\cite{BR2,BR3}.

\subsection{The Pretty Polytope $\Pi(d,e,k)$} \label{ss!pp}
\begin{figure}[ht]
\begin{picture}(200,160)(-100,20)
\put(0,104){$L$}
\put(169,104){$M$}
\put(10,100){\line(1,1){80}}
\put(10,100){\line(1,0){71}}
\put(103,100){\line(1,0){67}}
\put(10,100){\line(1,-1){80}}
\put(90,20){\line(1,1){80}}
\put(97,18){$B$}
\put(90,180){\line(1,-1){80}}
\put(98,178){$A$}
\put(90,180){\line(0,-1){160}}
\qbezier(90,180)(98,120)(106,60)
\qbezier(90,180)(85,120)(80,60)
\put(80,60){\line(1,0){26}}
\qbezier[45](10,100)(77,105)(144,110)
\qbezier[45](10,100)(77,93)(144,86)
\put(144,86){\line(0,1){24}}
\qbezier[20](106,60)(117,86)(144,86)
\qbezier[53](80,60)(100,95)(144,110)
\put(108,55){$x_0$}
\put(134,78){$x_k$}
\put(68,55){$y_0$}
\put(137,118){$y_l$}
\end{picture}
\caption{Pretty Polytope $\Pi$: Starting from simplex $ABLM$, pull out
$x_0$ on plane $ABL$, etc., with crosspiece $x_0y_0$ on the edge $AB$ in
ratio $1:d$, and $x_ky_l$ on the edge $LM$ in ratio $1:e$. $\Pi$ has 8
vertices and 12 triangular faces; $A,B,L,M$ have valency 5, and
$x_0,y_0,x_k,y_l$ valency 4.} \label{f!pp}
\end{figure}
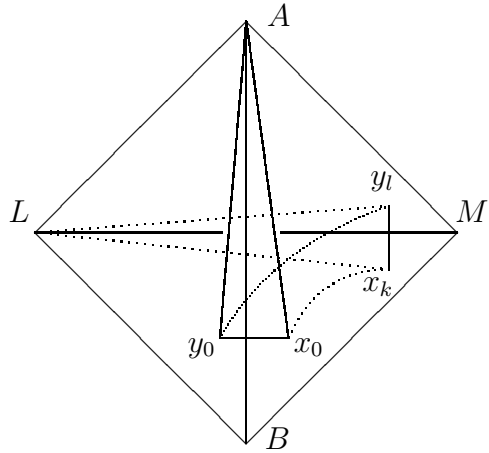

All our varieties $T,V_{AB},V_{LM},V_{ABLM}$ are equivariant under the
same torus $\T=\Gm^4$; write $\M=\Hom(\T,\Gm)$ for its character
lattice, identified with the monomial lattice of both $V_{AB}$ and
$V_{LM}$. The coordinate ring of $V_{ABLM}$ constructed in
Section~\ref{s!pf} is $\M$-graded (that is, $\T$-equivariant). Write
$f\simT g$ to mean that $f$ and $g$ are eigenfunctions with the same
$\T$-weight or eigenvalue in $\M$. This chapter mostly treats the
$\T$-weights of monomials; we mix additive and multiplicative notation,
and sometimes write $=$ for $f\simT g$, so that, for example, the first
equation of \eqref{eq!-1/d} means $x_0\simT L^{-1/d}A^\ga B^\de$.

The {\em Pretty Polytope} $\Pi$ of Figure~\ref{f!pp} combines the two
polytopes $\si_{AB}$ of \ref{ss!bet} and $\si_{LM}$ of \ref{ss!2Ph}.
While $V_{AB}$ and $V_{LM}$ each provided many possible $\Z$-bases of
$\M$, we use instead the {\em impartial}\/ $\Q$-basis $L,M,A,B$, writing
out the $\T$-weights of $x_{0\dots k}, y_{0\dots l}$ as follows:
\begin{equation}
x_0 = (-\dth,0,\ga,\de)
\quad\hbox{and}\quad
x_1 =(0,\eeth,\al,\be),
\label{eq!-1/d}
\end{equation}
where
\begin{equation}
\begin{pmatrix} \al & \be \\ \ga & \de \end{pmatrix}=
\begin{cases}
\begin{pmatrix} 0 & 1\\ -1 & e \end{pmatrix}\cdots
\begin{pmatrix} 0 & 1\\ -1 & e \end{pmatrix}
\begin{pmatrix} -\dth & 0 \\ 0 & \eeth \end{pmatrix}
&\hbox{if $k$ is even} \\[16pt]
\begin{pmatrix} 0 & 1\\ -1 & e \end{pmatrix}\cdots
\begin{pmatrix} 0 & 1\\ -1 & d \end{pmatrix}
\begin{pmatrix} -\eeth & 0 \\ 0 & \dth \end{pmatrix}
&\hbox{if $k$ is odd}
\end{cases}
\label{eq!maca}
\end{equation}
($k$ factors in each product). Compared to \eqref{eq!d-e}, we simply
remove the first and last tags ($d$ at $x_k$ and 0 at $x_0$), and put in
denominators $d,e$ corresponding to the index of the sublattice
$\M'=\Z\cdot(L,M,A,B)\subset\M$ (see Corollary~\ref{c!Q}).

The impartial basis gives $\M$ two projections
\begin{equation}
\pi_{AB}\colon\M \to \Q^2 \quad\hbox{and}\quad
\pi_{LM}\colon\M \to \Q^2
\label{eq!pi}
\end{equation}
that track the exponents of $A,B$ and of $L,M$. The image group $\Q^2$
is partially ordered, and we write $\pi_{LM}(m)\le0$ to mean that $m\in\M$ has
nonpositive $L,M$ exponents, etc.

\begin{prop} \label{p!xys}  In the impartial basis
$L,M,A,B$, the monomials $x_0,\dots,y_l$ have $\T$-weights of the form
(for even $k$):
\begin{equation}
\renewcommand{\arraystretch}{1.2}
\renewcommand{\arraycolsep}{.2em}
\begin{array}{rclccccc}
x_0 &=& (&-\dth&0&\ga&\de&) \\
x_1 &=& (&0&\eeth&\al&\be&) \\
x_2 &=& (&\dth&1&\lank&\lank&) \\
x_3 &=& (&1&d-\eeth&\lank&\lank&) \\
\end{array} \kern1.25cm
\begin{array}{rclccccc}
 &\cdots \\
x_{k-2} &=& (&\lank&\lank&\dth&1&) \\
x_{k-1} &=& (&\al&\be&0&\eeth&) \\
x_k &=& (&\ga&\de&-\dth&0&)
\end{array}
\label{eq!xs}
\end{equation}
and
\begin{equation}
\renewcommand{\arraystretch}{1.2}
\renewcommand{\arraycolsep}{.35em}
\begin{array}{rccclccccc}
y_0 &&=&& (&0&-\eeth&d\ga-\al&d\de-\be&)\\
y_1 &&=&& (&\dth&1-\eeth& \lank & \lank &)\\ % (&d+1&)\ga-\al&(&d+1&)\de-\be&)\\
 &&\cdots \\
y_{j+1} &&=&& b_jy_j - y_{j-1} \kern-2cm \\
 &&\cdots \\
y_{l-1} &&=&& (& \lank & \lank &\dth&1-\eeth&) \\ % (&d+1&)\ga-\al&(&d+1&)\de-\be
y_l &&=&& (&d\ga-\al&d\de-\be&0&-\eeth&)
\end{array}
\label{eq!ys}
\end{equation}
where the $b_j$ in \eqref{eq!ys} are the tags at $y_j$ (usually $2$ or $3$).

When $k$ is odd, the top-to-bottom symmetry swaps $d$ and $e$. At the
top, nothing changes (recall that we define $\al,\be,\ga,\de$ in
$x_1,x_0$ by the other choice in \eqref{eq!maca}); at the bottom we do
$d\bij e$ and modify $\al,\be,\ga,\de$ accordingly, giving
$x_k=(\ga',\de',-\eeth,0)$ and $y_l=(e\ga'-\al',e\de'-\be',0,-\dth)$.
\end{prop}

\begin{pf} The matrix product in \eqref{eq!maca} ensures that the
$k-1$ changes of basis of the form $x_2=x_1^ex_0\1$, etc., take the last
two entries $\left(\begin{smallmatrix} \al & \be \\ \ga & \de
\end{smallmatrix}\right)$ of $x_1,x_0$ into the last two entries
$\left(\begin{smallmatrix} -1/d & 0 \\ 0 & 1/e \end{smallmatrix}\right)$
of $x_k,x_{k-1}$. The first two columns then just record known data from
$V_{LM}$, and the last two from $V_{AB}$. \end{pf}

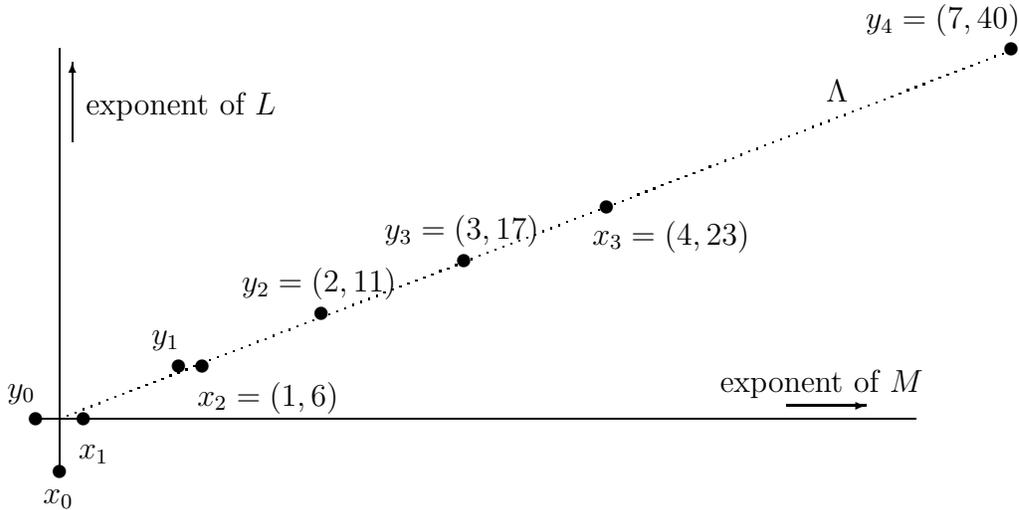
\begin{figure}[bht]
\begin{picture}(360,180)(-21,-28)
\put(-9,0){\line(1,0){333}}
\put(275,5){\vector(1,0){30}}
\put(250,10){exponent of $M$}
\put(0,-20){\line(0,1){160}}
\put(5,105){\vector(0,1){30}}
\put(10,115){exponent of $L$}
\put(0,-20){\circle*{5}} % x0
\put(-6.1,-32){$x_0$}
\put(-9,0){\circle*{5}} % y0
\put(-19.8,8){$y_0$}
\put(9,0){\circle*{5}} % x1
\put(7.2,-15){$x_1$}
\put(54,20){\circle*{5}} % x2
\put(52.2,5){$x_2=(1,6)$}
\put(45,20){\circle*{5}} % y1
\put(35,28){$y_1$}
\put(99,40){\circle*{5}} % y2
\put(69,48){$y_2=(2,11)$}
\put(153,60){\circle*{5}} % y3
\put(123,68){$y_3=(3,17)$}
\put(207,80){\circle*{5}} % x3
\put(201.6,66){$x_3=(4,23)$}
\put(360,140){\circle*{5}} % y4
\put(305,148){$y_4=(7,40)$}
\qbezier[120](0,0)(180,69.7)(360,139.4)
\put(290,120){$\La$}
\end{picture}
\caption{Scissors (compare the dots of Figure~\ref{f!Phi}.b). The
exponents of $L$ are in units of $1/4$ and those of $M$ in units of
$1/6$. The initial points are $x_0=(-1,0)$, $y_0=(0,-1)$, $x_1=(0,1)$,
$y_1=(1,5)$. \label{f!sciss}}
\end{figure}

\begin{cor} \label{c!ming}
\begin{enumerate}
\renewcommand{\labelenumi}{(\roman{enumi})}
\item Except for the explicit $-\dth$ and $-\eeth$ in $x_0,x_k,y_0,y_l$
at the four corners, all the entries are $\ge0$.

\item (From here on, we assume $d,e\ge2$.)
The $L$ and $M$ exponents $\pi_{LM}(x_i)$ and
$\pi_{LM}(y_j)$ increase monotonically with $i$ and $j$ (in fact,
increase {\em exponentially} if $de>4$, as illustrated in
Figure~\ref{f!sciss}), while $\pi_{AB}(x_i)$ and $\pi_{AB}(y_j)$
decrease.

\item No $x_{0\dots k}$ or $y_{0\dots l}$ is $\T$-equivalent to a
monomial in the other variables (all the $x_i$, $y_j$, $A,B,L,M$).
\end{enumerate}

\end{cor}

For (iii), notice that the $x_i$, $y_j$, $A$ and $B$ are minimal
generators of the coordinate ring of $V_{AB}$ by the results of
\ref{ss!bet}. So it is impossible to write even the first two entries of
$x_i$ or $y_j$ as a positive integral combination of the other variables.

\begin{exa}[Case $k=2$] \label{x!k=2}
Then
\begin{equation}
\begin{pmatrix} \al & \be \\ \ga & \de \end{pmatrix}=
\begin{pmatrix} 0 & 1 \\ -1 & e \end{pmatrix}
\begin{pmatrix} -\dth & 0 \\ 0 & \eeth \end{pmatrix}
=\begin{pmatrix} 0 & \eeth \\ \dth & 1 \end{pmatrix}
\notag
\end{equation}
The variables $x_{0\dots2}, y_{0\dots d}$ are
\begin{equation}
\renewcommand{\arraystretch}{1.2}
\renewcommand{\arraycolsep}{0.25em}
\begin{array}{rcl}
x_0 &=& (-\dth,0,\dth,1) \\
x_1 &=& (0,\eeth,0,\eeth) \\
x_2 &=& (\dth,1,-\dth,0)
\end{array}
\kern1.5cm
\begin{array}{rcl}
y_0 &=& (0,-\eeth,1,d-\eeth) \\
y_i &=& (\textstyle{\frac id },i-\eeth,1-\textstyle{\frac id },d-i-\eeth) \\
&& \quad\hbox{for $i=0,\dots,d$} \\
y_{d} &=& (1,d-\eeth,0,-\eeth)
\end{array}
\notag
\end{equation}
Check top-to-bottom symmetry. Check the two tag equations at $x_0$:
\begin{equation}
dx_0 + (1,0,0,0) = x_1+y_0; \enspace\hbox{and}\enspace
0x_0 + (0,0,1,d) = (0,0,1,d)
\notag
\end{equation}
corresponding to the corner tag equations $x_1y_0=x_0^dL$ in $V_{LM}$
and $x_1y_0=AB^d$ in $V_{AB}$. Check the tag equations at $y_0$:
$1y_0 + (0,1,0,0) = x_0+y_1$, and
\begin{equation}
(e-1)x_1 + (0,0,1,d-1) = (0,1-\eeth,1,d-\eeth)
\notag
\end{equation}
corresponding to $x_0y_1=y_0M$ in $V_{LM}$ and $x_0y_1=x_1^{e-1}AB^{d-1}$ in $V_{AB}$.
\end{exa}

\begin{exa}[Case $k=3$] \label{x!k=3}
Then
\begin{equation}
\begin{pmatrix} \al & \be \\ \ga & \de \end{pmatrix}=
\begin{pmatrix} 0 & 1 \\ -1 & e \end{pmatrix}
\begin{pmatrix} 0 & 1 \\ -1 & d \end{pmatrix}
\begin{pmatrix} -\eeth & 0 \\ 0 & \dth \end{pmatrix}
=\begin{pmatrix} \eeth & 1 \\ 1 & e-\dth \end{pmatrix}
\notag
\end{equation}
So $x_{0\dots3}, y_{0\dots d+e-2}$ are
\begin{equation}
\renewcommand{\arraystretch}{1.2}
\renewcommand{\arraycolsep}{0.2em}
\begin{array}{rcl}
x_0 &=& (-\dth,0,1,e-\dth) \\
x_1 &=& (0,\eeth,\eeth,1) \\
x_2 &=& (\dth,1,0,\dth) \\
x_3 &=& (1,d-\eeth,-\eeth,0)
\end{array}
\kern1cm
\begin{array}{rcl}
y_0 &=& (0,-\eeth,d-\eeth,de-2) \\
y_1 &=& (\dth,1-\eeth,d-1-\eeth,(d-1)e-2+\dth) \\
 && \qquad\dots \\
y_i &=& (\textstyle{\frac id },i-\eeth,d-i-\eeth,(d-i)e-2+\textstyle{\frac id }) \\
 && \qquad\hbox{for $i=0,\dots,d-1$} \\
 y_{d-2} &=& (1-{\textstyle\frac2d},d-2-\eeth,2-\eeth,2e-1-{\textstyle\frac2d})
\end{array}
\notag
\end{equation}
\begin{equation}
\renewcommand{\arraystretch}{1.2}
\renewcommand{\arraycolsep}{0.3em}
\begin{array}{rcl}
y_{d-1} &=& (1-\dth,d-1-\eeth,1-\eeth,e-1-\dth) \\
y_{d} &=& (2-\dth,2d-1-{\textstyle\frac2e},1-{\textstyle\frac2e},e-2-\dth) \\
 && \qquad\dots \\
y_{d-2+i} &=& (i-\dth,id-1-\textstyle{\frac ie },1-\textstyle{\frac ie },e-i-\dth) \\
 && \qquad\hbox{for $i=1,\dots,e$} \\
y_{d+e-3} &=& (e-1-\dth,d(e-1)-2+\eeth,\eeth,1-\dth) \\
y_{d+e-2} &=& (e-\dth,de-2,0,-\dth) \\
\end{array}
\notag
\end{equation}
\end{exa}
Same checks; note especially the effect of the tag 3 at $y_{d-1}$.

\begin{exa}[Case $d=4$, $e=6$, $k=6$] \label{x!46}
\begin{equation}
\renewcommand{\arraystretch}{0.92}
\renewcommand{\arraycolsep}{0.5em}
\setcounter{MaxMatrixCols}{14}
\begin{matrix}
 &&&&& L & M & A & B \\[3pt]
 && x_0 & \quad = \quad &(& -1/4 & 0 & 505/4 & 483 &)\\
6 && x_1 & \quad = \quad &(& 0 & 1/6 & 22 & 505/6 &)\\
4 && x_2 & \quad = \quad &(& 1/4 & 1 & 23/4 & 22 &)\\
6 && x_3 & \quad = \quad &(& 1 & 23/6 & 1 & 23/6 &)\\
4 && x_4 & \quad = \quad &(& 23/4 & 22 & 1/4 & 1 &)\\
6 && x_5 & \quad = \quad &(& 22 & 505/6 & 0 & 1/6 &)\\
 && x_6 & \quad = \quad &(& 505/4 & 483 & -1/4 & 0 &)
\end{matrix}
\end{equation}
and
\begin{equation}
\renewcommand{\arraystretch}{0.92}
\renewcommand{\arraycolsep}{0.3em}
\setcounter{MaxMatrixCols}{14}
\begin{matrix}
 &&&&& L & M & A & B \\[3pt]
 && y_0 & \quad = \quad &(& 0 & -1/6 & 483 & 11087/6 &)\\
2 && y_1 & \quad = \quad &(& 1/4 & 5/6 & 1427/4 & 8189/6 &)\\
2 && y_2 & \quad = \quad &(& 1/2 & 11/6 & 461/2 & 5291/6 &)\\
3 && y_3 & \quad = \quad &(& 3/4 & 17/6 & 417/4 & 2393/6 &)\\
2 && y_4 & \quad = \quad &(& 7/4 & 20/3 & 329/4 & 944/3 &)\\
2 && y_5 & \quad = \quad &(& 11/4 & 21/2 & 241/4 & 461/2 &)\\
2 && y_6 & \quad = \quad &(& 15/4 & 43/3 & 153/4 & 439/3 &)\\
3 && y_7 & \quad = \quad &(& 19/4 & 109/6 & 65/4 & 373/6 &)\\
2 && y_8 & \quad = \quad &(& 21/2 & 241/6 & 21/2 & 241/6 &)\\
3 && y_9 & \quad = \quad &(& 65/4 & 373/6 & 19/4 & 109/6 &)\\
2 && y_{10} & \quad = \quad &(& 153/4 & 439/3 & 15/4 & 43/3 &)\\
2 && y_{11} & \quad = \quad &(& 241/4 & 461/2 & 11/4 & 21/2 &)\\
2 && y_{12} & \quad = \quad &(& 329/4 & 944/3 & 7/4 & 20/3 &)\\
3 && y_{13} & \quad = \quad &(& 417/4 & 2393/6 & 3/4 & 17/6 &)\\
2 && y_{14} & \quad = \quad &(& 461/2 & 5291/6 & 1/2 & 11/6 &)\\
2 && y_{15} & \quad = \quad &(& 1427/4 & 8189/6 & 1/4 & 5/6 &)\\
 && y_{16} & \quad = \quad &(& 483 & 11087/6 & 0 & -1/6 &)
\end{matrix} \label{eq!bigm}
\end{equation}
We read this table in several ways. Omitting the $A$ and $B$ columns
describes $\si_{LM}$ in the impartial basis. Notice the tag equations
\begin{align*}
&\hbox{bottom: $x_1y_0=x_0^4L$ and $x_0y_1=y_0M$;} \\
&\hbox{sides: $x_0x_2=x_1^6$, $x_1x_3=x_2^4$ and so on;} \\
&\hbox{top: $x_5y_{16}=L^{505}M^{1932}$ and $x_6y_{15}=x_5^5L^{373}M^{1427}$.}
\end{align*}
\end{exa}

Figure~\ref{f!sciss} plots the first two columns of \eqref{eq!bigm} as
``scissors'' controlled by the points $x_0=(-\dth,0)$ and
$y_0=(0,-\eeth)$ and the origin $(0,0)$ (implicit but crucial). To
describe it in words, the sequence of $y_i$ starts from $y_0$ and tries
to grow along the line $\La$ of slope
$1/[4,6,4,6,\dots]\doteq0.261387212$, without crossing it. It first
tries $x_0$ (slope $-\infty$), then $x_1$ (slope $0$) and $x_2$ (slope
$1/4$, so under $\La$), then takes one step back to $y_1=x_2y_0$ (slope
$3/10$, so above $\La$). Now $y_0,y_1,y_2,y_3,x_3$ is an arithmetic
progression of length $5=d+1$ with increment $x_2$ (and
$y_{i+1}=y_ix_2$, so $0y_0y_1x_2$, $0y_1y_2x_2$, etc., are
parallelograms); but $x_3$ (slope $6/23$) is below $L$; so take one step
back to $y_3$ and construct the next arithmetic progression
$y_3,y_4,y_5,y_6,y_7,x_4$ of length $6=e$ with increment $x_3$, and so
on. Compare Figure~\ref{f!Phi}, where the scissors were more open.

\begin{remark} \label{r!sciss}
The abstract continued fraction $[e,d,\dots]$ and its complementary
continued fraction $[2,2,\dots,3,\dots]$ has two different ``scissors''
embeddings into the $L,M$-plane (as the dots of Figures~\ref{f!Phi}
and~\ref{f!sciss}) and into the $A,B$-planes, and the Pretty Polytope
$\Pi(d,e,k)$ is just the diagonal embedding into the product.
\end{remark}

\subsection{The quotient $Q$ and the Padded Cell} \label{ss!pad}
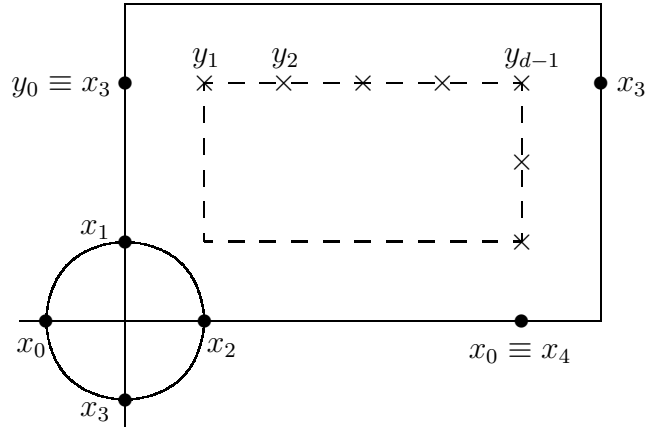
\begin{figure}[ht]
\begin{picture}(200,155)(-140,-35)
\put(-40,0){\line(1,0){220}}
\put(0,-40){\line(0,1){160}}
\put(0,120){\line(1,0){180}}
\put(180,0){\line(0,1){120}}
\put(30,30){\dashbox{6}(120,60){ }}
\put(25.4,87){$\times$} %y1
\put(25.4,98){$y_1$}
\put(55.4,87){$\times$}
\put(55.4,98){$y_2$}
\put(85.4,87){$\times$}
\put(115.4,87){$\times$}
\put(145.4,87){$\times$}
\put(143.4,98){$y_{d-1}$}
\put(145.4,57){$\times$}
\put(145.4,27){$\times$}
\put(0,30){\circle*{5}} % x1
\put(-17,33){$x_1$}
\put(0,-30){\circle*{5}} % x3 below line
\put(-17,-36){$x_3$}
\put(0,90){\circle*{5}} % x3
\put(-43,87){$y_0\equiv x_3$}
\put(180,90){\circle*{5}} % x3'
\put(186,87){$x_3$}
\put(30,0){\circle*{5}} % x2
\put(31,-12){$x_2$}
\put(-30,0){\circle*{5}} % x0 left of axis
\put(-41,-12){$x_0$}
\put(150,0){\circle*{5}} % x4
\put(130,-14){$x_0\equiv x_4$}
\qbezier(-30,0)(-28,28)(0,30)
\qbezier(30,0)(28,28)(0,30)
\qbezier(30,0)(28,-28)(0,-30)
\qbezier(-30,0)(-28,-28)(0,-30)
\end{picture}
\caption{The Padded Cell (with sides identified): the values of $x_i$
and $y_j$ in the torus $Q=\M/\M'$. The $x_i$ cycle around the 4 points
$(\pm\dth,0)$ and $(0,\pm\eeth)$ closest to the origin, while the $y_i$
walk around the path of Figure~\ref{f!pad}, performing $k-1$
quarter-circuits around the padding of the cell, starting from $x_3
\equiv y_0$. Each quarter-circuit takes place in steps of $x_i$ and has
endpoint $x_{i+1}$. \label{f!pad}}
\end{figure}

The exponents of $x_{0\dots k},y_{0\dots l}$ in Proposition~\ref{p!xys}
also behave in a characteristic way modulo the integers (see
Figure~\ref{f!pad}). To understand this, we write
$\M'=\Z\cdot(L,M,A,B)\subset \M$ for the sublattice generated by
$A,B,L,M$, and $Q=\M/\M'$ for the quotient. We think of $Q$ pictorially
as a fundamental domain in $\M$ for the translation lattice $\M'$, as in
Figure~\ref{f!pad}.

\begin{cor} \label{c!Q} \begin{enumerate}
\setcounter{enumi}{3}
\renewcommand{\labelenumi}{(\roman{enumi})}
\item $\Q\iso \Z/d \oplus \Z/e$, based by:
\begin{description}
\item{if $k=2\ka$ is even:}
\begin{equation}
x_0 \equiv (-\dth,0,\mp\dth,0);
\quad\hbox{and}\quad
y_0 \equiv (0,-\eeth,0,\pm\eeth);
\end{equation}

\item{if $k=2\ka+1$ is odd:}
\begin{equation}
x_0 \equiv (-\dth,0,0,\pm\dth)
\quad\hbox{and}\quad
y_0 \equiv (0,-\eeth,\pm\eeth,0),
\end{equation}
\item{where in either case $\pm=(-1)^\ka$.}
\end{description}
\item The classes in $Q$ of monomials $x_0,\dots,y_l$ are given as
follows (for even $k$):
\begin{equation}
\begin{gathered}
x_1 \equiv -y_0 \equiv (0,\eeth,0,\mp\eeth),
\quad
x_i \equiv -x_{i-2} \quad\hbox{for $i\ge2$} \\
\hbox{and}\quad y_{j+1}=y_j+x_{i(j)}
\end{gathered}
\end{equation}
for $j$ in the appropriate interval. In particular, in $Q$, the $x_i$
are periodic with period $4$, with $x_3\equiv y_0$.
\end{enumerate}
\end{cor}

Note that in $Q$, the different corner tags on the two long rectangles
say the same thing; thus
\begin{gather*}
x_1y_0=x_1^0A^\al B^\be=x_0^dL \quad\hbox{both give}\quad
x_1\equiv y_0\1 \in Q \\
x_0y_1=y_0^{-(e-1)}A^\ga B^\de=y_0M \quad\hbox{both give}\quad
y_0\equiv x_0y_1 \in Q
\end{gather*}
because $x_0^d,y_0^e\in\M$.

%%%%%%%%%%%%%%%%%%%%%%%%%%%%%%%%%%%%%%%%%%%%%%%%%%%%%%%%%

\section{Proof of Theorem~\ref{th!main}: main case} \label{s!pf}

We prove the existence of the diptych variety $V_{ABLM}$
for any pair of toric extensions of the tent $V_{AB}\supset T\subset V_{LM}$
arising from the Classification Theorem~\ref{th!d-e} {\bf under the
assumption that $d,e\ge 2$ and $de>4$}.

\subsection{Structure of the proof}

The proof of Theorem~\ref{th!main} builds a staircase: first, we
drop a chain of projections down from the top of $V_{AB}$
to eliminate the generator $x_{2\dots k}$ and $y_{2\dots l}$ one at a
time. This chain will serve as a guiding rail in the main construction;
it records the order of variables and the current state of the tags and
annotations as we eliminate them (Proposition~\ref{p!snu}): as each
$s_\nu=x_{i+1}$ or $y_{j+1}$ is eliminated from $V_{AB,\nu+1}$, it has
tag~1, and appears in an equation $s_\nu h_\nu=x_iy_j$ with its
neighbours, where $h_\nu=h_\nu(A,B)$ is the monomial in $A,B$ defined in
\ref{ss!hnu}.

We then build the 6-fold $V_{ABLM}$ up from the bottom, holding tight to
our guiding rail, the chain of projections of $V_{AB}$. Each step
$V_{\nu+1}\to V_\nu$ of the induction is a Kustin--Miller unprojection
(see \cite{PR}), and adjoins an unprojection variable $s_\nu=x_{i+1}$ or
$y_{j+1}$. The current $V_\nu$ is contained in the ambient space
$\Aff_\nu=\Aff^{i+j+6}_{\Span{x_{0\dots i},y_{0\dots j},A,B,L,M}}$. The
main point is to set up the unprojection divisor $D_\nu\subset V_\nu$;
we {\em define} it by the ideal
\begin{equation}
 I_{D_\nu} = (x_{0\dots i-1}, y_{0\dots j-1}, h_\nu),
\label{eq!IDnu}
\end{equation}
with $h_\nu(A,B)$ as in \ref{ss!hnu}, so that $D_\nu$ is the
hypersurface
\begin{equation}
D_\nu:(h_\nu(A,B)=0)\subset\Aff^6_{\Span{x_i,y_j,A,B,L,M}}.
\label{eq!Dnu}
\end{equation}
Thus $D_\nu$ is by definition the product of affine 4-space
$\Aff^4_{\Span{x_i,y_j,L,M}}$ with the monomial curve $h_\nu(A,B)=0$;
the elements $L,M$ form a regular sequence for $D_\nu$, and the section
$L=M=0$ in $D_\nu$ is the unprojection divisor $D_{AB,\nu}$ for $V_{AB,\nu+1}\to
V_{AB,\nu}$. The remaining issue is to prove that $D_\nu\subset V_\nu$,
or equivalently, that
\begin{equation}
I_{V_\nu}\subset I_{D_\nu}
= \bigl( x_{0\dots i-1}, y_{0\dots j-1},h_\nu \bigr).
\label{eq!DV}
\end{equation}
For this, rather than working with the actual equations of $V_\nu$ (that
we cannot always calculate in closed form, and include complicated
terms), we prove the stronger result: {\em any monomial in $x_{0\dots
i},y_{0\dots j},A,B,L,M$ with the same $\T$-weight as a generator of
$I_{V_\nu}$ is in $I_{D_\nu}=(x_{0\dots i-1}, y_{0\dots
j-1},h_\nu)$.}
Thus, every $\T$-homogeneous generator of $I_{V_\nu}$ is a sum of
monomials in $I_{D_\nu}$.

It turns out in the end, much to our regret, that our proof does here not
involve any explicit pentagrams or Pfaffians; however, they are important
in the constructions of~\cite{BR2} when $de=4$.

\subsection{The projection sequence of $V_{AB}$} \label{ss!pseq}

This section and the next set out facts and notation for the chains
of birational projections down from $V_{AB}$ and up from $V_{LM}$.
Either chain is provided by the blowdown of Proposition~\ref{prop!hj}(d)
applied to the conclusion $[a_2,\dots,b_1]=0$ of Corollary~\ref{cor!VAB}.

\begin{exa} \label{ex!decon}  Consider the long
rectangle of Figure~\ref{f!lrAB}. The concatenated continued fraction
$[4,2,1,3,2,2]=0$ is deconstructed as
\begin{equation}
[4,\underline{2,1,3},2,2] \mapsto [\underline{4,1,2},2,2] \mapsto
[\underline{3,1,2},2] \mapsto [\underline{2,1,2}] \mapsto [1,1]=0
\notag
\end{equation}
This is a recipe for a chain of birational projections, each
eliminating a monomial from $\si_{AB}$ with tag 1:
\begin{multline*}
\renewcommand{\arraycolsep}{.7em}
\begin{array}{rl}
&1 ^{\displaystyle \,B}\\
^{\displaystyle A \,} 2&3 \\
4&2 \\
2&2 \\
0&-1
\end{array}
\mapsto
\begin{array}{rl} \,\\
^{\displaystyle AB \,} 1&2^{\displaystyle \,B} \\
4&2 \\
2&2 \\
0&-1
\end{array}
\mapsto
\begin{array}{rl} \,\\
&1 ^{\displaystyle \,AB^2}\\
^{\displaystyle AB \,} 3&2 \\
2 &2 \\
0&-1
\end{array} \enspace\mapsto \\[10pt]
\renewcommand{\arraycolsep}{.7em}
\begin{array}{rl}
^{\displaystyle A^2B^3 \,}2&1 ^{\displaystyle \,AB^2} \\
2&2 \\
0&-1
\end{array}
\mapsto
\begin{array}{rl}
^{\displaystyle A^3B^5 \,}1& \\
2&1 ^{\displaystyle \,AB^2} \\
0&-1
\end{array}
\mapsto
\begin{array}{rl} \,\\
^{\displaystyle A^3B^5 \,}1&0 ^{\displaystyle \,A^4B^7} \\
0&-1
\end{array}
\end{multline*}
For example, on the second line, we read $x_1y_2=x_2^2A^2B^3$ and
$x_2y_1=y_2AB^2$ from the tags and annotation of the first rectangle,
that we can check against \eqref{eq!Laur}. Each rectangle is the
monomial cone $\si_{AB,\nu}$ of a Gorenstein affine toric variety
$V_{AB,\nu}$ with the given monomials in $A,B$ as annotations, and each
step $V_{AB,\nu+1}\to V_{AB,\nu}$ is a birational projection.
\end{exa}

\subsubsection{Order of monomials} \label{ss!order}

Our construction inverts this type of chain, up from a codimension~2
complete intersection in $x_0,x_1,y_0,y_1,A,B$, adding $x_2$, $y_2$ and
so on one at a time, to recover $V_{AB}$. For this, we order the $k+l-2$
steps {\em inverse to the elimination} of the monomials $x_{2\dots
k},y_{2\dots l}$; that is, we rename the $n$th eliminated monomial
$s_\nu$ with $\nu=k+l-2-n$, so that $s_0=x_2$ and $s_1=y_2$. We work by
induction on this $\nu$. At the same time, we name the annotation
$h_\nu$ on the monomial $s_\nu$ as it is eliminated; the chain starts
from the top with
\begin{equation}
s_{k+l-3}=y_l \quad\hbox{and}\quad h_{k+l-3}=B \quad \hbox{and} \quad b_l=1.
\end{equation}
(This uses the main case hypothesis $d,e\ge2$ so that $b_l=1$.)

Thus in Example~\ref{ex!decon},
$[s_0,s_1,s_2,s_3,s_4]=[x_2,y_2,y_3,x_3,y_4]$ and
\begin{equation}
[h_0,h_1,h_2,h_3,h_4]=[A^3B^5,AB^2,AB^2,AB,B].
\notag
\end{equation}
The scissors of Figure~\ref{f!sciss} strongly suggest this ordering of
the monomials, although there is a choice to make at the end between
$y_1$ and $x_2$, which both have tag 1; we always eliminate $s_1=x_2$.

\subsubsection{The projection $V_{AB,\nu+1}\to V_{AB,\nu}$ and the bar
$x_i\frac{\qquad}{}y_j$}
\label{ss!snu}
The projection sequence gives cones $\si_{AB,\nu}$ that depend on the
induction parameter $\nu$. The top corners of each $\si_{AB,\nu}$ are
monomials $x_i$ and $y_j$ with $i=i(\nu)$ and $j=j(\nu)$
(Table~\ref{t!spsh} keeps track of these functions), and we know the
equations of $V_{AB,\nu}$ including
\begin{equation}
x_{i-1}y_j=x_i^{\al_\nu}A_\nu
\quad\hbox{and}\quad
x_iy_{j-1}=y_j^{\be_\nu}B_\nu,
\label{eq!mod}
\end{equation}
given by the tags and annotations at $x_i$ and $y_j$ in $V_{AB,\nu}$
as in Figure~\ref{f!mod}.
We think of this action happening at the top of a sub-rectangle,
which we refer to as the {\em bar} $x_{i}\frac{\qquad}{}y_{j}$;
the bar cascades down the long rectangle as variables are projected away
(see Figure~\ref{f!pfpi}),
and the tag equations \eqref{eq!mod} at each bar provide the key pieces
of quantitative data about the convexity of $V_{AB}$ that we use
throughout the proof.
\begin{figure}[h] \[
\begin{picture}(120,35)
\put(17,-15){\line(0,1){40}}
\put(17,25){\line(1,0){80}}
\put(97,-15){\line(0,1){40}}
\put(17,1){\circle*{6}}
\put(97,1){\circle*{6}}
\put(97,25){\circle*{6}}
\put(17,25){\circle*{6}}
\put(-11,-3){$x_{i-1}$}
\put(22,10){$\al_\nu$}
\put(1,20){$x_i$}
\put(-4,32){$A_\nu$}
\put(104,32){$B_\nu$}
\put(103,20){$y_j$}
\put(103,-3){$y_{j-1}$}
\put(82,10){$\be_\nu$}
\end{picture} \]
\caption{The {\em bar} $x_i\frac{\qquad}{}y_j$ at the top of
$\si_{AB,\nu}$, with tag equations \eqref{eq!mod}.
\label{f!mod}}
\end{figure}

\begin{prop} \label{p!snu}  The chain of projections
$V_{AB,\nu+1}\to V_{AB,\nu}$ reduces $V_{AB}$ down to a codimension~$2$
complete intersection
$V_{AB,0}\subset\Aff^6_{\Span{x_0,x_1,y_0,y_1,A,B}}$. The step
$V_{AB,\nu+1}\to V_{AB,\nu}$ eliminates $s_\nu=x_{i+1}$ or $y_{j+1}$,
with two possible cases for the top of $\si_{AB,\nu+1}$:
\begin{equation}
\renewcommand{\arraycolsep}{1em}
\hbox{either}\quad
\begin{array}{rl}
s_\nu \\
x_i & y_j \\
x_{i-1} & y_{j-1}
\end{array}
\quad\hbox{or}\quad
\begin{array}{rl}
&s_\nu \\
x_i & y_j \\
x_{i-1} & y_{j-1}
\end{array}
\notag
\end{equation}
In the left case $s_\nu=x_{i+1}$, the top of $\si_{AB,\nu}$ and of
$\si_{AB,\nu+1}$ are related by
\begin{equation}
A_\nu=A_{\nu+1}, \enspace B_\nu=A_{\nu+1}B_{\nu+1}, \enspace
\al_{\nu+1}=1, \enspace
\al_\nu=a_i-1,
\enspace \be_\nu=\be_{\nu+1}-1,
\label{eq!nu+1l}
\end{equation}
and similarly in the right case by
\begin{equation}
A_\nu=A_{\nu+1}B_{\nu+1}, \enspace B_\nu=B_{\nu+1}, \enspace
\al_\nu=\al_{\nu+1}-1, \enspace \be_{\nu+1}=1, \enspace
\be_\nu=b_j-1.
\label{eq!nu+1r}
\end{equation}
\end{prop}

\subsubsection{Choice of $h_\nu(A,B)$ and the unprojection divisor
$D_{AB,\nu}\subset V_{AB,\nu}$} \label{ss!hnu}
Proposition~\ref{p!snu} described the projection $V_{AB,\nu+1}\to
V_{AB,\nu}$ that eliminates the variable $s_\nu$; inverting this, we
construct $V_{AB,\nu+1}$ as an unprojection from $V_{AB,\nu}$ adjoining
$s_\nu$. For this, we set $h_\nu=\hcf(A_\nu,B_\nu)$, equal to $A_\nu$ or
$B_\nu$ by (\ref{eq!nu+1l}--\ref{eq!nu+1r})
and define
$D_{AB,\nu}\subset\Aff^{i+j+4}_{\Span{x_{0\dots i},y_{0\dots j},A,B}}$
by the ideal $(x_{0\dots i-1},y_{0\dots j-1},h_\nu)$; thus $D_{AB,\nu}$
is the hypersurface $(h_\nu=0)\subset\Aff^4_{\Span{x_i,y_j,A,B}}$.

\begin{cla} The ideal of $V_{AB,\nu}$ is contained in the ideal
$(x_{0\dots i-1},y_{0\dots j-1},h_\nu)$ of $D_{AB,\nu}$, or in other
words, $D_{AB,\nu}\subset V_{AB,\nu}$. \end{cla}

\begin{pf} We know that $I_{V_{AB,\nu}}$ is generated by equations for
$x_{i'}y_{j'}$, $x_{i'}x_{i''}$ and $y_{j'}y_{j''}$ for suitable values
of the indexes $i',i'',j',j''$. Consider for example an equation
$x_iy_{j'}=x_i^\xi y_j^\eta A^\al B^\be$ for some $j'<j$. First $\xi=0$,
for otherwise dividing by $x_i$ gives a monomial expression for
$y_{j'}$ that contradicts Figure~\ref{f!siAB}, where $\Span{y_{0\dots
l}}$ is a 2-face of $\si_{AB}$. Substitute for $x_i$ from the tag equation
$x_iy_{j-1}=y_j^{\be_j}B_\nu$ to give
\begin{equation}
y_{j'}y_j^{\be_j-\eta}y_{j-1}\1 = A^\al B^\be B_\nu\1.
\label{eq!both1}
\end{equation}
Now both sides of \eqref{eq!both1} are 1, since the 4-dimensional vector
space $\M_\Q$ is the direct sum of the 2-dimensional subspace spanned by
$y_{0\dots l}$ and that spanned by $A,B$  (compare Figure~\ref{f!siAB}).
Therefore $A^\al B^\be=B_\nu$, and both sides of our equation are in the
ideal. The other equations are similar. \end{pf}

The initial case $n=0$ or $\nu=k+l-2$ is $V_{AB,\nu}=V_{AB}$; in our
construction of $V_{ABLM}$, it is the final goal: if we reach it, there
is nothing more to check. Then $A=A_\nu$, $B=B_\nu$, $h_\nu=1$, and
divisibility by $h_\nu$ is trivial.

\subsection{Crosses, pitchforks and pentagrams}
\subsubsection{The spreadsheet for $V_{AB}$}

Our construction of $V_{\nu+1}$ from $V_\nu$ reverses
the projection sequence down from the top of $V_{AB}$. Our proof also
needs information derived from the projection sequence up from the
bottom of $V_{LM}$. Thus in Extended Example~\ref{ex!intr}, we
deconstructed $V_{LM}$ by eliminating $y_0,y_1,y_2,x_0,x_1$ from the
bottom of Figure~\ref{f!lrLM}. Here we establish how the two projection
sequences interleave, as an exercise in patient bookkeeping.

Table~\ref{t!spsh} gives the function $i=i(\nu)$, $j=j(\nu)$ of
\ref{ss!snu} describing the top of $V_{AB,\nu}$ as in
Figure~\ref{f!mod}. The table repeats periodically with period $d+e-2$,
or alternate half periods of $d-1$, $e-1$. We set $v=\nu\mod d+e-2$ and
write $\nu=C(d+e-2)+v$.
\begin{table}[ht]
\begin{gather*}
\renewcommand{\arraycolsep}{.7em}
\renewcommand{\arraystretch}{1.25}
\begin{array}{|c|c|cl|}
\hline
v & i & j & \\
\hline
\hline
0 & 2C+1 & (d+e-4)C+1  & \\
\hline
a & 2C+2 & (d+e-4)C+a & \hbox{for $1\le a\le e-1$} \\
\hline
e+b-1 & 2C+3 & (d+e-4)C+e-2+b & \hbox{for $1\le b\le d-1$}\\
\hline
\hbox{Final} & k=2\ka+1 & l=(d+e-4)\ka+2 & \\
\hline
\end{array}
\\[5pt]
\renewcommand{\arraycolsep}{.9em}
\renewcommand{\arraystretch}{1.4}
\begin{array}{|c|c|cl|}
\hline
v & i & j & \\
\hline
\hline
0 & 2C+1 & (d+e-4)C+1  & \\
\hline
a & 2C+2 & (d+e-4)C+a & \hbox{for $1\le a\le d-1$} \\
\hline
d+b-1 & 2C+3 & (d+e-4)C+d-2+b & \hbox{for $1\le b\le e-1$}\\
\hline
a & 2C+2 & (d+e-4)C+a & \hbox{for $1\le a\le d-1$} \\
\hline
\hbox{Final} & k=2\ka & l & \\
\hline
\end{array}
\end{gather*}
\caption{Numbering the unprojection sequence for $V_{AB}$. The even case
$k=2\ka$ has one fewer half round. The final line is irregular: it adds a
final $y_l$ instead of $x_{k+1}$ with $l=(d+e-4)\ka+2$ or
$l=(d-2)\ka+(e-2)(\ka-1)+2$. \label{t!spsh}}
\end{table}

The starting point $\nu=0$ is $V_{AB,\nu}$ with $x_1,y_1$ at its top
bar. The table is split into two, the $k$ odd and $k$ even cases;
we describe the odd case $k=2\ka+1$. Set $C=0$ and enter the first
round: the line $v=a=1$ adds an $x_i$, then $a=2,\dots,e-1$ is a half
round that adds $e-2$ terms $y_j$; similarly, the line $v=e$ (so $b=1$)
adds an $x_i$ and then $b=2,\dots,d-1$ is a half round that adds $d-2$
terms $y_j$. We then increment $C\mapsto C+1$ and loop. Each half round
adds one $x_i$ and $d-2$ or $e-2$ terms $y_j$. There are $k-1$ half
rounds, ending with $\nu=(d+e-2)\ka$ if $k=2\ka+1$ or
$\nu=(d-1)\ka+(e-1)(\ka-1)$ if $k=2\ka$.

The above treatment assumes that we are in the main case $d,e\ge2$;
everything remains true when $d$ or $e$ or both are~2. Then the
intervals $2\le a\le d-1$ or $2\le b\le e-1$ are empty, so the
corresponding half periods add one $x_i$ and no $y_j$.

\subsubsection{Comparing the projection sequences for $V_{AB}$ and
$V_{LM}$}

We want to compare the bars $x_i,y_j$ at the top of $V_{AB,\nu}$ with
the corresponding thing at the bottom of $V_{LM}$ after a number of
projections. To see this, we divide the monomials $y_j$ up into
intervals according to the lines of Table~\ref{t!spsh}, writing
$Y_{i-1}$ for the $i$th half
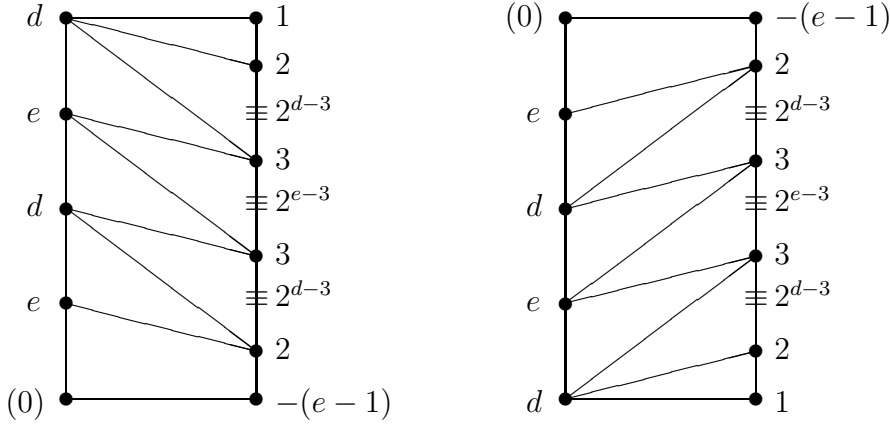
\begin{figure}[ht]
\begin{picture}(150,157)(-5,0)
% the box
\put(45,7){\line(0,1){144}}
\put(117,7){\line(0,1){144}}
\put(45,7){\line(1,0){72}}
\put(45,151){\line(1,0){72}}
% the left side tags
\put(45,7){\circle*{5}} \put(22,2){(0)}
\put(45,43){\circle*{5}} \put(30,40){$e$}
\put(45,79){\circle*{5}} \put(30,76){$d$}
\put(45,115){\circle*{5}} \put(30,112){$e$}
\put(45,151){\circle*{5}} \put(30,148){$d$}
% the right side tags
\put(117,7){\circle*{5}}
\put(124,2){$-(e-1)$}
\put(117,25){\circle*{5}}
\put(124,22){2}
\put(117,61){\circle*{5}}
\put(124,58){3}
\put(117,97){\circle*{5}}
\put(124,94){3}
\put(117,133){\circle*{5}}
\put(124,130){2}
\put(117,151){\circle*{5}}
\put(124,148){1}
\put(112.5,112){$\equiv$}
\put(124,112){$2^{d-3}$}
\put(112.5,78){$\equiv$}
\put(124,78){$2^{e-3}$}
\put(112.5,42){$\equiv$}
\put(124,42){$2^{d-3}$}
% diagonal lines
\put(45,151){\line(4,-1){72}}
\put(45,151){\line(4,-3){72}}
\put(45,115){\line(4,-1){72}}
\put(45,115){\line(4,-3){72}}
\put(45,79){\line(4,-1){72}}
\put(45,79){\line(4,-3){72}}
\put(45,43){\line(4,-1){72}}
\end{picture}
\begin{picture}(150,110)(-40,0)
% the box
\put(45,7){\line(0,1){144}}
\put(117,7){\line(0,1){144}}
\put(45,7){\line(1,0){72}}
\put(45,151){\line(1,0){72}}
% the left side tags
\put(45,7){\circle*{5}} \put(22,148){(0)}
\put(45,43){\circle*{5}} \put(30,40){$e$}
\put(45,79){\circle*{5}} \put(30,76){$d$}
\put(45,115){\circle*{5}} \put(30,112){$e$}
\put(45,151){\circle*{5}} \put(30,2){$d$}
% the right side tags
\put(117,7){\circle*{5}}
\put(124,148){$-(e-1)$}
\put(117,25){\circle*{5}}
\put(124,22){2}
\put(117,61){\circle*{5}}
\put(124,58){3}
\put(117,97){\circle*{5}}
\put(124,94){3}
\put(117,133){\circle*{5}}
\put(124,130){2}
\put(117,151){\circle*{5}}
\put(124,2){1}
\put(112.5,112){$\equiv$}
\put(124,112){$2^{d-3}$}
\put(112.5,78){$\equiv$}
\put(124,78){$2^{e-3}$}
\put(112.5,42){$\equiv$}
\put(124,42){$2^{d-3}$}
% diagonals
\put(45,7){\line(4,1){72}}
\put(45,7){\line(4,3){72}}
\put(45,43){\line(4,1){72}}
\put(45,43){\line(4,3){72}}
\put(45,79){\line(4,1){72}}
\put(45,79){\line(4,3){72}}
\put(45,115){\line(4,1){72}}
\end{picture}
\caption{Projecting $V_{AB}$ from the top and $V_{LM}$ from the bottom}
\label{f!pfpi}
\end{figure}
period. In more detail, for $k$ even, the line for even $i=2C+2$ gives
the interval
\begin{equation}
Y_{i-1}=\bigl\{y_j \bigm| \hbox{for $j\in [n_i+1,\dots,n_i+d-1]$} \bigr\}
\end{equation}
where $n_i=(d+e-4)\frac{i-2}2$; similarly, the line $i'=2C+3$ gives
\begin{equation}
Y_{i'-1} = \{y_j \bigm| \hbox{for $j\in [n_{i'}+1,\dots,n_{i'}+e-1]$} \bigr\}
\end{equation}
where $n_{i'}=(d+e-4)\frac{i'-3}2+d-2$.

Notice the adjacency between the intervals: the last entry $n_i+d-1$ of
$Y_{i-1}$ equals the first entry $n_{i'}+1$ of the following interval
$Y_{i'}$ with $i'=i+1$, and vice versa. For $d$ or $e=2$, the interval
$Y_i$ reduces to one element (which, in Figure~\ref{f!pfpi}, is tagged
with 4, rather than 3).

\begin{lem} \label{l!bar}
The bars at the top of $V_{AB,\nu}$ are precisely $x_{i+1},y_j$ with
$j\in Y_i$.

The bars at the bottom of $V_{LM,\nu'}$ (after projecting out $\nu'$
monomials from $V_{LM}$, starting with $y_0$) are precisely
$x_{i-1},y_j$ with $j\in Y_i$. See Figure~\ref{f!pfpi}. \end{lem}

The first clause is a more digestible rephrasing of the information
contained in Table~\ref{t!spsh} about the order of projection. The
projection sequence of $V_{LM}$ from the bottom is enumerated by a
symmetric spreadsheet, which proves the second clause.

The following simple consequence is a key point of our proof in
\ref{ss!pf}.

\begin{cor} \label{c!kil} Suppose that we project out $n_1$ monomials
from the top of $V_{AB}$ down to the top bar $x_i,y_j$ and $n_2$
monomials from the bottom of $V_{LM}$ up to the bottom bar
$x_{i'},y_{j'}$, where $n_1+n_2=k+l-2$, so that just $4$ monomials
remain. Then $i'<i$ and $j'\le j$.

Equivalently, either $i'=i-1$ and $j'=j-1$ or $i'=i-2$ and $j'=j$, so
that any such projection leads to a ``cross'' or ``pitchfork'' of the
shape
\begin{equation}
\label{eq!cross}
\vphantom{\begin{matrix} 0\\0\\0 \end{matrix}}
\setlength{\unitlength}{1mm}
\begin{picture}(32,9)(6,10)
\put(6,18){\circle*{1.44}} % x_1
\put(6,6){\circle*{1.44}} % x_0
\put(30,18){\circle*{1.44}} % y_1
\put(30,6){\circle*{1.44}} % y_0
\put(6,18){\line(2,-1){24}}
\put(6,6){\line(2,1){24}}
\put(0,17.2){$x_i$}
\put(-4,5.2){$x_{i-1}$}
\put(32.5,17.2){$y_j$}
\put(32.5,5.2){$y_{j-1}$}
\end{picture}
\qquad
\begin{picture}(5,9)(0,10)
\put(0,10.5){or}
\end{picture}
\qquad
\begin{picture}(38,9)(1,10)
\put(6,19){\circle*{1.44}} % x_2
\put(6,12){\circle*{1.44}} % x_1
\put(6,5){\circle*{1.44}} % x_0
\put(30,12){\circle*{1.44}} % y_1
\put(6,12){\line(1,0){24}}
\qbezier(6,19)(18,12)(6,5)
\put(0,18){$x_i$}
\put(-4,11.2){$x_{i-1}$}
\put(-4,3.6){$x_{i-2}$}
\put(32.5,11.4){$y_j$}
\end{picture}
\end{equation}
\end{cor}

The same phenomenon was already implicit in the cascade of pentagrams of
Example~\ref{ex!intr}; we include this, although it is not essential for
our proof.

\begin{cor}
Projecting out $n_1$ monomials from the top of $V_{AB}$ and $n_2$ from
the bottom of $V_{LM}$ with $n_1+n_2=k+l-3$ gives a pentagram of one of
the two shapes
\begin{equation}
\label{eq!x5}
\vphantom{\begin{matrix} 0\\0\\0\\0 \end{matrix}}
\setlength{\unitlength}{1mm}
\begin{picture}(38,9)(-2,6)
\put(6,16){\circle*{1.44}} % x_2
\put(6,8){\circle*{1.44}} % x_1
\put(6,0){\circle*{1.44}} % x_0
\put(30,8){\circle*{1.44}} % y_1
\put(30,0){\circle*{1.44}} % y_0
\put(6,8){\line(1,0){24}}
\put(6,8){\line(3,-1){24}}
\put(6,0){\line(3,1){24}}
\qbezier(6,16)(16,8)(6,0)
\qbezier(6,16)(21,9.67)(30,0)
\end{picture}
\qquad
\begin{picture}(5,11)(0,10)
\put(0,10.5){or}
\end{picture}
\qquad
\begin{picture}(38,13)(-3,10)
\put(6,24){\circle*{1.44}} % x_2
\put(6,16){\circle*{1.44}} % x_2
\put(6,8){\circle*{1.44}} % x_1
\put(6,0){\circle*{1.44}} % x_0
\put(30,12){\circle*{1.44}} % y_1
\put(6,16){\line(6,-1){24}}
\put(6,8){\line(6,1){24}}
\qbezier(6,24)(13,16)(6,8)
\qbezier(6,16)(13,8)(6,0)
\qbezier(6,24)(22,12)(6,0)
\end{picture}
\end{equation}
\end{cor}

\subsection{Proof by induction} \label{ss!pf}
We construct $V=V_{ABLM}$ by serial unprojection. The induction starts
from the codimension~2 complete intersection
\begin{equation}
V_0 \subset \Aff^8_{\Span{x_0,x_1,y_0,y_1,A,B,L,M}}
\notag
\end{equation}
defined by
\begin{equation}
x_1y_0=T_{x_0}(V_{AB}) + T_{x_0}(V_{LM})
\quad\hbox{and}\quad
x_0y_1=T_{y_0}(V_{AB}) + T_{y_0}(V_{LM})
\notag
\end{equation}
where $T_{x_0}(V_{AB})$ is the righthand side of the tag equation at
$x_0$ in $V_{AB}$, and similarly for the other three terms. Clearly
$V_0$ is Gorenstein and $A,B,L,M$ is a regular sequence, with the
regular section $L=M=0$ in $V_0$ the variety $V_{AB,0}$. We use the
following elementary fact about unprojection.

\begin{lem} \label{l!com}
Unprojection commutes with regular sequences: let $X,D$ be as in
\cite{PR}, Theorem~1.1 and $Y\to X$ the unprojection of $D$ in $X$.
Suppose that $z_1,\dots,z_r\in\Oh_X$ is a regular sequence for $X$ and
for $D$. Then $z_1,\dots,z_r$ is also a regular sequence for $\Oh_Y$,
and $Y_z$ is the unprojection of $D_z$ in $X_z$, where
$Y_z:(z_1=\cdots=z_r=0)\subset Y$ and similarly for $D_z$ and $X_z$. \qed \end{lem}

\begin{indass} \label{indas}  We own a variety
$V_\nu=V_{ABLM,\nu}$ having a $\T$-action, together with a regular
sequence $L,M$ made up of\/ $\T$-eigenfunctions such that $V_\nu\cap
(L=M=0)=V_{AB,\nu}$. \end{indass}

We start with $\nu=0$, and $V_{ABLM,\nu}=V_0$ as above. The induction
has $k+l-2$ steps, adjoining $x_{2\dots k}$ and $y_{2\dots l}$ in the
order determined in \ref{ss!order}. When $\nu$ reaches $k+l-2$ then
$V_{ABLM}=V_\nu$ and we are finished. Otherwise, if $\nu< k+l-2$, the
induction step consists of proving that $V_\nu$ has a divisor $D_\nu$ on
which $L,M$ is a regular sequence, and the section $D_\nu\cap(L=M=0)$ is
the divisor $D_{AB,\nu}\subset V_{AB,\nu}$.

If $\nu<k+l-2$, by \ref{ss!hnu}, the step $V_{AB,\nu+1}\to
V_{AB,\nu}$ of the chain down from $V_{AB}$ is the unprojection
adjoining the element $s_\nu$ with unprojection ideal $\bigl(x_{0\dots
i-1},y_{0\dots j-1},h_\nu\bigr)$, where $h_\nu$ is the monomial in $A,B$
defined in \ref{ss!hnu}. We seek to imitate this for the 6-fold
$V_\nu$; for this, define $D_\nu$ by
\begin{equation}
D_\nu \subset \Aff^8_{\Span{x_{0\dots i},y_{0\dots j},A,B,L,M}}
\quad\hbox{with ideal}\quad
I_{D_\nu}=\bigl(x_{0\dots i-1},y_{0\dots j-1},h_\nu\bigr).
\notag
\end{equation}
Clearly, it is the hypersurface
$D_\nu:(h_\nu=0)\subset\Aff^6_{\Span{x_i,y_j,A,B,L,M}}$, and is the product of
$\Aff^4_{\Span{x_i,y_j,L,M}}$ with the plane curve $h_\nu(A,B)=0$. The issue is
to prove that $D_\nu\subset V_\nu$.

\begin{prop}[Key point] \label{p!k}  $I_{V_\nu}\subset
I_{D_\nu}=(x_{0\dots i-1},y_{0\dots j-1},h_\nu)$ for every $\nu<k+l-2$.
\end{prop}

We prove this by a general argument on $\T$-weights of monomials that may
appear in a relation, without any need to analyse the actual equations
of $V_\nu$. We introduce the notation $R(\nu)$ for the $\T$-weights of
homogeneous generators of $I_{V_{AB,\nu}}$ or equivalently, of
$I_{V_\nu}$ (by $\T$-equivariance); we write $f\in R(\nu)$ to indicate
that $f$ is a homogeneous polynomial with $\T$-weight in $R(\nu)$. The
precise statement we prove is the following:

\begin{claim} \label{cl!div}
Any monomial $x_i^\xi y_j^\eta A^\al B^\be L^\la M^\mu \in
R(\nu)$ is divisible by $h_\nu$.
(We emphasise the prevailing hypotheses: $d,e\ge 2$ and $de>4$.)
\end{claim}

Recall that $h_\nu=\hcf(A_\nu,B_\nu)$; we usually prove divisibility by
$A_\nu$ or $B_\nu$. By definition, any alleged monomial in $R(\nu)$ is
$\T$-equivalent to a relation in $I_{V_{AB,\nu}}$ for $x_{i'}y_{j'}$ or
$x_{i'}x_{i''}$ or $y_{j'}y_{j''}$. The main mechanism of the proof is
to compare it with one of the two equations \eqref{eq!mod}, or more
precisely, with one of the model monomials
\begin{equation}
x_{i-1}y_j \simT x_i^{\al_\nu}A_\nu
\quad\hbox{and}\quad
x_iy_{j-1} \simT y_j^{\be_\nu}B_\nu,
\label{eq!mod2}
\end{equation}
coming from the top corners of $V_{AB,\nu}$ as in Figure~\ref{f!mod}.

\step1 {\em Claim~\ref{cl!div} holds for every monomial in $R(\nu-1)$.}
\addcontentsline{toc}{subsubsection}{Step 1. Monomials in $R(\nu-1)$}
Indeed, it is divisible by $h_{\nu-1}$ by induction, and by
(\ref{eq!nu+1l}--\ref{eq!nu+1r}) the $h_\nu$ increase as $\nu$ decreases.

\step2 \addcontentsline{toc}{subsubsection}{Step 2. Crossover monomials
$x_{i'}y_j$ and $x_iy_{j'}$} The first actual calculation in the proof:
{\em Claim~\ref{cl!div} holds for all the monomials $x_{i'}y_j$ with
$i'=0,\dots,i-1$ and $x_iy_{j'}$ with $j'=0,\dots,j-1$ appearing in {\em
cross-over} relations.}

\begin{pf} We write out the proof for $x_{i'}y_j$ in detail as a model
case. The method is to compare an alleged monomial
\begin{equation}
x_i^\xi y_j^\eta m\simT x_{i'}y_j \in R(\nu),
\quad\hbox{where $m$ is a monomial in $A,B,L,M$}
\notag
\end{equation}
with the known monomial $x_i^{\al_\nu}A_\nu\simT x_{i-1}y_j$ from
\eqref{eq!mod2}. (Of course, we never need consider monomials
which involve $x$ or $y$ variables from $I_{D_\nu}$.)
We have $\eta=0$: otherwise dividing both sides by
$y_j$ contradicts Corollary~\ref{c!ming}(iii). Consider
\begin{equation}
\frac{x_{i'}}{x_{i-1}} \simT x_i^{\xi-\al_\nu} \frac m{A_\nu}.
\end{equation}
By Corollary~\ref{c!ming}(ii) and the fact that $i'\le i-1$, the
lefthand side has $L,M$ exponents
$\pi_{LM}(\frac{x_{i'}}{x_{i-1}})\le0$ (see \eqref{eq!pi} for the
notation $\pi_{LM}$ and $\pi_{AB}$); thus $\al_\nu\ge\xi$, and the
equivalence takes the form
\begin{equation}
x_i^{\al_\nu-\xi}\frac{x_{i'}}{x_{i-1}} \simT \frac m{A_\nu}
\quad\hbox{with}\quad \al_\nu\ge \xi.
\label{eq!gt-1}
\end{equation}
Now for the same reason, $\pi_{AB}(\frac{x_{i'}}{x_{i-1}})\ge0$. The
same goes for $x_i^{\al_\nu-\xi}$, except for the case $x_i=x_k$, at the
top left of the rectangle for $V_{AB}$.
In the former case, we are done: $\pi_{AB}(m/A_\nu)\ge 0$ so $m$
is divisible by $A_\nu$ as required.

The initial case $x_i=x_k$ is important: $\pi_{AB}(x_k)=(-\frac1d,0)$ by
Proposition~\ref{p!xys}; because of the negative exponent, we cannot
get our conclusion by convexity alone. Instead we use a congruence
argument on the Padded Cell Figure~\ref{f!pad}: in
fact, the negative exponent is the smallest possible value $-\frac1d$,
and we claim that $\al_\nu$ is one of $d-1,d-2,\dots,1$. Indeed, if
$\nu=k+l-2$ we are at the end of the induction, and there is nothing to
prove. Otherwise, the tag at $x_k$ has decreased by at least one from
its pristine value $d$. It follows that the lefthand side of
\eqref{eq!gt-1} has $A$ exponent $>-1$ and $B$ exponent $\ge0$. On the
other hand, the righthand side of \eqref{eq!gt-1} is a Laurent monomial.
Therefore also in the initial case $m$ is divisible by $A_\nu$, as required.

The argument for $x_iy_{j'}$ is similar but slightly easier. Suppose
that
\begin{equation}\label{eq!topc}
x_iy_{j'} \simT x_i^\xi y_j^\eta m \quad\hbox{with}\quad
m=A^\al B^\be L^\la M^\mu.
\end{equation}
First $\xi=0$, because otherwise dividing through by $x_i$ would
contradict Corollary~\ref{c!ming}(iii). Next, dividing through by the
monomials in the second expression of \eqref{eq!mod2} gives
\begin{equation}
\frac {y_{j'}}{y_{j-1}}=y_j^{\eta-\be_\nu} \times \frac m {B_\nu}.
\end{equation}
As before, since $j'\le j-1$, Corollary~\ref{c!ming}(ii) gives that
$\pi_{LM}(\frac {y_{j'}}{y_{j-1}})\le0$. Therefore $\eta-\be_\nu\le0$.
Taking that term to the lefthand side gives
\begin{equation}
y_j^{\be_\nu-\eta}\frac{y_{j'}}{y_{j-1}}=\frac m {B_\nu}
\quad\hbox{with}\quad \be_\nu\ge\eta.
\end{equation}
Now $j'\le j-1$, so $\pi_{AB}(\frac{y_{j'}}{y_{j-1}})\ge0$; the same
goes for $y_j$ except if $j=l$ and $y_j$ is at the top of the long
rectangle, and we are finished, with $V_{\nu}=V_{ABLM}$. Therefore the
exponents of $A,B$ on the lefthand side are $\ge0$, and hence $m$ is
divisible by $B_\nu$.

This proves Step~2. \end{pf}

The proof of Step~2 used Corollary~\ref{c!ming}(ii) to compare the
exponents of $x_i/x_{i'}$ and $y_j/y_{j'}$, with typical implication
$i>i'\Rightarrow\pi_{LM}(x_i)>\pi_{LM}(x_{i'})$. For Step~3 we need a
similar comparison for monomials $x_i/y_{j'}$ and $y_j/x_{i'}$. Care is
needed here to distinguish the order of monomials in the projection
sequences from the top of $V_{AB}$ and from the bottom of $V_{LM}$: the
$L,M$ exponents behave monotonically in the projection sequences of
$V_{AB}$, and vice versa.

\begin{lem} \label{l!cmp}
Given two monomials $m_1,m_2\in\{x_{0\dots k},y_{0\dots
l}\}$, suppose that the projection sequence for $V_{AB}$ eliminates
$m_1$ before $m_2$; then
\begin{equation}
\pi_{LM}(m_1)\ge\pi_{LM}(m_2).
\end{equation}
Similarly, if the projection sequence for $V_{LM}$ eliminates $m_1$
before $m_2$ then
\begin{equation}
\pi_{AB}(m_1)\ge\pi_{AB}(m_2).
\end{equation}
\end{lem}

See Scissors, Figure~\ref{f!sciss} for a picture. The proof is simply to
observe that when a variable is introduced in an unprojection sequence,
it appears linearly in the new tag equation at its corner. Example~\ref{x!46}
provides a numerical sanity check, with the respective orders of
elimination
\begin{align*}
V_{AB}: &\ y_{16},y_{15},y_{14},x_6,y_{13},y_{12},y_{11},y_{10},x_5,
y_9,y_8,x_4,y_7,y_6,y_5,y_4,x_3,y_3,y_2,x_2; \\
V_{LM}: &\ y_0,y_1,y_2,x_0,y_3,y_4,y_5,y_6,x_1,y_7,y_8,x_2,y_9,
y_{10},y_{11},y_{12},x_3,y_{13},y_{14},x_4.
\end{align*}

\step3 \addcontentsline{toc}{subsubsection}{Step 3. Monomials $y_jy_a$
with $a=0,\dots,j-2$} {\em Claim~\ref{cl!div} holds for all monomials
$y_jy_a$ with $a=0,\dots,j-2$.}

First, Corollary~\ref{c!kil} implies that the $V_{LM}$ projection
sequence eliminates $y_a$ before $x_{i-1}$. Indeed, $x_{i-1}$ is joined
to $y_j$ in a cross or pitchfork involving at most $y_j$ and $y_{j-1}$,
and these are parties to which no $y_a$ with $a\le j-2$ is invited.
Therefore Lemma~\ref{l!cmp} gives
\begin{equation}
\label{eq!piAB}
\pi_{AB}\Bigl( \frac{y_a}{x_{i-1}} \Bigr) \ge 0.
\end{equation}

As before, comparing the alleged monomial $x_i^\xi y_j^\eta m\simT y_jy_a \in R(\nu)$
with the first of
\eqref{eq!mod2} gives
\begin{equation}
\frac{y_a}{x_{i-1}} \simT x_i^{\xi-\al_\nu}\frac{m}{A_\nu}. \label{eq!ya}
\end{equation}

The proof divides into two cases.

\case1 The projection sequence for $V_{AB}$ eliminates $x_i$ before $y_a$.

Lemma~\ref{l!cmp} says that $\xi-\al_\nu\ge1$ is impossible in
\eqref{eq!ya} (the lefthand side would have $\pi_{LM}$ strictly smaller
than the right). Thus
\begin{equation}
x_i^{\al_\nu-\xi}\frac{y_a}{x_{i-1}} \simT \frac{m}{A_\nu}
\quad\hbox{with} \quad \al_\nu\ge\xi, \label{eq!yaw}
\end{equation}
and \eqref{eq!piAB} implies that $A_\nu$ divides $m$.

\case2 The projection sequence for $V_{AB}$ eliminates $y_a$ before
$x_i$. This means that $y_a,y_j$ are both contained in the interval
$Y_{i-1}$ of Lemma~\ref{l!bar}, and that $y_a$ is not at the bottom:
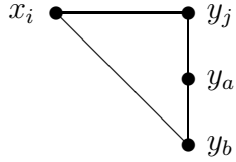
\[
\begin{picture}(76,55)(0,0)
\put(-12,48){$x_i$}
\put(6,50){\circle*{5}}
\put(6,50){\line(1,0){50}}
\put(56,50){\circle*{5}}
\put(63,48){$y_j$}
\put(56,25){\circle*{5}}
\put(63,23){$y_a$}
\put(56,0){\circle*{5}}
\put(63,-2){$y_b$}
\put(6,50){\line(1,-1){50}}
\put(56,50){\line(0,-1){50}}
\end{picture}
\]
Suppose that $x_i$ is tagged with $d$ (or simply replace $d\bij e$ in
what follows), and write $Y_{i-1}=[b,b+d-2]$ for the interval of
Lemma~\ref{l!bar}. Our conclusion in this case is that $\xi-\al_\nu<
a-b+1$ and $\equiv a-b+1 \mod d$, and so $\xi\le\al_\nu$, and the
argument of Case~1 works as before.

The proof goes as follows:
\begin{enumerate}
\renewcommand{\labelenumi}{(\alph{enumi})}
\item For $y_a\in Y_{i-1}$
\begin{equation}
\pi_{LM}(y_a)=(a-b)\pi_{LM}(x_i)+\pi_{LM}(y_b)<(a-b+1)\pi_{LM}(x_i).
\label{eq!piLM}
\end{equation}

\item On the other hand, taking $\pi_{LM}$ in \eqref{eq!ya} gives
\begin{equation}
(\xi-\al_\nu)\pi_{LM}(x_i) \le \pi_{LM}(y_a)-\pi_{LM}(x_{i-1})<\pi_{LM}(y_a)
\end{equation}
Therefore $\xi-\al_\nu<a-b+1\le d-2$.

\item Moreover modulo $\M'$, we have
\begin{equation}
y_a \equiv \frac{x_i^{a-b+1}}{x_{i-1}} \in Q.
\end{equation}

\item Therefore in \eqref{eq!ya}, $\xi-\al_\nu\equiv a-b+1\mod d$.
\end{enumerate}

\begin{pf} (a) follows from the tag equations for the toric variety
$V_{AB}$ at the successive $y_\al$: as $y_{\al+1}$ is eliminated it has
tag~1 and tag equation
\begin{equation}
x_iy_\al=y_{\al+1}A^{u_\al}B^{v_\al}. \label{eq!yal}
\end{equation}
Applying $\pi_{LM}$ gives the equality in \eqref{eq!piLM}, and the
inequality comes from Lemma~\ref{l!cmp}. (b) explains itself.

When we reach the bottom of this interval, we eliminate $x_i$, with
tag~1 and tag equation
\begin{equation}
x_{i-1}y_b=x_iA^{u_b}B^{v_b}. \label{eq!yb}
\end{equation}
Viewing this equation modulo $\M'$ gives $y_b\equiv x_i/x_{i-1}$, and
together with \eqref{eq!yal} this gives the value of $y_a$ in $Q$
as
\begin{equation}
y_a \equiv y_bx_i^{(a-b)} \equiv \frac {x_i^{a-b+1}}{x_{i-1}},
\label{eq!ysuba}
\end{equation}
which proves (c).

In the coordinates of the Padded Cell $Q$, we know that $x_{i-1}$ is
$(0,\pm\eeth)$ and $x_i$ is $(\pm \dth, 0)$. The alleged monomial tells
us that $y_a\equiv x_i^{\xi-\al}/x_{i-1}$ modulo $\M'$, and (d)
follows.
\end{pf}

\step4 \addcontentsline{toc}{subsubsection}{Step 4. Monomials $x_ix_a$
with $a=0,\dots,i-2$} {\em Claim~\ref{cl!div} holds for all monomials
$x_ix_a$ with $a=0,\dots,i-2$.}

The prevailing assumption that $d,e\ge2$ and $de>4$ is necessary here;
when $d=e=2$, Claim~\ref{cl!div} fails on equations of this type.

\begin{pf} We compare an alleged monomial $x_ix_a\simT y_j^\eta m$ with
the second expression of \eqref{eq!mod2} as usual; move the $y_j$ term across, this
time regardless of sign, obtaining
\begin{equation}
y_j^{\be_\nu-\eta}\frac{x_a}{y_{j-1}}\simT\frac{m}{B_\nu}.
\end{equation}
Our conclusion in this case is that $\be_\nu-\eta-1>-2$ and $d$ divides
$\be_\nu-\eta-1$; this implies that
$\pi_{AB}(m/B_\nu)\ge\pi_{AB}(x_ay_j/y_{j-1})\ge0$, so that $B_\nu$
divides $m$ as required.

The proof breaks up into cases as follows; we suppose that the pristine tag on $x_i$ is $d$:

\paragraph{The case $d>2$, $e>2$}
The alleged monomial is $x_ix_a \simT y_j^\eta m$, and we divide by
the second of \eqref{eq!mod2} to give
\begin{equation}\label{eq!star2}
y_j^{\be_\nu - \eta} \frac{x_a}{y_{j-1} }\simT \frac{m}{B_\nu}.
\end{equation}
This equation is the key at the end of the argument, but first we
rewrite it trivially as
\begin{equation}\label{eq!star3}
y_j^{\be_\nu - \eta-1} x_a\frac{y_j}{y_{j-1} }\simT \frac{m}{B_\nu}.
\end{equation}
Since $e>2$, then as long as $i>2$ the tag equation in $V_{AB}$ at $y_j$
when $y_j$ is about to be eliminated is
\[
x_{i-1}y_{j-1} \simT y_j A^\bullet B^\bullet
\]
for nonnegative powers of $A, B$ that will not concern us, which we
rewrite as
\begin{equation}\label{eq!yjj-1}
\frac{y_j}{y_{j-1} }\simT \frac{x_{i-1}}{A^\bullet B^\bullet}.
\end{equation}
Now \eqref{eq!star3} and \eqref{eq!yjj-1} together give
\begin{equation}\label{eq!star4}
y_j^{\be_\nu - \eta-1} x_a\frac{x_{i-1}}{A^\bullet B^\bullet}\simT \frac{m}{B_\nu}.
\end{equation}
Since $\pi_{LM}(m/B_\nu)\ge 0$, we have
\[
(\be_\nu-\eta-1)\pi_{LM}(y_j)\ge -\pi_{LM}(x_a) - \pi_{LM}(x_{i-1}),
\]
and since $y_j$ is eliminated in the projection sequence of $V_{AB}$ before
$x_a$ and $x_{i-1}$, we have $\be_\nu - \eta - 1 > -2$, or in other words that
\[
\be_\nu - \eta \ge 0.
\]
(The case $i=2$ is simpler: it must have $a=0$, so $\pi_{LM}(x_ix_a)=(0,*)$
by Proposition~\ref{p!xys} and thus $\eta=0$; in particular $\be_\nu-\eta\ge0$.)

Since $d>2$, the variable $x_a$ is eliminated before $y_{j-1}$ in the projection
sequence for $V_{LM}$ if and only if $a < i-2$.
When $a<i-2$, this projection implies that $\pi_{AB}(x_a) > \pi_{AB}(y_{j-1})$, so
that $\be_\nu - \eta \ge 0$ already gives
\[
0 \le \pi_{AB}(\hbox{LHS\eqref{eq!star2}}) \le \pi_{AB}(m) - \pi_{AB}(B_\nu).
\]
In other words, $B_\nu$ divides $m$, and we are done.
On the other hand, if $a=i-2$, then we can use the tag equation in $V_{LM}$
at the point that $y_{j-1}$ is eliminated (again using that $d>2$): namely
\begin{equation}\label{eq!tagyj-1}
x_{i-2}y_j \simT y_{j-1}L^\bullet M^\bullet
\end{equation}
(for nonnegative powers of $L, M$ that will not concern us). Writing this as
\[
\frac{x_a}{y_{j-1}} \simT \frac{L^\bullet M^\bullet}{y_j}
\]
and substituting into \eqref{eq!star2} gives
\[
y_j^{\be_\nu-\eta-1}L^\bullet M^\bullet \simT \frac{m}{B_\nu}.
\]
Since $y_j$ lies in a corner of the Padded Cell, this implies that
$d$ divides $\be_\nu-\eta-1$ (which is known from above to be $\ge -1$),
so $\be_\nu-\eta\ge 1$. With this, the tag equation \eqref{eq!tagyj-1}
implies again that $\pi_{AB}$ of the left-hand side of \eqref{eq!star2} is $\ge 0$,
and we conclude as before.

\paragraph{The case $d>2$, $e=2$}
The argument proceeds almost identically, except that the tag equation in $V_{AB}$
when $y_j$ is eliminated is now
\[
x_{i-2}y_{j-1} \simT y_j A^\bullet B^\bullet
\]
which we rewrite, to replace \eqref{eq!yjj-1} in the argument, as
\[
\frac{y_j}{y_{j-1} }\simT \frac{x_{i-2}}{A^\bullet B^\bullet}.
\]
The only effect is to replace occurrences $x_{i-1}$ by $x_{i-2}$, and the
conclusion $\be_\nu-\eta\ge0$ still holds.

\paragraph{The case $d=2$, $e>2$}
This case differs from the others by having a cross $x_i,x_{i-1},y_j,y_{j-1}$
at this projection bar, rather than the usual pitchfork.
Nevertheless, the proof follows without change to show that $\be_\nu-\eta\ge 0$.

But now it is easier: the cross (rather than pitchfork) implies that
$x_a$ is eliminated before $y_{j-1}$ for any $a\le i-2$, and the
proof follows as before.
\end{pf}

\section{Final remarks} \label{s!fin}

%%% \subsection{Relation with Mori flips of type A}

This paper arose out of a study of Mori's remarkable ``continued
division'' Euclidean algorithm \cite{Mo} in the divisor class group of
an extremal 3-fold neighbourhood $C\subset X$ (see also \cite{R}). As
Mori has explained to us over a couple of decades, the main result of
\cite{Mo} corresponds to a 2-step recurrent continued fraction
$[d,e,d,\dots]$ as in our Classification Theorem~\ref{th!d-e}. To
paraphrase his argument: an extremal neighbourhood of type~A has an
exceptional curve $C\iso\PP^1(r_1,r_2)_{\Span{x_k,y_l}}$ that is cut
transversally by two divisors $\divi x_k$, $\divi y_l$ through the
terminal points of type~A. Mori's algorithm replaces these two variables
successively by $x_k$, $x_{k-1}=Ax_k^d/y_l$, then
$x_{k-1},x_{k-2}=x_{k-1}^e/x_k$, continuing down the $d,e,d,\dots$ side
of a long rectangle until it reaches $x_0,y_0$, that are detected by a
sign reversal in their degrees. At this point, the two divisors $\divi
x_0$, $\divi y_0$ in $X$ intersect set-theoretically only in the curve
$C\subset X$. It follows that they define a pencil $X\broken\PP^1$, and,
in the flipping case, the flip $C^+\subset X^+$ as its normalised graph.
Our take on this is that the canonical cover of a Mori flip of Type~A
arises as a regular pullback from a diptych variety. We return to this
in \cite{BR4}.

%%% where we explain how a Mori flip of Type~A determines a diptych
%%% variety, as an application of the Main Theorem~\ref{th!main} here,
%%% and conversely we describe semistable Mori flips of each
%%% combinatorial type that arises.

%%% \subsection{Discussion of results}

We comment here on the equations of $V_{ABLM}$, since the proof by
unprojection in Section~\ref{s!pf} deduces them by unprojection, and
does not write them all out explicitly. The equations lift the toric
equations of $T$ or of $V_{AB}$, that are of the form $v_iv_j=\cdots$
for any pair of nonadjacent monomials on the boundary of $\si_{AB}$;
there are thus $\binom{k+l+1}2-1$ of them. Our favourites among them are
the Pfaffian equations coming from magic pentagrams, which are
trinomials, comparable to the binomial equations of toric geometry. These
determine everything, and, as in the extended example, the full set of
equations can be obtained if required. In practice, this amounts to
taking a colon ideal against powers of the top monomial $x_ky_lAB$, or
of the bottom monomial $x_0y_0LM$.

In toric geometry, we assume out of habit that Jung--Hirzebruch
continued fractions $[a_1,\dots,a_k]$ have entries $a_i\ge2$. In fact,
if a tag equation $x_{i-1}x_{i+1}=x_i^{a_i}$ has $a_i=1$ then $x_i$ is a
redundant generator. Our Classification Theorem~\ref{th!d-e}, already
implicit in Mori \cite{Mo}, Lemma~3.3, has output including $d$ or $e=1$
as regular cases. It turns out that $de=4$ should be treated separately
(see \cite{BR2}). We used the {\em main case} assumption $d,e\ge2$ and
$de>4$ in an essential way at several point in the proof of Main
Theorem~\ref{th!main} in Sections~\ref{s!pp}--\ref{s!pf}. In particular,
the order of unprojecting variables from top and bottom was determined
by the Scissors of Figure~\ref{f!sciss}. In the cases $d$ or $e=1$ with
$de>4$, we prove Theorem~\ref{th!main} in \cite{BR3} by an argument that
keeps the variable $x_i$ tagged with 1. Allowing $d=1$ (say) changes the
shape of Scissors, and the nature of the Pretty Polytope and the Padded
Cell of Section~\ref{s!pp}. In fact, the $x_i$ marked with 1 are
redundant generators, and should be projected out before their
neighbours marked with $e>4$. This obliges us to start the proof again
from scratch adopting a new order of projection; the proof then goes
through in parallel with the main case.

%%% \subsection{Relation with GHK}

We treat the extension $T\subset V_{AB}$ in closed form, rather than
via the infinitesimal deformations of Altmann \cite{alt}. Paul Hacking
observes that our diptych varieties $V_{ABLM}$ can be seen as a variant
of the construction of Gross, Hacking and Keel \cite{GHK} in a special
case. Their starting point is the {\em vertex of degree} $n$, the tent
in $\Aff^n$ that is the $n$-cycle formed by the coordinate planes
$\Aff^2_{\Span{x_i,x_{i+1}}}$. They construct a formal scheme that
deforms this tent; their basic idea is to smooth $T$, replacing the
local equation $x_{i-1}x_{i+1}=0$ of $T$ in a neighbourhood of each
punctured $x_i$-axis by the tag equation
\begin{equation}
x_{i-1}x_{i+1}=A_ix_i^{a_i},
\label{eq!GHKtag}
\end{equation}
where the tags $a_i=-D_i^2$ arise from a cycle of rational curves
$D=D_1+\cdots+D_n$ on a mirror log Calabi--Yau surface $Y,D$ and the
deformation parameters $A_i$ play a similar role to our annotations
$A,B,L,M$. The main difference between the two constructions is
summarised by the slogan ``perturbative versus nonperturbative''.
Whereas we work in closed form with varieties and birational
unprojections, \cite{GHK} proceed by successive infinitesimal steps:
their affine pieces \eqref{eq!GHKtag} are glued by formal power series
expansions in affine linear transformations, specified by Gromov--Witten
theory of $Y,D$ and the Kontsevich--Soibelman and Gross--Siebert
scattering diagram.

 \label{end!bib}

\bigskip

\noindent Gavin Brown, \\
Department of Mathematical Sciences, Loughborough University, \\
LE11 3TU, England

\noindent {\it e-mail}: G.D.Brown@lboro.ac.uk

\bigskip

\noindent Miles Reid, \\
Mathematics Institute, University of Warwick, \\
Coventry CV4 7AL, England

\noindent {\it e-mail}: Miles.Reid@warwick.ac.uk

\end{document}